\documentclass[pdflatex,sn-mathphys-num]{sn-jnl}

\usepackage{color}
\usepackage{geometry}
\usepackage{subfiles}
\usepackage{graphicx}
\usepackage{xcolor}
\usepackage[section]{placeins}
\usepackage{hyperref}
\usepackage{stmaryrd}
\usepackage{amsmath,amssymb,amsfonts,amsthm,textcomp}
\usepackage{mathtools}
\usepackage{tensor}
\usepackage{multicol}
\usepackage{subfig}
\usepackage{enumitem}
\usepackage{setspace}
\usepackage{nicefrac}
\usepackage{algpseudocode}
\usepackage{algorithm}
\usepackage{booktabs}
\usepackage{tabularx}
\usepackage{bbm}
\usepackage{cancel}
\usepackage{soul}

\usepackage{multirow}%
\usepackage{mathrsfs}%
\usepackage[title]{appendix}%
\usepackage{xcolor}%
\usepackage{textcomp}%
\usepackage{manyfoot}%
\usepackage{listings}%

\newfont{\tenbfsl}{cmbxti9 scaled 1200}
\newfont{\tenbbb}{msbm10}
\newfont{\svnbbb}{msbm8}

\newcommand{\inner}[2]{\left\langle#1,#2\right\rangle}
\newcommand{\norm}[1]{\left\|#1\right\|}

\newcommand{\bs}[1]{\boldsymbol{#1}}
\newcommand{\cl}[1]{\mathcal{#1}}

\newcommand{\bb}[1]{\mathbb{#1}}

\newtheorem{thm}{Theorem}
\newtheorem{cor}{Corollary}
\newtheorem{lem}{Lemma}
\newtheorem{prop}{Proposition}
\theoremstyle{remark}
\newtheorem{rmk}{Remark}
\newtheorem{set}{Definition}
\newtheorem{ass}{Assumption}

\newcommand{\rset}{\bb{R}}

\newcommand{\Exp}[1]{\bb{E}\left[#1\right]}
\newcommand{\Expc}[2]{\bb{E}\left[\left.#1\,\right|\,#2\right]}

\newcommand{\Probc}[2]{\bb{P}\left[\left.#1\,\right|\,#2\right]}
\newcommand{\opt}{\bs{\xi}^*}
\newcommand{\grad}{\nabla_{\bs{\xi}}}
\newcommand{\bigO}[1]{\cl{O}\left( #1 \right)}
\newcommand{\xiset}[1]{\left\{ \bs{\xi}_{\ell'} \right\}_{\ell' \in \cl{L}_{#1}}}

\newcommand{\xik}{\bs{\xi}_k}
\newcommand{\xio}{\bs{\xi}_0}
\newcommand{\xikp}{\bs{\xi}_{k+1}}

\newcommand{\xil}{\bs{\xi}_\ell}
\newcommand{\xilpr}{\bs{\xi}_{p_k(\ell)}}
\newcommand{\xilp}{\bs{\xi}_{p_k(\ell)}}
\newcommand{\Fk}{\nabla \cl{F}_k}

\newcommand{\tht}{\bs{\theta}}

\newcommand{\bias}{\bs{b}_k}
\newcommand{\muhat}[1]{\hat{\bs{\mu}}_{#1}}

\newcommand{\cost}[2][k]{(1+\mathbbm{1}_{\mskip-2mu\underline{\cl{L}_{#1}}}(#2))}
\newcommand{\V}[2][k]{V_{#2, #1}}
\newcommand{\M}[2][k]{M_{#2, #1}}

\newcommand{\gnorm}[1]{\norm{\grad F\left(\bs{\xi}_{#1}\right)}}

\newcommand{\Add}[1][ ]{\texttt{Add}#1}
\newcommand{\Drop}[1][ ]{\texttt{Drop}#1}
\newcommand{\Restart}[1][ ]{\texttt{Restart}#1}
\newcommand{\Clip}[1][ ]{\texttt{Clip}#1}

\newcommand{\MICE}[1][ ]{\texttt{MICE}#1}
\newcommand{\SGD}[1][ ]{\texttt{SGD#1}}

\newcommand{\Adam}[1][ ]{\texttt{Adam}#1}

\begin{document}

\title[Multi-Iteration Stochastic Optimizers]{Multi-Iteration Stochastic Optimizers}

\author*[1]{\fnm{André} \sur{Carlon}}\email{carlon@uq.rwth-aachen.de}
\author[2]{\fnm{Luis} \sur{Espath}}
\author[3]{\fnm{Rafael} \sur{Holdorf}}
\author[1,4,5]{\fnm{Raúl} \sur{Tempone}}

\affil[1]{\orgdiv{Department of Mathematics}, \orgname{RWTH Aachen University}, \orgaddress{\street{Pontdriesch 14-16}, \city{Aachen}, \postcode{52062}, \country{Germany}}}

\affil[2]{\orgdiv{School of Mathematical Sciences}, \orgname{University of Nottingham}, \orgaddress{\street{NG7 2RD}, \city{Nottingham}, \country{United Kingdom}}}

\affil[3]{\orgdiv{School of Engineering}, \orgname{Federal University of Santa Catarina}, \orgaddress{\street{Rua Jo\~{a}o Pio Duarte da Silva}, \city{Florian\'{o}polis}, \postcode{88040-970}, \state{SC}, \country{Brazil}}}

\affil[4]{\orgdiv{Computer, Electrical and Mathematical Sciences \& Engineering Division}, \orgname{King Abdullah University of Science \& Technology}, \orgaddress{\city{Thuwal}, \postcode{23955-6900},\country{Saudi Arabia}}}

\affil[5]{\orgdiv{Alexander von Humboldt Professor in Mathematics for Uncertainty Quantification}, \orgname{RWTH Aachen University}, \orgaddress{\country{Germany}}}

\abstract{
  We introduce Multi-Iteration Stochastic Optimizers, a novel class of first-order stochastic methods that control the relative $L^2$ error using successive control variates along the iteration path. By exploiting correlations between iterates, these control variates reduce the estimator's variance, making an accurate mean gradient estimation computationally affordable. Our approach centers on the Multi-Iteration stochastiC Estimator (MICE), which can be seamlessly coupled with any first-order stochastic optimizer due to its non-intrusive design. The algorithm adaptively selects which iterates to include in its index set. We provide both an error analysis of MICE and a convergence analysis for Multi-Iteration Stochastic Optimizers across various problem classes, including some non-convex cases. In the smooth, strongly convex setting, we demonstrate that to approximate a minimizer within a tolerance $tol$, SGD-MICE requires, on average, $O(tol^{-1})$ stochastic gradient evaluations, compared to $O(tol^{-1}\log(tol^{-1}))$ for SGD with adaptive batch sizes. In numerical experiments, SGD-MICE achieved the desired tolerance with fewer than 3\% of the gradient evaluations required by adaptive batch SGD. Additionally, MICE offers a straightforward stopping criterion based on the gradient norm, validated through consistency tests. To assess its efficiency, we present examples using both SGD-MICE and Adam-MICE, including a stochastic adaptation of the Rosenbrock function and logistic regression on various datasets. Compared to SGD, SAG, SAGA, SVRG, and SARAH, our approach consistently reduces the gradient sampling cost without the need for extensive parameter tuning.

\noindent\textbf{AMS subject classifications:}
$\cdot$
62L20 
$\cdot$
65K05 
$\cdot$
90C15 
$\cdot$
65C05 
$\cdot$
\\
}


\maketitle

\section{Introduction}
\label{sec:intro}

We focus on the stochastic optimization problem of minimizing the objective function $\bb{E}\left[f(\bs{\xi},\bs{\theta})|\bs{\xi}\right],$ where $f$ is a given real-valued function, $\bs{\xi}$ is the design variable vector, $\bs{\theta}$ is a random vector, and $\bb{E}\left[\cdot|\bs{\xi}\right]$ is the expectation conditioned on $\bs{\xi}$.
Stochastic optimization problems \cite{Marti,Pardalos,Birge2011} are relevant to different fields, such as machine learning \cite{Lan2020}, Stochastic Optimal Control \cite{wallace2005applications,fleming2012deterministic}, Computational Finance \cite{Ziemba1975,Conejo,bayer2018pricing}, Economics \cite{Chang}, Insurance \cite{Azcue}, Communication Networks \cite{ding2012stochastic}, Queues and Supply Chains \cite{Yao}, and  Bayesian Optimal Design of Experiments \cite{ryan2016review,carlon2020nesterov}, among many others.

In the same spirit and inspired by the work by Heinrich \cite{heinrich2000multilevel} and Giles \cite{giles:MLMC} on \emph{Multilevel Monte Carlo} methods, we propose the \emph{Multi-Iteration stochastiC Estimator}---\texttt{MICE}---to obtain a computationally efficient approximation of the mean gradient at iteration $k$, $\nabla_{\bs{\xi}} \bb{E}\left[f(\bs{\xi},\bs{\theta})| \bs{\xi} = \bs{\xi}_k\right]$, which may be coupled with any first order stochastic optimizer in a non-intrusive fashion.
Combining \texttt{MICE} with any stochastic optimizer furnishes \emph{Multi-Iteration Stochastic Optimizers}, a novel class of efficient and robust stochastic optimizers.
In this class of stochastic optimizers, the mean gradient estimator's relative variance is controlled using successive control variates based on previous iterations' available information.
This procedure results in a more accurate yet cost-effective estimation of the mean gradient.
In approximating the mean gradient, \texttt{MICE} constructs an index set of iterations and performs control variates for every pair of nested elements of this index set.
As the stochastic optimization evolves, we increase the number of samples along the index set while keeping the previously sampled gradients, i.e., we use gradient information from previous iterations to reduce the variance in the current gradient estimate, which is a crucial feature to make \texttt{MICE} competitive.
We design \texttt{MICE} to achieve a given relative error for the mean gradient with minimum additional gradient sampling cost.
Indeed, in the \texttt{MICE} index set constructed along the stochastic optimization path $\{\bs{\xi}_\ell\}_{\ell=0}^k$, our generic optimizer optimally decides whether to drop a particular iteration $\ell$ out of the index set or restart it to reduce the total optimization work.
Moreover, it can decide if it is advantageous, from the computational work perspective, to clip the index set at some point $\ell$, discarding iterations before $\ell$.
Since we control the gradients' error using an estimate of the gradient norm, we propose a resampling technique to get a gradient norm estimate, reducing the effect of sampling error and resulting in robust optimizers.
We note in passing that \texttt{MICE} can be adjusted to the case of finite populations; see \eqref{eq:finite_sum_problem}, for optimization problems arising in supervised machine learning.

Generally speaking, in first-order stochastic optimization algorithms that produce convergent iterates, the mean gradient converges to zero as the number of iterations, $k$, goes to infinity, that is $\norm{\bb{E}\left[\nabla_{\bs{\xi}}f(\bs{\xi}_k,\bs{\theta})\right]}\to0$; however, the gradient covariance,
$\bb{C}[\nabla_{\bs{\xi}}f(\bs{\xi}_k,\bs{\theta}),\nabla_{\bs{\xi}}f(\bs{\xi}_k,\bs{\theta})]$,
does not converge to zero.
Thus, to ensure convergence of the iterates $\bs{\xi}_k$, in the literature it is customary to use decreasing step-size (learning rate) schedules, reducing the effect of the statistical error in the gradient onto the iterates $\bs{\xi}_k$.
However, this approach also results in sublinear convergence rates \cite{Rup88}. Another approach to deal with the gradient's statistical error is to increase the sample sizes (batch sizes) while keeping the step-size fixed, thus avoiding worsening the convergence.
Byrd et al. \cite{Byr12} propose to adaptively increase the sample sizes to guarantee that the trace of the covariance matrix of the mean gradient is proportional to its norm.
This approach forces the statistical error to decrease as fast as the gradient norm.
Balles et al. \cite{balles2016coupling} use a similar approach; however, instead of setting a parameter to control the statistical error, they set a step-size and find the parameter that guarantees the desired convergence.
Bollapragada et al. \cite{bollapragada2018adaptive} propose yet another approach to control the variance of gradient estimates in stochastic optimization, which they call the inner product test.
Their approach ensures that descent directions are generated sufficiently often.
  Relatedly, for certain over-parameterized (interpolating) models it is possible to establish convergence of stochastic methods with constant step-size under the \emph{strong growth condition} \cite{vaswani2019fast}, which bounds the second moment of the (raw) stochastic gradient by a multiple of $\|\nabla F(\bs{\xi})\|^2$.
  This assumption is conceptually different from our approach: rather than postulating a property of the intrinsic gradient noise, we adaptively allocate samples so that the \texttt{MICE} mean-gradient estimator satisfies a relative accuracy requirement with user-chosen tolerance $\epsilon$ in $L^2$.

Instead of increasing the sample size, some methods rely on using control variates with respect to previously sampled gradients to reduce the variance in current iterations and thus be able to keep a fixed step-size.
Pioneering ideas of control variates in stochastic optimization, by Johnson \& Zhang \cite{Joh13}, profit on an accurate mean gradient estimation at the initial guess $\bs{\xi}_0$, $\nabla_{\bs{\xi}}\bb{E}\left[f(\bs{\xi},\bs{\theta})|\bs{\xi} = \bs{\xi}_{0}\right]$, to update and compute, via single control variates, an inexpensive and accurate version of the mean gradient at the iteration $k$, $\nabla_{\bs{\xi}}\bb{E}\left[f(\bs{\xi},\bs{\theta})|\bs{\xi} = \bs{\xi}_k\right]$.
Instead of doing control variates with respect to one starting full-gradient, \texttt{SARAH}, by Nguyen et al. \cite{nguyen2017sarah}, computes an estimate of the gradient at the current iteration by using control variates with respect to the last iteration.
An `inexact' version of \texttt{SARAH} is presented in \cite{nguyen2018inexact}, where \texttt{SARAH} is generalized to the minimization of expectations.
In the spirit of successive control variates, \texttt{SPIDER} by Fang et al. \cite{fang2018spider} uses control variates between subsequent iterations; however, it employs the \emph{normalized} gradient descent instead of plain gradient descent.
In a different approach, \texttt{SAGA}, by Defazio et al. \cite{defazio2014saga}, keeps in the memory the last gradient $\nabla_{\bs{\xi}} f$ observed for each data point and computes $\nabla_{\bs{\xi}}\bb{E}\left[f(\bs{\xi}_k,\bs{\theta})|\bs{\xi}_k\right]$ using control variates with respect to the average of this memory.
Lastly, many algorithms try to `adapt' the initial batch size of the index set of batches using predefined rules, such as exponential or polynomial growth, as presented by Friedlander \& Schmidt \cite{friedlander2012hybrid}, or based on statistical bounds as discussed by De et al. \cite{de2016big} and Ji et al. \cite{ji2019faster}, to mention a few.

Although our proposed \emph{Multi-Iteration Stochastic Optimizers} share similarities with \texttt{SVRG}  \cite{Joh13}, \texttt{SARAH}, and \texttt{SPIDER}, our stochastic optimizers distinctly control the relative variance in gradient estimates.
We achieve this control by sampling the entire index set of iterations, optimally distributing the samples to minimize the gradient sampling cost.
While the previously mentioned methods are devised for finite sum minimization, \texttt{MICE} can tackle both finite sum and expectation minimization.
Moreover, we provide additional flexibility by including dropping, restart, and clipping operations in the \texttt{MICE} index set updates.

For strongly-convex and $L$-smooth objective functions, Polyak, in his book \cite[Theorem 5, pg 102]{polyak1987introduction}, shows a convergence rate in the presence of random relative noise. 
The theorem states a linear (geometric) convergence $cq^k$ in terms of the number of iterations. 
However, the dependency on the relative noise level, $\epsilon$, of the constants $c$ and $q$ is not made explicit. 
This work presents the explicit form of these constants and their dependency on $\epsilon$. Using this, we can estimate the total average computational work in stochastic gradient evaluations and optimize it with respect to the controllable relative noise $\epsilon$. 
Finally, we conclude that to generate an iterate $\xik$ such that  $\norm{\grad F(\xik)}^2 < tol$, \texttt{SGD-MICE} requires, on average, $\cl{O}(tol^{-1})$ stochastic gradient evaluations, while \texttt{SGD} with adaptive batch sizes requires the larger $\cl{O}(tol^{-1}\log(tol^{-1}))$, correspondingly, as we establish formally in \S\ref{sec:sampling_cost} (Corollaries~\ref{cor:cost_sgd_mice} and~\ref{cor:cost_sgd_a}).
While the reuse of previous data causes the MICE estimator to be conditionally biased, we present an analysis for the conditional bias and characterize the $L^2$ error, including bias and statistical error, which is controlled to achieve convergence of \SGD[-MICE].

Since \texttt{MICE} is non-intrusive and designed for both continuous and discrete random variables, it can be coupled with most available optimizers with ease.
For instance, we couple \texttt{MICE} with \texttt{SGD} \cite{robbins1951stochastic} and \texttt{Adam} \cite{kingma2014adam}, showing the robustness of our approach.
The \texttt{Adam} algorithm by Kingma \& Ba \cite{kingma2014adam} does not exploit control variates techniques for variance reduction.
Instead, it reduces the gradient estimator's variance based on iterate history by adaptive estimates of lower-order moments, behaving similarly to a filter.
Thus, the coupling \texttt{Adam-MICE} profits from the information available in the optimizer path in more than one way.

More generally, \texttt{MICE} can also be coupled with momentum schemes (e.g., heavy-ball or Nesterov-type momentum) by replacing the stochastic gradient oracle with the \texttt{MICE} estimator.

Finally, the reader is referred to the books by Spall \cite{spall2005introduction} and Shapiro, Dentcheva, and Ruszczy{\'n}ski \cite{shapiro2014lectures} for comprehensive overviews on stochastic optimization.

To assess \texttt{MICE}'s applicability, we numerically minimize expectations of continuous and discrete random variables using analytical functions and logistic regression models.
Also, we compare \texttt{SGD-MICE} with \texttt{SVRG}, \texttt{SARAH}, \texttt{SAG}, and \texttt{SAGA} in training the logistic regression model with datasets with different sizes and numbers of features.

\subsection{Optimization of expectations and stochastic optimizers}
\label{sec:so_definitions}

To state the stochastic optimization problem, let $\bs{\xi}$ be the design variable in dimension $d_{\bs{\xi}}$ and $\bs{\theta}$ a vector-valued random variable in dimension $d_{\bs{\theta}}$, whose probability distribution $\pi$ may depend on  $\bs{\xi}.$ Throughout this work we assume that we can produce as many independent identically distributed samples from $\pi$ as needed.
Here, $\bb{E}\left[\cdot|\bs{\xi}\right]$ and $\bb{V}\left[\cdot|\bs{\xi}\right]$ are respectively the expectation and variance operators conditioned on $\bs{\xi}$. 
Aiming at optimizing expectations on $\bs{\xi}$, we state our problem as follows. Find $\bs{\xi}^*$ such that
\begin{equation}\label{eq:underlying.problem}
\bs{\xi}^* = \underset{\bs{\xi} \in \rset^{d_{\bs{\xi}}}}{\arg\min} \,  \bb{E}\left[f(\bs{\xi},\bs{\theta})\right],
\end{equation}
where
$f \colon \rset^{d_{\bs{\xi}}} \times \rset^{d_{\bs{\theta}}} \to \bb{R}$. Through what follows, let the objective function in our problem be denoted by $F(\bs{\xi}')\coloneqq\bb{E}\left[f(\bs{\xi},\bs{\theta})|\bs{\xi}=\bs{\xi}'\right]$.
In general, function $F$ might not have a unique minimizer, in which case we define $\Xi^*$ as the set of all $\bs{\xi}^*$ satisfying \eqref{eq:underlying.problem}.
The case of minimizing a finite sum of functions is of special interest given its importance for training machine learning models in empirical risk minimization tasks,
\begin{equation} \label{eq:finite_sum_problem}
        \bs{\xi}^* = \underset{\bs{\xi} \in \rset^{d_{\bs{\xi}}}}{\arg\min} \,  
        \frac{1}{N}\sum_{n=1}^N 
        f(\bs{\xi},\bs{\theta}_n),
\end{equation}
where $N$ is usually a large number.
Note that the finite sum case is a special case of the expectation minimization, i.e., let $\bs{\theta}$ be a random variable with probability mass function
\begin{equation}
    \bb{P}(\bs{\theta} = \bs{\theta}_n) = \frac{1}{N}.
\end{equation}
We recall that throughout this work we assume the sampling distribution $\pi$ of $\bs{\theta}$ is independent of $\bs{\xi}$ (see Assumption~\ref{as:theta_distribution} and Remark~\ref{rmk:theta_distribution}); in the finite-sum setting this corresponds to sampling from a fixed population.

In minimizing \eqref{eq:underlying.problem} with respect to the design variable $\bs{\xi} \in \rset^{d_{\bs{\xi}}}$, \texttt{SGD} is constructed with the following updating rule
\begin{equation}\label{eq:sgd}
\bs{\xi}_{k+1} = \bs{\xi}_{k} - \eta_k \bs{\upsilon}_k,
\end{equation}
where $\eta_k>0$ is the step-size at iteration $k$ and $\bs{\upsilon}_k$ is an estimator of the gradient of $F$ at $\bs{\xi}_k$.
For instance, an unbiased estimator $\bs{\upsilon}_k$ of the gradient of $F$ at $\bs{\xi}$ at the iteration $k$ may be constructed by means of a Monte Carlo estimator, namely
\begin{equation}\label{eq:grad.def.2}
\nabla_{\bs{\xi}} F(\bs{\xi}_k) =  \bb{E}\left[\nabla_{\bs{\xi}} f(\bs{\xi},\bs{\theta})|\bs{\xi} = \bs{\xi}_k\right] \approx \bs{\upsilon}_k :=\dfrac{1}{M} \sum_{\alpha\in\cl{I}} \nabla_{\bs{\xi}} f(\bs{\xi}_k,\bs{\theta}_{\alpha}),
\end{equation}
with $M$ independent and identically distributed (iid) random variables $\bs{\theta}_\alpha\sim\pi$ given $\bs{\xi}_k$,  $\alpha\in\cl{I}$, with $\cl{I}$ being an index set with cardinality $M\coloneqq|\cl{I}|$. Bear in mind that an estimator of the type \eqref{eq:grad.def.2} is, in fact, a random variable and its use in optimization algorithms gives rise to the so-called \emph{Stochastic Optimizers}. The challenge of computing the gradient of $F$ in an affordable and accurate manner motivated the design of several gradient estimators.

For the sake of brevity, the following  review on control variates techniques for stochastic optimization is not comprehensive. To motivate our approach, we recall the control variates proposed by Johnson \& Zhang \cite{Joh13} (and similarly, by Defazio et al.\ \cite{defazio2014saga}) for the optimization of a function defined by a finite sum of functions. The idea of control variates is to add and subtract the same quantity, that is, for any $\bs{\xi}_0$,
\begin{equation}\label{eq:grad.svrg}
\nabla_{\bs{\xi}} F (\bs{\xi}_k) = \bb{E}\left[\nabla_{\bs{\xi}} f(\bs{\xi},\bs{\theta}) - \nabla_{\bs{\xi}} f(\bs{\xi}_0,\bs{\theta})| \bs{\xi} = \bs{\xi}_k \right] + \bb{E}\left[ \nabla_{\bs{\xi}} f(\bs{\xi}_0,\bs{\theta})\right],
\end{equation}
rendering the following sample-based version
\begin{equation}\label{eq:grad.svrg.discrete}
\nabla_{\bs{\xi}} F (\bs{\xi}_k) \approx \dfrac{1}{M_k} \sum_{\alpha\in\cl{I}_k} (\nabla_{\bs{\xi}} f(\bs{\xi}_k,\bs{\theta}_{\alpha}) - \nabla_{\bs{\xi}} f(\bs{\xi}_0,\bs{\theta}_{\alpha})) + \dfrac{1}{M_0-M_k}\sum_{\alpha\in\cl{I}_0 \!\setminus \cl{I}_k} \nabla_{\bs{\xi}} f(\bs{\xi}_0,\bs{\theta}_{\alpha}),
\end{equation}
where $M_0 \gg M_k$ and $\bs{\theta}_{\alpha}$ are iid samples from the $\pi$ distribution, which does not depend on $\bs{\xi}$ in their setting.
In the original work by Johnson \& Zhang \cite{Joh13}, $M_0$ is the total population and $M_k=1$. Later, Nitanda \cite{Nit14} and Kone\v{c}n\`y et al.\ \cite{Kon16} also used the total populations $M_0$ at $\bs{\xi}_0$, but with $M_k=2,4,8\ldots$, to study the efficiency of the algorithm. Additionally, the work \cite{Joh13} restarts the algorithm after a pre-established number of iterations by setting $\bs{\xi}_0\leftarrow\bs{\xi}_k$. The efficiency of this algorithm relies on the correlation between the components of the gradients $\nabla_{\bs{\xi}} F (\bs{\xi}_0)$ and $\nabla_{\bs{\xi}} F (\bs{\xi}_k)$. If this correlation is high, the variance of the mean gradient estimator \eqref{eq:grad.svrg.discrete} is reduced.

\subsection{Paper outline}

The remainder of this work is as follows.
In \S\ref{sec:intro}, we describe the stochastic optimization problem, classical stochastic optimization methods and motivate variance reduction in this context.
In \S\ref{sec:mice}, we construct the \texttt{MICE} statistical estimator \S\ref{sec:mice.construction}; analyze its error \S\ref{sec:mice.error}; compute the optimal number of samples for the current index set \S\ref{sec:mice.optimal.setting}; present the operators used to build \MICE[]'s index set and derive a work-based criterion to choose one \S\ref{sec:opt_index_set}.
In \S\ref{sec:mice_conv_and_work}, we present a convergence analysis of $L^2$ error-controlled \texttt{SGD}, which includes \texttt{SGD-MICE}, showing these converge polynomially for general $L$-smooth problems, and exponentially if the objective function is gradient-dominated \S\ref{sec:convergence}.
In \S\ref{sec:sampling_cost} we present gradient sampling cost analyses for \texttt{SGD-MICE} and \texttt{SGD-A} (\texttt{SGD} with adaptive increase in the sample sizes) on expectation minimization~\S\ref{sec:sampling_cost_expectation} and finite sum minimization~\S\ref{sec:sampling_cost_finite}.
In \S\ref{sec:mice_algorithm}, practical matters related to implementation of the \texttt{MICE} estimator are discussed.
In \S\ref{sec:numerics}, to assess the efficiency of \emph{Multi-Iteration Stochastic Optimizers}, we present some numerical examples, ranging from analytical functions to the training of a logistic regression model over datasets of size up to $11 \times 10^6$.
In Appendix~\ref{sec:algorithms}, are presented detailed pseudocodes for the \emph{Multi-Iteration Stochastic Optimizers} used in this work.

\section{Multi-iteration stochastic optimizers}\label{sec:mice}

\subsection{Multi-iteration gradient estimator}\label{sec:mice.construction}

We now construct an efficient estimator of the mean gradient at the current iteration $k$,
$\nabla_{\bs{\xi}} F(\bs{\xi}_k) = \bb{E}\left[\nabla_{\bs{\xi}}f(\bs{\xi}_k,\bs{\theta})|\bs{\xi}_k\right]$, which we name
\emph{Multi-Iteration stochastiC Estimator}---\texttt{MICE}.
Profiting from available information already computed in previous iterations, \texttt{MICE} uses multiple control variates between pairs of, possibly non-consecutive, iterations along the optimization path to approximate the mean gradient at iteration $k$.
Bearing in mind that stochastic optimization algorithms, in a broad sense, create an $L^2$ convergent path where $\Exp{\norm{\bs{\xi}_k-\bs{\xi}_\ell}^2} \to 0$ as $\ell,k \to \infty$,
the gradients evaluated at $\bs{\xi}_{\ell}$ and $\bs{\xi}_k$ should become more and more correlated for $k,\ell \to \infty$. In this scenario, control variates with respect to previous iterations become more efficient, in the sense that one needs fewer and fewer new samples to accurately estimate the mean gradient.

  To motivate our \texttt{MICE} estimator, first consider the following special case:
  \begin{rmk}[Add-only special case]\label{rmk:add_only}
  A simple and special case of \MICE[] is the \emph{add-only} case:
  \begin{equation}
    \grad\cl{F}_k 
    =
    \frac{1}{M_{0,k}} \sum_{\alpha \in \cl{I}_{0,k}} \grad f(\xio, \bs{\theta}_\alpha)
    +
    \sum_{\ell=1}^k
      \frac{1}{M_{\ell,k}} \sum_{\alpha \in \cl{I}_{\ell,k}} 
      \left(
        \grad f(\xil, \bs{\theta}_\alpha)
        -
        \grad f(\bs\xi_{\ell-1}, \bs{\theta}_\alpha)
      \right),
  \end{equation}
  which is identical to \texttt{SARAH} but with \emph{cumulative sampling}.
  Note that the sample sizes $M_{\ell, k}$ change both with $\ell$ and $k$.
  \end{rmk}
  Now, we generalize our estimator with an adaptive selection of which iterates $\xil$ we want to consider in the \MICE[] estimator.
First let us establish some notation. Let $\cl{L}_k$ be an index set, such that, $\cl{L}_k\subset\{0,\ldots,k\}$, where $k$ is the current iteration and $k\in\cl{L}_k$.
This index set is $\cl{L}_0= \{0\}$ at the initial iteration, $k=0$, and for later iterations it contains the indices of the iterations \texttt{MICE} uses to reduce the computational work at the current iteration, $k>0$, via control variates.
For the special case in Remark~\ref{rmk:add_only} the index set is $\cl{L}_k = \{0, 1, \ldots, k\}$.
Next, for any $\min\{\cl{L}_k\}<\ell\in \cl{L}_k$, let $p_k(\ell)$ be the element \emph{previous} to $\ell$ in $\cl{L}_k$,

\begin{equation}\label{eq:previous.operator}
p_k(\ell)\coloneqq\max\{\ell^\prime\in\cl{L}_k\colon\ell^\prime<\ell\}.
\end{equation}
Then, the mean gradient at $\bs{\xi}_k$ conditioned on the sequence of random iterates, $\bs{\xi}$, indexed by the set $\cl{L}_k$ can be decomposed as
\begin{equation}\label{eq:mice.continuous}
\grad F(\xik) = \bb{E}\left[\nabla_{\bs{\xi}} f(\bs{\xi}_k, \bs{\theta})|\{\bs{\xi}_\ell\}_{\ell\in\cl{L}_k}\right]=\sum_{\ell\in\cl{L}_k}\bs{\mu}_{\ell,k},\qquad\bs{\mu}_{\ell,k}\coloneqq\bb{E}\left[\Delta_{\ell,k}|\bs{\xi}_\ell,\bs{\xi}_{p_k(\ell)}\right],
\end{equation}
with the gradient difference notation
\begin{equation}\label{eq:Delta.continuous}
\Delta_{\ell,k} \coloneqq
\left\{
\aligned
\nabla_{\bs{\xi}} f(\bs{\xi}_\ell, \bs{\theta}) - \nabla_{\bs{\xi}} f(\bs{\xi}_{p_k(\ell)}, \bs{\theta}),& \text{ if } \ell > \min\{\cl{L}_k\}, \\
\nabla_{\bs{\xi}} f(\bs{\xi}_\ell, \bs{\theta}),& \text{ if  } \ell = \min\{\cl{L}_k\}.
\endaligned
\right.
\end{equation}
Thus, the conditional mean $\bs{\mu}_{\ell,k}$ defined in \eqref{eq:mice.continuous} is simply
\begin{equation}
  \bs{\mu}_{\ell,k} =
  \left\{
  \begin{aligned}
  \grad F(\bs{\xi}_\ell) -
  \grad F(\bs{\xi}_{p_k(\ell)}) ,& \text{ if } \ell > \min\{\cl{L}_k\}, \\
  \grad F(\bs{\xi}_\ell)  ,& \text{ if } \ell =\min\{\cl{L}_k\}.
  \end{aligned}
  \right.
\end{equation}

For readability, we make the assumption that the distribution of $\bs{\theta}$ does not depend on $\bs{\xi}.$
Observe that this assumption is more general than it may seem, see the discussion on Remark~\ref{rmk:theta_distribution}.
\begin{ass}[Simplified probability distribution of $\bs{\theta}$]\label{as:theta_distribution}
The probability distribution of $\bs{\theta}$, $\pi$, does not depend on $\bs{\xi}$.
\end{ass}
Now we are ready to introduce the \texttt{MICE} gradient estimator.
\begin{set}[\texttt{MICE} gradient estimator]\label{df:general.mice.estimator}
Given an index set $\cl{L}_k$ such that $k\in\cl{L}_k\subset\{0,\ldots,k\}$ and positive integer numbers $\{M_{\ell,k}\}_{\ell\in\cl{L}_k}$, we define the \texttt{MICE} gradient estimator for $\grad F(\bs{\xi}_k)$ at iteration $k$ as
\begin{equation}\label{eq:mice}
\grad\cl{F}_k 
= \sum_{\ell\in\cl{L}_k} \hat{\bs{\mu}}_{\ell,k}
,
\qquad
\hat{\bs{\mu}}_{\ell,k} \coloneqq \frac{1}{M_{\ell,k}}\sum_{\alpha\in\cl{I}_{\ell,k}}\Delta_{\ell,k,\alpha},
\end{equation}
where, for each index $\ell\in\cl{L}_k$, the set of samples, $\cl{I}_{\ell,k}$, has cardinality $M_{\ell,k}$.
Finally, denote as before the difference to the previous gradient as
\begin{equation}\label{eq:Delta.sampled}
  \Delta_{\ell,k,\alpha} \coloneqq
  \left\{
  \aligned
  \grad f(\bs{\xi}_\ell, \bs{\theta}_\alpha) - \grad f(\bs{\xi}_{p_k(\ell)}, \bs{\theta}_\alpha),& \text{ if } \ell > \min\{\cl{L}_k\}, \\
  \grad f(\bs{\xi}_\ell, \bs{\theta}_\alpha),& \text{ if  } \ell = \min\{\cl{L}_k\}.
  \endaligned
  \right.
\end{equation}
\end{set}
For each $\ell \in \cl{L}_k$, we might increase the sample sizes $M_{\ell, k}$ with respect to $M_{\ell, k-1}$, hence the dependence on both $\ell$ and $k$ on the notation.
Definition~\ref{df:general.mice.estimator} allows us to manipulate the \texttt{MICE} index set to improve its efficiency; one can pick which $\ell$ to keep in $\cl{L}_k$.
For example, $\cl{L}_k = \{0, k\}$ furnishes an \texttt{SVRG}-like index set, $\cl{L}_k = \{0, 1,...,k\}$ furnishes a \texttt{SARAH}-like index set, and $\cl{L}_k = \{k\}$ results in \texttt{SGD}.
A description of these baseline variance reduction methods and their connection to \texttt{MICE} are presented in Appendix~\ref{sec:var_red_methods}.
The construction of the index set $\cl{L}_k$ is discussed in \S\ref{sec:opt_index_set}.

\begin{rmk}[Cumulative sampling in \texttt{MICE}]
\label{rmk:cum_sampling}
As the stochastic optimization progresses, new additional samples of $\bs{\theta}$ are taken and others, already available from previous iterations, are reused to compute the  \texttt{MICE} estimator at the current iteration,
\begin{equation}
  \grad \cl{F}_k
  =
  \underbrace{
  \sum_{\ell \in \cl{L}_k \cap \cl{L}_{k-1}} \frac{M_{\ell, k-1}}{M_{\ell, k}} \hat{\bs\mu}_{\ell, k-1}
  }_{\text{sunken cost}}
  \;
  +
  \;
  \underbrace{
    \sum_{\ell \in \cl{L}_k \cap \cl{L}_{k-1}} 
    \frac{1}{M_{\ell, k}} 
    \sum_{\alpha \in \cl{I}_{\ell, k} \setminus \cl{I}_{\ell, k-1}} \Delta_{\ell, k, \alpha}
    +
    \hat{\bs{\mu}}_{k, k}.
  }_{\text{Additional \texttt{MICE} cost incurred at iteration $k$}}
\end{equation}
This sampling procedure is defined by the couples $(M_{\ell,k},\bs{\xi}_{\ell})_{\ell\in\cl{L}_k}$, making $\bs{\xi}_{k+1}$ a deterministic function of all the samples in the index set $\cl{L}_k$.

\begin{rmk}[Conditional bias of the \texttt{MICE} estimator]\label{rmk:cond_bias}
As described in Remark~\ref{rmk:cum_sampling}, \texttt{MICE} reuses samples generated at previous iterations, and because the iterate path is itself a function of past samples, the estimator $\Fk$ is generally \emph{not} unbiased when conditioning on the iterate history $\xiset{k}$.
In other words, the conditional bias
\begin{equation}\label{eq:cond_bias_def}
  \bias \coloneqq \Expc{\Fk}{\xiset{k}} - \grad F(\xik)
\end{equation}
may be nonzero even though $\Fk$ is built from unbiased per-sample gradients.
Accordingly, the squared $L^2$ error admits the decomposition
\begin{equation}\label{eq:mice_err_bias_var}
  \Exp{\norm{\Fk - \grad F(\xik)}^2}
  =
  \underbrace{\Exp{\norm{\Fk - \Expc{\Fk}{\xiset{k}}}^2}}_{\text{statistical error}}
  +
  \underbrace{\Exp{\norm{\bias}^2}}_{\text{bias contribution}}.
\end{equation}
In Appendix~\ref{sec:error_decomposition} we derive an explicit expression for $\bias$ (and $\Exp{\norm{\bias}^2}$) in terms of the previously computed level estimators.
Moreover, restarting the index set $\cl{L}_k$ resets the conditional bias.
\end{rmk}

\begin{rmk}[About \texttt{MICE} and \texttt{MLMC}]
Note that \texttt{MICE} resembles the estimator obtained in the \emph{Multilevel Monte Carlo} method---\texttt{MLMC} \cite{Heinrich:MLMC,giles:MLMC,Gil15}.
For instance, if $\cl{L}_k=\{0, 1, \ldots, k\}$, \texttt{MICE} reads
\begin{align}
  \grad \cl{F}_k 
  &= \frac{1}{\M{0}} \sum_{\alpha \in \cl{I}_{0, k}} \grad f(\xio, \bs{\theta}_\alpha)
  + \sum_{\ell=1}^{k} 
  \frac{1}{\M{\ell}}
  \sum_{\alpha \in \cl{I}_{\ell, k}} \grad f(\xil, \bs{\theta}_\alpha) - \grad f(\bs{\xi}_{\ell-1}, \bs{\theta}_\alpha).
\end{align}

Indeed, we may think that in \texttt{MICE}, the iterations play the same role as the levels of approximation in  \texttt{MLMC}.
However, there are several major differences with  \texttt{MLMC}, namely \emph{i}) \texttt{MICE} exploits sunk cost of previous computations, computing afresh only what is necessary to have enough accuracy on the current iteration \emph{ii}) there is dependence in \texttt{MICE} across iterations and \emph{iii}) in \texttt{MICE}, the sample cost for the gradients is the same in different iterations while in  \texttt{MLMC} one usually has higher cost per sample for deeper, more accurate levels.

Indeed, assuming the availability of a convergent hierarchy of approximations and following the \texttt{MLMC} lines, the work \cite{multilevel_RobMon} proposed and analyzed multilevel stochastic approximation algorithms, essentially recovering the classical error bounds for multilevel Monte Carlo approximations in this more complex context.
In a similar  \texttt{MLMC} hierarchical approximation framework, the work by Yang, Wang, and Fang \cite{yang2019multilevel} proposed a stochastic gradient algorithm for solving optimization problems with nested expectations as objective functions.
Last, the combination of \texttt{MICE} and the \texttt{MLMC} ideas like those in  \cite{multilevel_RobMon} and \cite{yang2019multilevel}  is thus a natural research avenue to pursue.
\end{rmk}

\subsection{\texttt{MICE} estimator mean squared error}\label{sec:mice.error}

To determine the optimal number of samples per iteration $\ell\in\cl{L}_k$, we begin by defining the square of the error, $\cl{E}$, as the squared $L^2$-distance between \texttt{MICE} estimator \eqref{eq:mice} and the true gradient conditioned on the iterates generated up to $k$, which leads to
\begin{align}\label{eq:mse.def}
\left(\cl{E}_k\right)^2
& \coloneqq
\Expc{\norm{
\nabla_{\bs{\xi}}\cl{F}_k - \grad F(\xik)
}^2}{\xiset{k}}.
\end{align}

Next we prove that the expected $L^2$ error of the \MICE estimator is identical to the expectation of the contribution of the statistical error of each element of the index set.
Before we start, let's prove the following Lemma.
\begin{lem} \label{lem:inner_product}
Let $\muhat{\ell, k}$, as defined in \eqref{eq:mice}, be generated by a multi-iteration stochastic optimizer using \MICE as a gradient estimator.
Then, for $j \neq \ell$,
\begin{equation}
  \Exp{\inner{\muhat{\ell, k} - \bs{\mu}_{\ell, k}}{\muhat{j, k} - \bs{\mu}_{j, k}}} = 0.
\end{equation}
\end{lem}

\begin{proof}
  First, let us assume $j>\ell$ without loss of generality.
  Note that the samples of $\tht$ used to compute $\muhat{\ell, k}$ and $\muhat{j, k}$ are independent.
  However, the iterates $\{\bs{\xi}_m\}_{m=\ell+1}^j$ depend on $\muhat{\ell, k}$, thus, $\muhat{\ell, k}$ and $\muhat{j, k}$ are not independent.
  To prove Lemma~\ref{lem:inner_product}, let us use the law of total expectation to write the expectation above as the expectation of an expectation conditioned on $\xiset{j}$.
  Since $\muhat{\ell, k}$ and $\muhat{j, k}$ are then  conditionally independent,
  \begin{align}
    \Exp{\inner{\muhat{\ell, k} - \bs{\mu}_{\ell, k}}{\muhat{j, k} - \bs{\mu}_{j, k}}}
    &=
    \Exp{
      \Expc{
      \inner{
        \muhat{\ell, k} - \bs{\mu}_{\ell, k}
      }{
        \muhat{j, k} - \bs{\mu}_{j, k}
      }
      }{
        \xiset{j}
      }
    } \nonumber\\
    &=
    \Exp{
      \inner{
        \underbrace{
        \Expc{
          \muhat{\ell, k} - \bs{\mu}_{\ell, k}
        }{
          \xiset{j}
        }
        }_{\neq 0}
      }{
        \underbrace{
        \Expc{
          \muhat{j, k} - \bs{\mu}_{j, k}
        }{
          \xiset{j}
        }
        }_{=0}
      }
    },
  \end{align}
  concluding the proof.
\end{proof}

Although the level estimators are not independent unconditionally (since later iterates depend on earlier samples), Lemma~\ref{lem:inner_product} shows that the cross terms vanish after conditioning on the iterate history: the samples used at different levels are independent and each level error has zero conditional mean, hence $\Exp{\inner{\muhat{\ell, k} - \bs{\mu}_{\ell, k}}{\muhat{j, k} - \bs{\mu}_{j, k}}} = 0$ for $j\neq \ell$.

Let $\Delta_{\ell,k}^{(i)}$ be the $i$-th component of the $d_{\bs{\xi}}$ dimensional vector $\Delta_{\ell,k}$.
Then, we define
\begin{equation}\label{eq:Vl}
  V_{\ell,k}\coloneqq
  \sum_{i=1}^{d_{\bs{\xi}}}
  \bb{V}\left[\Delta_{\ell,k}^{(i)}\bigg\vert\{\bs{\xi}_{\ell'}\}_{\ell'\in\cl{L}_k} \right].
\end{equation}

\begin{lem}[Expected squared $L^2$ error of the \MICE estimator for expectation minimization]
  \label{lem:mice_expected_error}
  The expected mean squared error of the \MICE estimator is
  given by 
  \begin{align}
    \label{eq:err_expectation}
    \Exp{\left(\cl{E}_k\right)^2}
    =
    \Exp{
  \sum_{\ell \in \cl{L}_k}
  \frac{
      V_{\ell, k}
      }{M_{\ell, k}}
  },
  \end{align}
  where $\V{\ell}$ is as in \eqref{eq:Vl}.
\end{lem}
\begin{proof}
The mean squared error of the MICE estimator is
\begin{align}
  \Exp{\norm{\Fk - \grad F(\xik)}^2}
  &= \Exp{\norm{
    \sum_{\ell \in \cl{L}_k} 
    (\muhat{\ell, k} - \bs{\mu}_{\ell, k})
  }^2} \nonumber\\
  &= 
  \Exp{
    \sum_{\ell \in \cl{L}_k}
    \norm{\muhat{\ell, k} -  \bs{\mu}_{\ell, k}}^2
    } + 
  2 \sum_{\ell \in \cl{L}_k} \sum_{j \in \cl{L}_k: j>\ell}
  \Exp{
    \inner{
      \muhat{\ell, k} -  \bs{\mu}_{\ell, k}
    }{
      \muhat{j, k} -  \bs{\mu}_{j, k}
    }
  }.
\end{align}
Thus, using Lemma \eqref{lem:inner_product} and the law of total expectation,
\begin{align}
  \Exp{\norm{\Fk - \grad F(\xik)}^2}
  &= \Exp{
      \sum_{\ell \in \cl{L}_k}
      \Expc{
        \norm{\muhat{\ell, k} -  \bs{\mu}_{\ell, k}}^2
        }{\xiset{\ell}}
  }.
\end{align}
Since
\begin{align}
  \Expc{
        \norm{\muhat{\ell, k} -  \bs{\mu}_{\ell, k}}^2
        }{\xiset{\ell}}
    &=
    \label{eq:err_lvl_l}
    \Expc{
      \norm{
        \frac{1}{M_{\ell,k}}\sum_{\alpha\in\cl{I}_{\ell,k}}\Delta_{\ell,k,\alpha}
        -
        \Expc{\Delta_{\ell, k}}{\xil, \xilp}
        }^2
      }{\xil, \xilp} \\
    &=
    \frac{
      \sum_{i=1}^{d_{\bs{\xi}}}
  \bb{V}\left[\Delta_{\ell,k}^{(i)}\bigg\vert\{\bs{\xi}_{\ell'}\}_{\ell'\in\cl{L}_k} \right]
    }{
      \M{\ell}
    },
\end{align}
using \eqref{eq:Vl} concludes the proof.
\end{proof}

\begin{rmk}[Expected squared $L^2$ error of the \MICE estimator for finite sum minimization]\label{rmk:err_finite_sum}
  When minimizing a finite sum of functions as in \eqref{eq:finite_sum_problem},
  we sample the random variables $\bs{\theta}$ without replacement.
  Thus, the variance of the estimator should account for the ratio between the actual number of samples $M_{\ell,k}$ used in the estimator and the total population $N$~\cite[Section 3.7]{davison1997bootstrap}.
  In this case, the error analysis is identical to the expectation minimization case up to \eqref{eq:err_lvl_l}, except in this case we include the correction factor $(N - \M{\ell}) N^{-1}$ in the sample variance due to the finite population having size $N$, resulting in
  \begin{equation}
  \label{eq:err_finite}
  \Exp{\left(\cl{E}_k\right)^2}  = \Exp{
    \sum _{\ell\in\cl{L}_k} \frac{V_{\ell,k}}{M_{\ell,k}} \left(\frac{N-M_{\ell,k}}{N}\right)
  }.
  \end{equation}
\end{rmk}

Note that, in practice, the terms $V_{\ell,k}$ are computed using sample approximations for each $\ell\in\cl{L}_k$. 
In the convergence analysis in \S\ref{sec:mice_conv_and_work}, we assume that they are computed exactly. 
The squared $L^2$ error of \MICE can be decomposed in bias and statistical error, which are analyzed in Appendix~\ref{sec:error_decomposition}.

\subsection{Multi-iteration optimal setting for gradient error control} \label{sec:mice.optimal.setting}

First, let the gradient sampling cost and the total \texttt{MICE} work be defined as follows.
The number of gradient evaluations is $1$ for $\Delta_{\ell, k, \alpha}$ when $\ell=\min\{\cl{L}_k\}$ and $2$ otherwise.
For this reason, we define the auxiliary index function
\begin{equation}
\mathbbm{1}_{\mskip-2mu\underline{\cl{L}_k}}(\ell)\coloneqq
\begin{cases}
0 & \text{if}\quad\ell=\min\{\cl{L}_k\},\\
1 & \text{otherwise}
\end{cases},
\end{equation}
and define the gradient sampling cost in number of gradient evaluations as
\begin{equation}\label{eq:grad_sampling_cost}
  \cl{C}(\{M_{\ell,k}\}_{\ell\in\cl{L}_k})
  \coloneqq
  \sum_{\ell\in\cl{L}_k} \cost{\ell} M_{\ell,k}.
\end{equation}

Motivated by the analysis of \texttt{SGD-MICE} in \S\ref{sec:mice_conv_and_work}, here we choose the number of samples for the index set $\cl{L}_k$ by approximate minimization of the gradient sampling cost \eqref{eq:grad_sampling_cost} subject to a given tolerance $\epsilon>0$ on the relative error in the mean gradient approximation, that is
\begin{equation}\label{eq:M_opt_problem}
\begin{aligned}
\{M^*_{\ell,k}\}_{\ell\in\cl{L}_k}
 = & \underset{\{M_{\ell,k}\}_{\ell\in\cl{L}_k}}{\text{arg min}}
\cl{C}\left(
\{M_{\ell,k}\}_{\ell\in\cl{L}_k}
\right) \\
& \text{subject to} \quad (\cl{E}_k)^2 \le \epsilon^2\norm{\nabla_{\bs{\xi}}F(\bs{\xi}_k)}^2.
\end{aligned}
\end{equation}

\subsubsection{Expectation minimization}

As a consequence of \eqref{eq:M_opt_problem} and Lemma~\ref{lem:mice_expected_error}, we define the sample sizes as the solution of the following constrained optimization problem,
\begin{equation}\label{eq:M_relstatprob}
  \begin{aligned}
  \{M^*_{\ell,k}\}_{\ell\in\cl{L}_k}
   = & \underset{\{M_{\ell,k}\}_{\ell\in\cl{L}_k}}{\text{arg min}}
  \cl{C}\left(
  \{M_{\ell,k}\}_{\ell\in\cl{L}_k}
  \right) \\
  & \text{subject to} \quad \sum_{\ell \in \cl{L}_k} \frac{V_{\ell, k}}{M_{\ell, k}} \le 
  \epsilon^2\norm{\nabla_{\bs{\xi}}F(\bs{\xi}_k)}^2.
  \end{aligned}
  \end{equation}

An approximate integer-valued solution based on Lagrangian relaxation to problem \eqref{eq:M_relstatprob} is
\begin{equation}\label{eq:optimal.sampling.1}
M^*_{\ell,k} = \left\lceil
\frac{1}{\epsilon^2 \norm{\grad F(\xik)}^2}
\left( \sum _{\ell^\prime\in\cl{L}_k} \sqrt{V_{\ell^\prime\!,k}\cost{\ell^\prime}} \right) 
\sqrt{\dfrac{V_{\ell,k}}{\cost{\ell}}}
\, \right\rceil,\qquad \forall \ell\in\cl{L}_k.
\end{equation}
In general, in considering the cost of computing new gradients at the iteration $k$, the expenditure already carried out up to the iteration $k-1$ is sunk cost and must not be included, as described in Remark~\ref{rmk:cum_sampling}, that is, one should only consider the incremental cost of going from $k-1$ to $k$. 
Moreover, in the variance constraint of problem \eqref{eq:M_relstatprob}, since we do not have access to the  norm of the  mean gradient, $\norm{\nabla_{\bs{\xi}}F(\bs{\xi}_k)}$, we use a resampling technique combined with the \texttt{MICE} estimator as an approximation; see Remark~\ref{rmk:grad_resampling}.

\subsubsection{Finite sum minimization}

In view of Remark~\ref{rmk:err_finite_sum}, we define the sample sizes for \MICE[] as the solution of the following optimization problem,
\begin{align}
\text{find } \{\M{\ell}^*\}
= \underset{\{\M{\ell}\}_{\ell \in \cl{L}_k}}{\arg \min}
&\cl{C}\left(
  \{M_{\ell,k}\}_{\ell\in\cl{L}_k}
  \right) \\
\nonumber
\text{subject to } & 
\begin{cases}
    \sum_{\ell \in \cl{L}_k} 
    \frac{\V{\ell}}{\M{\ell}} - \frac{\V{\ell}}{N} \le 
    \epsilon^2 \gnorm{k}^2 \\
    M_{min} \le \M{\ell} \le N \quad \forall \ell \in \cl{L}_k.
\end{cases}
\end{align}
This problem does not have a closed form solution, but can be solved in an iterative process by noting that any $\ell$ such that $\M{\ell} = N$ does not contribute to the error of the estimator.
Then, letting
\begin{equation}
  \cl{G}_k = \{\ell \in \cl{L}_k : M_\ell < N\},
\end{equation}
we derive a closed form solution for the sample sizes as
\begin{equation}
    \label{eq:sample_size_finite}
    \M{\ell}^* = 
    \left\lceil
    \,
    \frac{
        \sum_{\ell'\in\cl{G}_k} \sqrt{\cost[\ell']{k} \V{\ell}}
    }{
        \epsilon^2 \gnorm{k}^2
        + N^{-1}\sum_{\ell'' \in \cl{G}_k} \V{\ell''}
    }
    \sqrt{\frac{\V{\ell}}{\cost{\ell}}}
    \,
    \right\rceil
    .
\end{equation}
However, it is not possible to know directly the set $\cl{G}_k$.
So, we initialize $\cl{G}_k = \cl{L}_k$ and iteratively remove elements that do not satisfy the condition $M_\ell < N$ as presented in Algorithm~\ref{alg:finite_sum}.
\end{rmk}

\begin{algorithm}[h]
  \setstretch{1.5}
  \caption{Computing sample size of SGD-MICE for the finite sum case.
  }\label{alg:finite_sum}
  \begin{algorithmic}[1]
  \State{$\cl{G}_k \gets \cl{L}_k$}
  \State{Set $\M{\ell}$ using \eqref{eq:sample_size_finite} for all $\ell \in \cl{G}_k$}
  \While{\text{any} $\{\ell \in \cl{G}_k: \M{\ell} \ge N\}$}
  \For{$\ell \in \{\ell \in \cl{G}_k: \M{\ell} \ge N\}$}
  \State $\M{\ell} \gets N$
  \State $\cl{G}_k \gets \cl{G}_k \setminus \{\ell\}$
  \EndFor
  \State{Set $\M{\ell}$ using \eqref{eq:sample_size_finite} for all $\ell \in \cl{G}_k$}
  \EndWhile
  \State \text{\bf{Return}} $\{\lceil\M{\ell}\rceil\}_{\ell \in \cl{L}_k}$
  \end{algorithmic}
  \end{algorithm}

\subsection{Optimal index set operators}
\label{sec:opt_index_set}

As for the construction of the \texttt{MICE} index set at iteration $k$, that is, $\cl{L}_k$, from the previous one, $\cl{L}_{k-1}$, we use one of the following index set operators:
\begin{set}\label{df:opt_index_set_operators}[Construction of the index set $\cl{L}_k$]
For $k=0$, let $\cl{L}_0 = \{0\}$. If $k\ge 1$,
After this step, there are four possible cases to finish the construction of $\cl{L}_k$:
\begin{center}
\begin{tabular}{llll}
\Add:         &$\cl{L}_k$   &$\leftarrow   \cl{L}^{\text{add}}_{k}$    & $\gets \cl{L}_{k-1}\cup\{k\}$ \\
\Drop:        &$\cl{L}_k$   &$\leftarrow   \cl{L}^{\text{drop}}_{k}$   & $\gets \cl{L}_{k-1}\cup\{k\}\!\setminus\{k-1\}$ \\
\Restart:     &$\cl{L}_k$   &$\leftarrow   \cl{L}^{\text{rest}}_{k}$   & $\gets \{k\}$ \\
\Clip at $\ell^*$: &$\cl{L}_k$   &$\leftarrow   \cl{L}^{\text{clip,}\ell^*}_{k}$ & $\gets \cl{L}_{k-1}\cup\{k\}\!\setminus\{\ell \in \cl{L}_{k-1}: \ell < \ell^*\}$
\end{tabular}
\end{center}
\end{set}
The \Add operator simply adds $k$ to the current index set.
The \Drop operator does the same but also removes $k-1$ from the index set.
As the name suggests, \Restart resets the index set at the current iterate.
Finally, \Clip adds $k$ to the current index set and removes all components previous to $j$.
For more details, see \S\ref{sec:mice_algorithm} for an algorithmic description.

In the previous section, the sample sizes for each element of the index set are chosen as to minimize the gradient sampling cost while satisfying a relative error constraint.
However, to pick one of the operators to update the index set, we must use the work including the overhead of aggregating the index set elements.
Let the gradient sampling cost increment at iteration $k$ be
\begin{equation}
  \label{eq:cost_of_iter_k}
  \Delta \cl{C}_k(\cl{L})
  =
  \sum_{\ell \in \cl{L} \cap \cl{L}_{k-1}} \cost{\ell}(\M[k]{\ell}^* - \M[k-1]{\ell})
  + \cost{k} \M{k}^*,
\end{equation}
with $\M{\ell}^*$ as in \eqref{eq:optimal.sampling.1} or Algorithm~\ref{alg:finite_sum}.

The total work of a \texttt{MICE} evaluation is then the sum of the cost of sampling the gradients and the cost of aggregating the gradients as
\begin{equation}
  \cl{W}\left(
  \{M_{\ell,k}\}_{\ell\in\cl{L}_k}
  \right)
  \coloneqq \cl{C}(
  \{M_{\ell,k}\}_{\ell\in\cl{L}_k}
  ) C_\nabla
  + |\cl{L}_k| C_\text{aggr},
\end{equation}
where $C_\nabla$ is the work of sampling $\nabla_{\bs{\xi}} f$ and $C_\text{aggr}$ is the work of averaging the $\Delta_{\ell, k}$ to construct $\cl{F}_k$.
Then, the work done in iteration $k$ to update \texttt{MICE} is
\begin{equation}\label{eq:work_update_MICE}
  \Delta \cl{W}(\cl{L})
  \coloneqq
  \Delta \cl{C}_k(\cl{L})
   C_\nabla
  + |\cl{L}| C_\text{aggr}.
\end{equation}
We choose the index set operator for iteration $k$ by a greedy work-based policy.
We then accept drop/restart operations if they are not significantly more expensive than continuing, and accept clipping if it strictly reduces the work:
\begin{equation}
  \label{eq:work_incr}
  \begin{aligned}
  &\text{Initialize}\quad \cl{L}_k \gets \cl{L}^{\text{add}}_{k},\\
  & \text{(Add) Let } \Delta\cl{W}_{\text{add}}\coloneqq \Delta \cl{W}_k(\cl{L}^{\text{add}}_{k}) \\
  &\text{(Drop) if}\quad \Delta \cl{W}_k(\cl{L}^{\text{drop}}_{k})
  \le (1+\delta_{\text{drop}})\, \Delta\cl{W}_{\text{add}}
  \quad\text{then set}\quad \cl{L}_k \gets \cl{L}^{\text{drop}}_{k},\\
  &\text{(Restart) if}\quad \Delta \cl{W}_k(\cl{L}^{\text{rest}}_{k})
  < (1+\delta_{\text{rest}})\, \Delta \cl{W}_k(\cl{L}_k)
  \quad\text{then set}\quad \cl{L}_k \gets \cl{L}^{\text{rest}}_{k},\\
  &\text{(Clip) if}\quad \min_{\ell \in \cl{L}_{k-1}} \Delta \cl{W}_k(\cl{L}^{\text{clip,}\ell}_{k})
  < \Delta \cl{W}_k(\cl{L}_k)
  \quad\text{then clip at the minimizing }\ell,
  \end{aligned}
\end{equation}
and otherwise we keep $\cl{L}_k$ as selected by the checks above (initialized as $\cl{L}^{\text{add}}_{k}$).
Here $\cl{L}^{\text{clip,}\ell}_k$ is discussed in more detail in \S\ref{sec:mice.clipping}, and $\delta_{\text{drop}}, \delta_{\text{rest}} \ge 0$ are slack parameters: increasing them makes the corresponding operator more likely to be selected, which tends to keep the index set smaller and reduce \MICE[]'s overhead.
In our experiments we use values in $[0,1]$.

\subsubsection{Dropping iterations of the \texttt{MICE} index set}\label{sec:mice.dropping}

Given our estimator's stochastic nature, at the current iteration $k$, we may wonder if the iteration $k -1$ should be kept or dropped out from the \texttt{MICE} index set since it may not reduce the computational work.
The procedure we follow here draws directly from an idea introduced by Giles \cite{Gil15} for the \texttt{MLMC} method.
Although the numerical approach is the same, we construct the algorithm in a greedy manner. We only check the case of dropping the previous iteration in the current index set. In this approach, we never drop the initial iteration $\min\{\cl{L}_k\}$.

\subsubsection{Restarting the \texttt{MICE} index set}\label{sec:mice.restart}
As we verified in the previous section on whether we should keep the iteration $\ell = k -1$ in the \texttt{MICE} index set, we also may wonder if restarting the estimator may be less expensive than updating it.
Usually, in the literature of control variates techniques for stochastic optimization, the restart step is performed after a fixed number of iterations; see, for instance, \cite{Joh13,Nit14,Kon16}.
Moreover, restarting the index set resets the conditional bias discussed in Remark~\ref{rmk:cond_bias}.

\subsubsection{Clipping the \texttt{MICE} index set} \label{sec:mice.clipping}

In some cases, it may be advantageous to discard only some initial iterates indices out of the index set instead of the whole index set.
We refer to this procedure as clipping the index set.
We propose two different approaches to decide when and where to clip the index set.
\begin{description}
  \item[Clipping ``A''] $\cl{L}^{\text{clip,}\ell^*}_{k}$ is as in Definition~\ref{df:opt_index_set_operators} with
  \begin{align}
    \ell^*
    =
    \underset{
      \ell \in \cl{L}_{k-1}
    }{
      \arg \min
    }
    \;
    \Delta \cl{W}_k(\cl{L}^{\text{clip,}\ell}_{k}).
  \end{align}
  This clipping technique can be applied in both the continuous and discrete cases.
  \item[Clipping ``B''] This technique is simpler but can only be used in the finite sum case. It consists in clipping $\cl{L}_{k-1}$ at $\ell^* = \max\{\ell \in \cl{L}_{k-1}: \M[k-1]{\ell} = N\}$.
\end{description}

Clipping ``A'' adds an extra computation overhead when calculating $M_{\ell, k}$ for each $\ell \in \cl{L}_k$ each iteration $k$.
Thus, in the finite sum case, we suggest using Clipping ``B''.
Clipping shortens the index set, thus possibly reducing the general overhead of \texttt{MICE}.
Moreover, clipping the index set may reduce the frequency of restarts and the bias of the \texttt{MICE} estimator.

\section{\texttt{SGD-MICE} convergence and gradient sampling cost analysis} \label{sec:mice_conv_and_work}

In this section, we will analyze the convergence of stochastic gradient methods with fixed step size as
\begin{equation}
  \xikp = \xik - \eta \bs{\upsilon}_k,
\end{equation}
with gradient estimates controlled as
\begin{align}
  \label{eq:err_condition}
  \Exp{\norm{\bs\upsilon_k - \grad F(\xik)}^2}
  \le \epsilon^2 \Exp{\norm{\grad F(\xik)}^2}.
\end{align}
In special, we are interested in \texttt{SGD-MICE}, where $\bs\upsilon_k = \grad \cl{F}_k$ as defined in~\eqref{eq:mice}, and \texttt{SGD-A}, where $\bs\upsilon_k = M_k^{-1} \sum_{i=1}^{M_k} \grad f(\xik, \bs{\theta}_i)$.
Here, \texttt{SGD-A} is SGD where the sample sizes are increased to control the statistical error condition in~\eqref{eq:err_condition} and can be seen as a special case of \texttt{SGD-MICE} where \Restart is used every iteration.
For \MICE[], this condition is satisfied by the choice of the sample sizes in \S\ref{sec:mice.optimal.setting}.

Let us lay some assumptions.
\begin{ass}[Lipschitz continuous gradient]
  \label{as:lipschitz}
  If the gradient of $F \colon \bb{R}^{d_{\bs{\xi}}}\mapsto\bb{R}$ is Lipschitz continuous, then, for some $L > 0$,
  \begin{equation}\label{eq:lipschitz}
  \norm{\nabla_{\bs{\xi}} F(\bs{x}) - \nabla_{\bs{\xi}} F(\bs{y})}{} \le L \norm{\bs{y} - \bs{x}}{}, \qquad \forall \bs{x}, \bs{y} \in \rset^{d_{\bs{\xi}}}.
  \end{equation}
\end{ass}

\begin{ass}[Convexity]\label{as:convex}
  If $F$ is convex, then,
  \begin{equation}\label{eq:convex}
  F(\bs{y}) \ge F(\bs{x}) + \inner{\grad F(\bs{x})}{\bs{y} - \bs{x}}, \qquad \forall\,\bs{x}, \bs{y} \in \rset^{d_{\bs{\xi}}}.
  \end{equation}
\end{ass}

\begin{ass}[Strong convexity]\label{as:strongly_convex}
  If $F$ is $\mu$-strongly convex, then, for some $\mu>0$,
  \begin{equation}\label{eq:strongly_convex}
  F(\bs{y}) \ge F(\bs{x}) + \langle \nabla_{\bs{\xi}} F(\bs{x}), \bs{y} - \bs{x} \rangle + \dfrac{\mu}{2} \norm{\bs{y} - \bs{x}}^2, \qquad \forall\,\bs{x}, \bs{y} \in \rset^{d_{\bs{\xi}}}.
  \end{equation}
\end{ass}

\begin{ass}[Polyak--{\L}ojasiewicz]\label{as:polyak}
  If $F$ is gradient dominated, it satisfies the Polyak--{\L}ojasiewicz inequality
  \begin{equation}
      \frac{1}{2} \norm{\grad F(\bs{x})}^2 \ge \mu (F(\bs{x}) - F^*), \qquad 
      \forall\,\bs{x} \in \rset^{d_{\bs{\xi}}}, 
  \end{equation}
  for a constant $\mu>0$, where $F^*$ is the minimum value of $F$.
\end{ass}
Assumption~\ref{as:polyak} is weaker than Assumption~\ref{as:strongly_convex}, holding even for some non-convex problems~\cite{karimi2016linear}.

\subsection{Optimization convergence analysis}
\label{sec:convergence}
\begin{prop}[Local convergence in expectation of gradient-controlled \texttt{SGD} on $L$-smooth problems]\label{prp:smooth_convergence}
  Let $F:\bb{R}^{d_\xi} \rightarrow \bb{R}$ be a differentiable function satisfying Assumption~\ref{as:lipschitz} with constant $L>0$.
  Then, \texttt{SGD} methods with relative gradient error control $\epsilon<1$ in the $L^2$-norm sense and step-size $\eta = 1/L$ reduce the optimality gap in expectation as
  \begin{equation} \label{eq:conv_convex}
      \Exp{F(\xikp)}
      \le \Exp{F(\xik)}
       - \left(\frac{1 - \epsilon^2}{2L}\right) \Exp{\norm{\grad F(\xik)}^2}.
 \end{equation}

\end{prop}
\begin{proof}
  Let $\bs{e}_k \coloneqq \bs{\upsilon}_k - \grad F(\xik)$.
  From $L$-smoothness,
  \begin{align}
      F(\xikp)
      &\le 
      F(\xik) - \eta \inner{\grad F(\xik)}{\grad F(\xik) + \bs{e}_k}
      + \frac{L\eta^2}{2} \norm{\grad F(\xik) + \bs{e}_k}^2 \\
      & =
      F(\xik)
      + \left(\frac{L \eta^2}{2} - \eta\right) \norm{\grad F(\xik)}^2
      + (L \eta^2 - \eta) \inner{\grad F(\xik)}{\bs{e}_k}
      + \frac{L \eta^2}{2} \norm{\bs{e}_k}^2.
  \end{align}
  Taking expectation on both sides and then using the Cauchy--Schwarz inequality,
    \begin{align}
    \Exp{F(\xikp)}
    &\le
    \!\begin{multlined}[t][10.5cm]
        \Exp{F(\xik)}
        + \left(\frac{L \eta^2}{2} - \eta\right) \Exp{\norm{\grad F(\xik)}^2}
        + |L \eta^2 - \eta| \sqrt{\Exp{\norm{\grad F(\xik)}^2}\Exp{\norm{\bs{e}_k}^2}} \\
        + \frac{L \eta^2}{2} \Exp{\norm{\bs{e}_k}^2}.
    \end{multlined} \\
    &\le \Exp{F(\xik)}
    + \left(
        \frac{L \eta^2}{2} - \eta
        + \epsilon |L \eta^2 - \eta|
        + \epsilon^2 \frac{L \eta^2}{2}
    \right) \Exp{\norm{\grad F(\xik)}^2},
    \end{align}
    where \eqref{eq:err_condition} is used to get the last inequality.
    Here, the step size that minimizes the term inside the parenthesis is $\eta=1/L$. 
    Substituting the step size in the equation above and taking full expectation on both sides concludes the proof.
\end{proof}
If the function $F$ is also unimodal, as in the case of $F$ satisfying Assumptions~\ref{as:convex} or \ref{as:polyak}, then the convergence presented in Proposition~\ref{prp:smooth_convergence} is also global, i.e., $\Exp{F(\xikp) - F(\opt)} \rightarrow 0$.

\begin{prop}[Global convergence of gradient-controlled \texttt{SGD} in gradient-dominated problems]
  \label{prp:pl_convergence}
  Let all Assumptions of Proposition~\ref{prp:smooth_convergence} be satisfied. Moreover, let $F$ satisfy Assumption~\ref{as:polyak} with constant $\mu>0$.
  Then, gradient-controlled \texttt{SGD} with step-size $\eta=1/L$ converges linearly,
  \begin{equation}
      \label{eq:linear_convergence_opt_gap}
      \Exp{F(\xikp) - F(\opt)}
      \le 
      \left(1 - (1-\epsilon^2)\frac{\mu}{L}\right)^{k+1} \Exp{F(\xio) - F(\opt)}.
  \end{equation}
\end{prop}
\begin{proof}
  From \eqref{eq:conv_convex}, using Assumption~\ref{as:polyak} and unrolling the recursion concludes the proof.
\end{proof}

\begin{cor}\label{cor:grad_norm_convergence}
   If conditions of Proposition~\ref{prp:pl_convergence} hold and $F$ also satisfies Assumption~\ref{as:convex}, the squared $L^2$-norm of the gradient of the objective function is bounded as
   \begin{equation}
      \Exp{\norm{\grad F(\xikp)}^2}
      \le 2 L
      \left(1 - (1-\epsilon^2)\frac{\mu}{L}\right)^{k+1} \Exp{F(\xio) - F(\opt)}.
   \end{equation}
\end{cor}
\begin{proof}
    From \cite[Theorem 2.1.5]{nesterov2018lectures}, if $F$ is convex and $L$-smooth,
    \begin{equation}
        \label{eq:grad_bound}
        \norm{\grad F(\bs{\xi})}^2 \le 2L (F(\bs{\xi}) - F(\opt)).
    \end{equation}
    Substituting this inequality for $\xikp$ into \eqref{eq:linear_convergence_opt_gap} finishes the proof.
\end{proof}

\begin{cor}\label{cor:dist_to_opt_convergence}
  If conditions of Proposition~\ref{prp:pl_convergence} are satisfied and $F$ also satisfies Assumption~\ref{as:strongly_convex},
  \begin{equation}
      \Exp{\norm{\xikp - \opt}^2} \le \frac{2}{\mu}
      \left(1 - (1-\epsilon^2)\frac{\mu}{L}\right)^{k+1} \Exp{F(\xio) - F(\opt)}.
  \end{equation}
\end{cor}
\begin{proof}
    From the definition of strong-convexity in Assumption~\ref{as:strongly_convex},
    \begin{equation}
        \norm{\bs{\xi} - \opt}^2 \le \frac{2}{\mu} (F(\bs{\xi}) - F(\opt)).
    \end{equation}
    Substituting into \eqref{eq:linear_convergence_opt_gap} finishes the proof.
\end{proof}

\subsubsection{High-probability convergence in the PL regime}
\label{sec:hp_main}
We complement the expectation-based results above with a high-probability statement.
The proof separates into two steps:
(\emph{i}) show linear convergence on any event where the relative error holds uniformly in $k$ and
(\emph{ii}) show that, under a simple tail assumption (Assumption~\ref{as:hp_coord_subg}), such an event has probability at least $1-\delta$ in an add-only \texttt{MICE} regime (Appendix~\ref{sec:hp_mice}).

\begin{ass}[Coordinatewise sub-Gaussian level increments]\label{as:hp_coord_subg}
In an add-only regime, assume that each level $\ell\in\bb{N}$ is anchored to a fixed iterate pair between restarts and its contribution to $\Fk$ is a sample mean of i.i.d.\ centered increments $Z_{\ell,1},Z_{\ell,2},\ldots\in\bb{R}^{d_\xi}$ with $\bb{E}[Z_{\ell,i}]=0$.
Assume that for each coordinate $j\in\{1,\ldots,d_\xi\}$ there exists $\sigma_\ell^2\ge 0$ such that for all $\lambda\in\bb{R}$,
\begin{equation}
  \label{eq:hp_subg_mgf}
  \bb{E}\!\left[\exp\!\big(\lambda (Z_{\ell,i})_j\big)\right]
  \le
  \exp\!\Big(\tfrac{\lambda^2\sigma_\ell^2}{2}\Big).
\end{equation}
\end{ass}

\begin{thm}[Linear convergence on a uniform relative error event]
\label{thm:hp_pl_event}
Assume $F$ satisfies Assumptions~\ref{as:lipschitz} and~\ref{as:polyak} with constants $L>0$ and $\mu>0$.
Fix $\eta>0$ and $\epsilon\in[0,1)$.
Define
\begin{equation}
  \label{eq:c_eta_eps}
  c(\eta,\epsilon)
  \coloneqq
  \eta(1-\epsilon) - \frac{L\eta^2}{2}(1+\epsilon)^2.
\end{equation}
Assume $c(\eta,\epsilon)>0$ and $2\mu\,c(\eta,\epsilon)<1$.
Let $\Omega$ be any event on which the realized gradient estimation errors satisfy
\begin{equation}
  \label{eq:hp_rel_err_event}
  \norm{\bs{\upsilon}_k - \grad F(\xik)} \le \epsilon \norm{\grad F(\xik)},
  \qquad \forall\,k\ge 0.
\end{equation}
Then, on $\Omega$, the iterates of \texttt{SGD} with update $\xikp=\xik-\eta\bs{\upsilon}_k$ satisfy, for all $k\ge 0$,
\begin{equation}
  \label{eq:hp_linear_rate}
  F(\xik) - F(\opt)
  \le
  r_{\mathrm{hp}}^k \big(F(\xio) - F(\opt)\big),
  \qquad
  r_{\mathrm{hp}} \coloneqq 1-2\mu\,c(\eta,\epsilon)\in(0,1).
\end{equation}
\end{thm}
\begin{proof}
Work on $\Omega$ and fix $k\ge 0$.
From $L$-smoothness and the update $\xikp=\xik-\eta(\grad F(\xik)+(\bs{\upsilon}_k-\grad F(\xik)))$, we have
\begin{align}
  F(\xikp)
  &\le
  F(\xik)
  -\eta \langle \grad F(\xik), \bs{\upsilon}_k\rangle
  +\frac{L\eta^2}{2}\norm{\bs{\upsilon}_k}^2 \nonumber\\
  &\le
  F(\xik)
  -\eta\norm{\grad F(\xik)}^2
  +\eta\norm{\grad F(\xik)}\,\norm{\bs{\upsilon}_k-\grad F(\xik)}
  +\frac{L\eta^2}{2}\big(\norm{\grad F(\xik)}+\norm{\bs{\upsilon}_k-\grad F(\xik)}\big)^2 \nonumber\\
  &\le
  F(\xik)
  -c(\eta,\epsilon)\,\norm{\grad F(\xik)}^2.
  \label{eq:hp_descent_step}
\end{align}
where in the last step we used \eqref{eq:hp_rel_err_event}.
Applying the PL inequality and unrolling the recursion yields \eqref{eq:hp_linear_rate}.
\end{proof}

\begin{cor}[High-probability linear convergence for add-only \texttt{MICE}]
\label{cor:hp_pl_mice}
Assume the conditions of Theorem~\ref{thm:hp_pl_event}.
Consider \texttt{SGD-MICE} in an add-only regime and assume Assumption~\ref{as:hp_coord_subg}.
Assume further that, conditional on the iterate history, the sample sets used at distinct levels are independent, and that the sample sizes $\{M_{\ell,k}\}_{\ell\in\cl{L}_k}$ are chosen before drawing the iteration-$k$ samples.
Fix a summable schedule $(\delta_k)_{k\ge 0}$ with $\sum_{k\ge 0}\delta_k\le \delta$.
If the (predictable) sample sizes satisfy the scalar constraint
\begin{equation}
  \label{eq:hp_variance_sum_constraint}
  \sum_{\ell\in\cl{L}_k}\frac{\sigma_\ell^2}{M_{\ell,k}}
  \le
  \frac{\epsilon^2}{2d_\xi\,\log(2d_\xi/\delta_k)}\,\norm{\grad F(\xik)}^2,
  \qquad \forall\,k\ge 0,
\end{equation}
then, with probability at least $1-\delta$, the bound \eqref{eq:hp_linear_rate} holds for all $k\ge 0$.
\end{cor}

\begin{proof}
  Define the gradient estimation error $e_k\coloneqq \bs{\upsilon}_k-\nabla F(\xik)$.
  Under Assumption~\ref{as:hp_coord_subg}, the conditional independence across levels, and the predictability of the sample sizes,
  Lemma~\ref{lem:hp_mice_single_event} (in Appendix~\ref{sec:hp_mice}) yields an event $\Omega^{\mathrm{MICE}}_\delta$ such that
  $\bb{P}[\Omega^{\mathrm{MICE}}_\delta]\ge 1-\delta$ and, on $\Omega^{\mathrm{MICE}}_\delta$, for all $k\ge 0$,
  \begin{equation}
  \label{eq:hp_mice_err_bd_recalled}
  \norm{e_k}
  \le
  \sqrt{2d_\xi\log\!\Big(\frac{2d_\xi}{\delta_k}\Big)}
  \left(\sum_{\ell\in\cl{L}_k}\frac{\sigma_\ell^2}{M_{\ell,k}}\right)^{1/2}.
  \end{equation}
  (The summability $\sum_k\delta_k\le \delta$ is used in Lemma~\ref{lem:hp_mice_single_event} via a union bound over $k$.)
  
  Assume now that the sample sizes satisfy \eqref{eq:hp_variance_sum_constraint}. Substituting this bound into
  \eqref{eq:hp_mice_err_bd_recalled} gives, on $\Omega^{\mathrm{MICE}}_\delta$, for all $k\ge 0$,
  \begin{equation}    
    \norm{e_k}\le \epsilon \norm{\nabla F(\xik)}.
  \end{equation}
  Hence the uniform relative-error condition \eqref{eq:hp_rel_err_event} of Theorem~\ref{thm:hp_pl_event} holds on
  $\Omega^{\mathrm{MICE}}_\delta$. Applying Theorem~\ref{thm:hp_pl_event} concludes that \eqref{eq:hp_linear_rate}
  holds for all $k\ge 0$ with probability at least $1-\delta$.
  \end{proof}

The high-probability result above is intentionally stated for the add-only regime and for sample sizes $M_{\ell,k}$ chosen predictably with respect to the past.
Extending it to fully adaptive index-set operations (Drop/Restart/Clip) and within-iteration stopping rules is left for future work.

Having established convergence rates for gradient-controlled SGD methods, we now quantify the total number of gradient evaluations required to achieve a given tolerance, comparing \texttt{SGD-MICE} with adaptive batch-size approaches.

\subsection{Gradient sampling cost analysis}
\label{sec:sampling_cost}
To analyze the gradient sampling cost, we focus on the analysis in expectation with $L^2$ control on the error.
A discussion on cost in high probability is presented after Corollary~\ref{cor:cost_sgd_mice} in Remark~\ref{rmk:cost_expect_vs_hp}.
Assuming the assumptions of Proposition~\ref{prp:pl_convergence} hold, the optimality gap converges with rate $r \coloneqq 1 - (1 - \epsilon^2) \mu / L$.
Then, we have the following inequalities that will be used throughout this section,
\begin{align}
    \label{eq:rate_ineq_1}
    \frac{1}{\log(r)}
    \le
    \frac{1}{1 - r} 
    =
    \frac{\kappa}{1-\epsilon^2},
\end{align}
where $\kappa = L / \mu$.
Moreover,
\begin{align}
    \label{eq:rate_ineq_2}
    \frac{1}{1 - \sqrt{r}}
    \le
    \frac{2 \kappa}{1-\epsilon^2}.
\end{align}

For the sake of simplicity and given the cumulative nature of the computational gradient sampling cost in \texttt{MICE}, we analyze the total gradient sampling cost on a set of iterations $\{\bs{\xi}_\ell\}_{\ell=0}^{k^*}$ converging to $\bs{\xi}^\ast$ as per Proposition~\ref{prp:pl_convergence}.
Observe that in this simplified setting, the number of iterations required to stop the iteration, $k^*= k^*(tol)$, and both the sequences $(\bs{\xi}_\ell)$ and $(M_{\ell,k})$  are still random. Indeed, we define
\begin{equation}
  \label{eq:stopping_criterion}
  k^*= \min\{k\ge0\colon \norm{\grad F(\xik)}^2\le tol\}.
\end{equation}

\begin{cor}[Number of iterations]\label{cor:convergence.seq}
If the assumptions of Corollary~\ref{cor:grad_norm_convergence} hold then, letting
\begin{equation}
k_1 \coloneqq \frac{\log(tol^{-1} 2 L\Exp{F(\xio) - F(\opt)})}{\log(1/r)},
\end{equation}
we have
\begin{equation}\label{eq:distrib}
\bb{P}\left[k^* \ge k\right] \le \left\{
\begin{aligned}
1, & \text{ if }k< k_1 \\
r^{ k-k_1}& \text{ otherwise.}
\end{aligned}
\right.
\end{equation}
Moreover, we have
\begin{equation}\label{eq:expect}
\bb{E}\left[k^*\right] \le \frac{1}{1-r} +
    \max\left\{0,\frac{
            \log(tol^{-1} 2 L \Exp{F(\xio) - F(\opt)})
        }{
            \log(1/r)
        }\right\}.
\end{equation}
\end{cor}
\begin{proof}
First observe that
\begin{equation}
\bb{P}\left[k^* \ge k\right] \le \bb{P}\left[\norm{\grad F(\xik)}^2 \ge tol\right].
\end{equation}
Then apply Markov's inequality and the exponential convergence in $L^2$-norm presented in Corollary~\ref{cor:grad_norm_convergence}, yielding
\begin{equation}\label{eq:distrib1}
\bb{P}\left[k^* \ge k\right] \le \min\left\{1,
tol^{-1} 2L r^k \Exp{F(\xio) - F(\opt)} \right\}.
\end{equation}
The result \eqref{eq:distrib} follows then directly.
To show \eqref{eq:expect}, simply use \eqref{eq:distrib1} and that
\begin{equation}\label{eq:distrib2}
\bb{E}\left[k^*\right] = \sum_{k\ge 0 }\bb{P}\left[k^* \ge k\right] \le \max\{0,k_1\}+\frac{1}{1-r}.
\end{equation}
\end{proof}

The expected value of $k^*$ can be bounded using \eqref{eq:rate_ineq_1} as
\begin{equation}
    \bb{E}\left[k^*\right]
    \le
    \max\left\{0,
    \frac{\kappa}{1-\epsilon^2} \log(tol^{-1} 2 L\Exp{F(\xio) - F(\opt)})\right\}
    + \frac{\kappa}{1-\epsilon^2}.
\end{equation}

\begin{ass}[Bound on second moments of gradient differences]
    \label{as:variance}
    \begin{equation}
        \Expc{\norm{\grad f(\bs{x}, \bs{\theta}) - \grad f(\bs{y}, \bs{\theta})}^2}{\bs{x}, \bs{y}}
        \le \sigma^2 \norm{\grad F(\bs{x}) - \grad F(\bs{y})}^2.
    \end{equation}
\end{ass}

Assumption~\ref{as:variance} is a convenient way to control the second moments of gradient differences by differences of true gradients, and it holds for many smooth models with light-tailed gradient noise.
However, in heavy-tailed regimes (or when gradients are contaminated by occasional large outliers), second moments may be very large or even infinite, so Assumption~\ref{as:variance} may fail or yield overly conservative bounds.
In such cases, a practical safeguard is to clip large gradients or gradient differences (not to be confused with our \Clip[] operator), which limits the influence of outliers at the cost of introducing a controlled bias.

If $f$ satisfies Assumption~\ref{as:variance} for $\ell>0$,
\begin{align}
    \V{\ell} &\le \Expc{\norm{\grad f(\xil, \bs{\theta}) - \grad f(\xilp, \bs{\theta})}^2}{\xil, \xilp} \\
    &\le \sigma^2 \norm{\grad F(\xil) - \grad F(\xilp)}^2.
\end{align}
For $\ell=0$,
\begin{align}
    \sqrt{\V{0}} &\le \sqrt{\Expc{\norm{\grad f(\xio, \bs{\theta})}^2}{\xio}} \\
    &= \sqrt{\Expc{\norm{\grad f(\xio, \bs{\theta}) - \grad f(\opt, \bs{\theta}) + \grad f(\opt, \bs{\theta})}^2}{\xio}} \\
    &\le
    \sqrt{\Expc{\norm{\grad f(\xio, \bs{\theta}) - \grad f(\opt, \bs{\theta})}^2}{\xio}}
    +
    \underbrace{
    \sqrt{\Exp{\norm{\grad f(\opt, \bs{\theta})}^2}}
     }_{\sqrt{V_*}} \\
    &\le
    \sigma \norm{\grad F(\xio)} + \sqrt{V_*}.
\end{align}

Let the total gradient sampling cost to reach iteration $k'+1$ be
\begin{align}
    \cl{C}_{k'} = \sum_{k=0}^{k'} \Delta \cl{C}_k(\cl{L}_k),
\end{align}
where $\Delta \cl{C}$ is defined as in \eqref{eq:cost_of_iter_k}.
In this section, we present limited analyses of \SGD[-MICE] to reach $k^*$ where we assume only the \Add operator is used, thus, using the equation above,
\begin{align}
    \cl{C}_{k^*-1} 
    &= \sum_{k=0}^{k^*-1} \sum_{\ell=0}^{k} \cost{\ell} (\M[k]{\ell} - \M[k-1]{\ell}) \\
    &= \sum_{\ell=0}^{k^*-1} \cost[k^*-1]{\ell} \M[k^*-1]{\ell}.
\end{align}
As will be shown in \S\ref{sec:numerics}, the other index set operators, \Drop[], \Restart[], and \Clip[] greatly improve the convergence of \SGD[-MICE].
As a consequence, these analyses considering only the \Add[] operator are pessimistic.

\subsubsection{Expectation minimization problems}
\label{sec:sampling_cost_expectation}
\begin{cor}[Expected gradient sampling cost of \texttt{SGD-MICE} with linear convergence]
    \label{cor:cost_sgd_mice}
    Let the Assumptions of Corollary~\ref{cor:grad_norm_convergence} and Assumption~\ref{as:variance} hold.
    Moreover, let $k^*$ be the smallest $k$ such that $\norm{\grad F(\bs{\xi}_{k})}^2 < tol$ and all sample sizes at the last iteration be larger than $M_{min}$.
    Then, the expected number of gradient evaluations needed to generate $\bs{\xi}_{k^*}$ is
    \begin{multline} \label{eq:expected_cost_expectation}
    \Exp{\cl{C}_{k^*-1}} \le \epsilon^{-2} tol^{-1}
    \left(
        4 \sigma
        \sqrt{L \Exp{F(\xio) - F(\opt)}}
        \left(\frac{2 \kappa}{1 - \epsilon^2}\right)
        + \sqrt{V_*}
    \right)^2 + \\
    2 M_{min}
    \left(
    \max\left\{0,
    \frac{\kappa}{1-\epsilon^2} \log(tol^{-1}2 L\Exp{F(\xio) - F(\opt)})\right\}
    + \frac{\kappa}{1-\epsilon^2}
    \right).
    \end{multline}
    Moreover, the relative gradient error that minimizes the expected gradient sampling cost is $\epsilon=\sqrt{1/3}$.
\end{cor}

\begin{proof}
We know that $k^*-1$ iterations are needed to generate $\bs{\xi}_{k^*}$.
Thus, the whole optimization cost is
\begin{align}
    \cl{C}_{k^*-1}
    &\le \epsilon^{-2} \norm{\grad F(\bs{\xi}_{k^*-1})}^{-2}
        \left(\sum_{\ell' \in \cl{L}_{k^*-1}} \sqrt{\V[k^*-1]{\ell'} \cost[k^*-1]{\ell'}}\right)^2
        + \sum_{\ell' \in \cl{L}_{k^*-1}} \cost[k^*-1]{\ell'} M_{min}
        \\
    &\le \epsilon^{-2} tol^{-1}
    \left(\sum_{\ell' \in \cl{L}_{k^*-1}} \sqrt{\V[k^*-1]{\ell'} \cost[k^*-1]{\ell'}}\right)^2
    + 2 |\cl{L}_{k^*-1}| M_{min}.
\end{align}
Let us analyze the following sum
\begin{align}
    \sum_{\ell' \in \cl{L}_{k^*-1}} \sqrt{\V[k^*-1]{\ell'} \cost[k^*-1]{\ell'}}
    &=
    \sqrt{\V{0}} + \sqrt{2} \sum_{1 \le \ell' \le k^*-1} \sqrt{\V{\ell'}} \\
    &\le
    \sigma \norm{\grad F(\xio)} + \sqrt{V_*}
    + \sqrt{2} \sigma \sum_{1 \le \ell' \le k^*-1}  \norm{\grad F(\bs{\xi}_{\ell'}) - \grad F(\bs{\xi}_{p_{k^*-1}(\ell')})}\\
    &\le
    \sigma \norm{\grad F(\xio)} + \sqrt{V_*} + 
    \sqrt{2} \sigma \sum_{1 \le \ell' \le k^*-1}  \norm{\grad F(\bs{\xi}_{\ell'})} + \norm{\grad F(\bs{\xi}_{p_{k^*-1}(\ell')})} \\    
    &\le
    2 \sqrt{2} \sigma \sum_{0 \le \ell' \le k^*-1}  \norm{\grad F(\bs{\xi}_{\ell'})} + \sqrt{V_*}.
\end{align}
Taking expectation of the summation above squared,
\begin{align}
    \Exp{\left(\sum_{\ell' \in \cl{L}_{k^*-1}} \sqrt{\V[k^*-1]{\ell'} \cost[k^*-1]{\ell'}}\right)^2}
    &\le
    \!\begin{multlined}[t][8cm]
    8 \, \sigma^2
    \sum_{\ell' \in \cl{L}_{k^*-1}} 
    \sum_{\ell \in \cl{L}_{k^*-1}}
    \Exp{
        \norm{\grad F(\bs{\xi}_{\ell'})}
        \norm{\grad F(\bs{\xi}_{\ell})}
    } \\
    + 4 \sigma \sqrt{2 V_*}
    \sum_{\ell'' \in \cl{L}_{k^*-1}}
    \Exp{\norm{\grad F(\bs{\xi}_{\ell''})}}
    + V_*
    \end{multlined}
    \\
    &\le
        \!\begin{multlined}[t][8cm]
            8 \, \sigma^2
            \sum_{\ell' \in \cl{L}_{k^*-1}} 
            \sum_{\ell \in \cl{L}_{k^*-1}} 
            \sqrt{
            \Exp{
                \norm{\grad F(\bs{\xi}_{\ell'})}^2
            }
            \Exp{
                \norm{\grad F(\bs{\xi}_{\ell})}^2
            }} \\
            + 4 \sigma \sqrt{2 V_*}
            \sum_{\ell'' \in \cl{L}_{k^*-1}}
            \sqrt{\Exp{\norm{\grad F(\bs{\xi}_{\ell''})}^2}}
            + V_*
        \end{multlined}
    \\
    &=
    \left(
        2 \sqrt{2} \sigma
        \sum_{\ell \in \cl{L}_{k^*-1}} 
        \sqrt{\Exp{\norm{\grad F(\bs{\xi}_{\ell})}^2}}
        + \sqrt{V_*}
    \right)^2 \\
    &\le
    \left(
        4 \sigma
        \sqrt{L \Exp{F(\xio) - F(\opt)}}
        \left(\sum_{\ell \in \cl{L}_{k^*-1}} r^{\ell/2} \right)
        + \sqrt{V_*} 
    \right)^2 \\
    &\le
    \left(
        4 \sigma
        \sqrt{L \Exp{F(\xio) - F(\opt)}}
        \left(\frac{1}{1 - \sqrt{r}}\right)
        + \sqrt{V_*} 
    \right)^2.
\end{align}
Substituting back to the expected cost,
\begin{align}
    \Exp{\cl{C}_{k^*-1}}
    \le \epsilon^{-2} tol^{-1}
    \left(
        4 \sigma
        \sqrt{L \Exp{F(\xio) - F(\opt)}}
        \left(\frac{1}{1 - \sqrt{r}}\right)
        + \sqrt{V_*} 
    \right)^2
    +
    2 \Exp{k^*} M_{min}.
\end{align}
Substituting the expected number of iterations from Corollary~\ref{cor:convergence.seq} and using \eqref{eq:rate_ineq_2} results in \eqref{eq:expected_cost_expectation}.

Since the term $(1 / (1 - \sqrt{r}))^2$ is $\bigO{(1 - \epsilon^2)^{-2}\kappa^2}$, it dominates convergence as $\kappa \rightarrow \infty$, thus the expected work of \texttt{SGD-MICE} without restart or dropping is $\bigO{\epsilon^{-2} (1 - \epsilon^2)^{-2} \kappa^2 tol^{-1}}$.
Therefore, the relative gradient error that minimizes the total gradient sampling cost is $\epsilon=\sqrt{1/3}$.
\end{proof}

  \begin{rmk}[Expected vs.\ high-probability sampling-cost bounds]\label{rmk:cost_expect_vs_hp}
      The results in \S\ref{sec:sampling_cost} bound the gradient sampling cost in \emph{expectation}, under the mean-square relative error control \eqref{eq:err_condition}.
  In contrast, the high-probability PL result in \S\ref{sec:hp_main} requires the stronger condition that the realized error satisfies $\|e_k\|\le \epsilon\|\nabla F(\bs{\xi}_k)\|$ \emph{uniformly in $k$} with probability at least $1-\delta$.
  Under Assumption~\ref{as:hp_coord_subg} and the conditional independence structure in Appendix~\ref{sec:hp_mice}, this uniform control is obtained by imposing the variance-sum constraint \eqref{eq:hp_variance_sum_constraint}, which inflates the required sample sizes by an explicit factor $\log(2d_\xi/\delta_k)$.
  With a standard summable schedule, for example $\delta_k \propto \delta/(k+1)^2$, one has $\log(2d_\xi/\delta_k)=\log(2d_\xi/\delta)+\bigO{\log(k+1)}$.
  Since $k^*(tol)=\bigO{\log(1/tol)}$ in the PL regime, the corresponding high-probability sampling-cost bounds preserve the leading $tol^{-1}$ scaling, at the cost of logarithmic dependence on $\delta$ (and at most a mild $\log\log(1/tol)$ correction through the schedule).
  Concretely, Appendix~\ref{sec:hp_mice} builds a single ``good'' event $\Omega^{\mathrm{MICE}}_\delta$ via a union bound over iterations (Lemma~\ref{lem:hp_mice_single_event}), shows that the relative error control holds uniformly on that event (Corollary~\ref{cor:hp_mice_rel_err}), and then derives explicit bounds on the iteration complexity (Lemma~\ref{lem:hp_iter_to_tol}) and on the total gradient sampling cost (Corollary~\ref{cor:hp_cost_to_tol}).
  \end{rmk}

\begin{cor}[Expected gradient sampling cost of \texttt{SGD-A}]
    \label{cor:cost_sgd_a}
    If Assumptions of Corollary~\ref{cor:grad_norm_convergence} hold and Assumption~\ref{as:variance} also holds, \texttt{SGD-A} generates an iterate $\bs{\xi}_{k^*}$ satisfying $\norm{\grad F(\bs{\xi}_{k^*})}^2 \le tol$ with an expected gradient sampling cost
    \begin{equation}
      \label{eq:expected_cost_sum}
      \begin{aligned}
        \Exp{\cl{C}_{k^*-1}} \le {} &
        \left(
            \frac{3 (\sigma^2+1)}{\epsilon^2}
            + \frac{2 V_*}{\epsilon^2 tol}
            + M_{min}
        \right)
        \\
        &\times
        \Bigl(
        \max\left\{0,
        \frac{\kappa}{1-\epsilon^2}\log\!\left(tol^{-1}2L\Exp{F(\xio) - F(\opt)}\right)\right\}
        + \frac{\kappa}{1-\epsilon^2}
        \Bigr).
      \end{aligned}
    \end{equation}
\end{cor}
\begin{proof}
    Let the gradient sampling cost of SGD-A be
    \begin{equation}
        \cl{C}_{k^*-1} = \sum_{k=0}^{k^*-1} \M{k}.
    \end{equation}
    The sample sizes are
    \begin{equation}
        \M{k} \le 
        \frac{\V{k}}{\epsilon^2 \norm{\grad F(\xik)}^2} + M_{min}
    \end{equation}
    We can bound $\V{k}$ as
    \begin{align}
        V_{\ell}
        & =
        \Expc{\norm{\grad f(\bs{\xi}_\ell, \bs{\theta})
          - \grad F(\bs{\xi}_\ell)}^2}{\bs{\xi}_\ell}\\[4pt]
        & \le
        2 \Expc{
        \norm{
        \grad f(\bs{\xi}_\ell, \bs{\theta})
        - \grad f(\bs{\xi}^*, \bs{\theta})
        }^2
        }{\bs{\xi}_\ell}
        + 2
        \Expc{
        \norm{
        \grad f(\bs{\xi}^*, \bs{\theta})
        - \grad F(\bs{\xi}^*)
        }^2
        }{\bs{\xi}_\ell}\\
        \label{eq:sgd_a_vl}
        & \le
        2 \sigma^2
        \norm{
         \grad F(\xil)
        }^2 + 2 V^*.
      \end{align}
      \begin{align}
        \cl{C}_{k^*-1}
        &\le
        \sum_{k=0}^{k^*-1} \left(
            \frac{\V{k}}{\epsilon^2 \norm{\grad F(\xik)}^2} + M_{min}
        \right) \\
        &\le
        \sum_{k=0}^{k^*-1} \left(
            \frac{2 \sigma^2}{\epsilon^2}
            + \frac{2 V^*}{\epsilon^2 \norm{\grad F(\xik)}^2}
            + M_{min}
        \right) \\
        &\le
        \sum_{k=0}^{k^*-1} \left(
            \frac{2 \sigma^2}{\epsilon^2}
            + \frac{2 V^*}{\epsilon^2 tol}
            + M_{min}
        \right) \\
        &\le
        \left(
            \frac{2 \sigma^2}{\epsilon^2}
            + \frac{2 V^*}{\epsilon^2 tol}
            + M_{min}
        \right)
        k^*.
      \end{align}
      Taking expectation and substituting $\Exp{k^*}$ from Corollary~\ref{cor:convergence.seq} finishes the proof.
\end{proof}

Table~\ref{tab:main_costs} presents, for different methods, the cost to reach a desired tolerance in expectation minimization.

\begin{table}[h]
  \caption{Summary of the main gradient sampling cost bounds in expectation minimization. The bounds depend on the relative error parameter $\epsilon$; we report the leading-order scaling in $\kappa$ and $tol$ (see Corollaries~\ref{cor:cost_sgd_mice} and~\ref{cor:cost_sgd_a} for full expressions).}
\label{tab:main_costs}
\centering
{\footnotesize
\setlength{\tabcolsep}{4.5pt}
\begin{tabularx}{\linewidth}{l X X X}
  \toprule
  Method & Gradient estimator / mechanism & Assumptions (expectation case) & Expected gradient sampling cost to reach $\norm{\grad F(\bs{\xi}_{k^*})}^2\le tol$ \\ \midrule
  \emph{vanilla} \SGD[] & Fixed-batch Monte Carlo gradients & Standard (see \cite{Rup88}) & Standard sublinear baseline \\
  \SGD[-A] & Adaptive batch sizes to enforce relative error & Cor.~\ref{cor:grad_norm_convergence} + Ass.~\ref{as:variance} & $\bigO{\kappa\,tol^{-1}\log(tol^{-1})}$ \\
  \emph{vanilla} \SGD[-MICE] & \MICE[] with \Add[] only (no operators) & Cor.~\ref{cor:grad_norm_convergence} + Ass.~\ref{as:variance} & $\bigO{\kappa^2\,tol^{-1}}$ \\
  \SGD[-MICE] + operators & \MICE[] with \Add[]/\Drop[]/\Restart[]/\Clip[] & Empirical improvement (see \S\ref{sec:numerics}) & $\bigO{\kappa\,tol^{-1}}$\footnote{Empirically observed.} \\
  \bottomrule
\end{tabularx}
}
\end{table}

\begin{rmk}[Stopping criterion]
  In practice, applying the stopping criterion \eqref{eq:stopping_criterion} requires an approximation of the mean gradient norm at each iteration.
  A natural approach is to use the  \texttt{MICE} estimator as such an approximation, yielding
  \begin{equation}\label{eq:stopping_real}
  \norm{\grad \cl{F}_{k^*}}^2 <   \, tol,
  \end{equation}
  provided that the error in the mean gradient is controlled in a relative sense.
  This quality assurance requires a certain number of gradient samples.
  For example, let us consider the ideal case of stopping when we start inside the stopping region, near the optimal point $\bs{\xi}^\ast$.
  To this end, suppose that the initial iteration point, $\bs{\xi}_0$, is such that $\norm{\grad F(\xio)}^2\le tol$. What is the cost needed to stop by sampling gradients at $\bs{\xi}_0$ without iterating at all? Observing that we need a tolerance $  \, tol$, we thus need a number of samples $M$ that  satisfies
  \begin{equation}
  \frac{\bb{E}\left[\norm{\nabla_{\bs{\xi}} f(\bs{\xi}_0,\theta)}^2\right]}{ \, tol} \le  M.
  \end{equation}
  Compare the last estimate with \eqref{eq:expected_cost_expectation} and \eqref{eq:expected_cost_sum}.
\end{rmk}

\subsubsection{Finite sum minimization problems}
\label{sec:sampling_cost_finite}

\begin{cor}[Cost analysis of \texttt{SGD-MICE} on the finite sum case]
If Assumptions of Corollary~\ref{prp:pl_convergence} hold, \texttt{SGD-MICE} achieves a stopping criterion with expected gradient sampling cost
\begin{align}
    \Expc{\cl{C}_{k^*-1}}{\xio}
    &\le
    \!\begin{multlined}[t][10.5cm]
    \frac{
        (N-1)
        \left(
            8
            \frac{\kappa}{1-\epsilon^2}
            \sigma
            \sqrt{L (F(\xio) - F(\opt))}
            + \sqrt{V_*} 
        \right)^2
    }{
        \V[k^*-1]{0}
    }
    \log\left(
        \frac{
            \V[k^*-1]{0}
        }{
            tol (N-1) \epsilon^2
        }
        +1
    \right) \\
    +M_{min}
    \left(
        \max\left\{0,
        \frac{\kappa}{1-\epsilon^2} \log(tol^{-1}2 L
        (F(\xio) - F(\opt)))\right\}
        + \frac{\kappa}{1-\epsilon^2}
    \right)
    \end{multlined}
\end{align}
\end{cor}

\begin{proof}
\begin{align}
    \cl{C}_{k^*-1} 
    &= \sum_{\ell=0}^{k^*-1} \cost[k^*-1]{\ell} \M[k^*-1]{\ell} \\
    &\le 
    \frac{N}{N-1}
    \frac{
        \left(\sum_{\ell=0}^k \sqrt{\cost[k^*-1]{\ell}\V[k^*-1]{\ell}}\right)^2
    }{
        \epsilon^2 \gnorm{k^*-1}^2 + (N-1)^{-1} \sum_{\ell'=0}^k \V[k^*-1]{\ell'}
    }
    + (k^*-1) M_{min} \\
    &\le
    2
    \frac{
        \left(\sum_{\ell=0}^k \sqrt{\cost[k^*-1]{\ell}\V[k^*-1]{\ell}}\right)^2
    }{
        \epsilon^2 tol + (N-1)^{-1} \V[k^*-1]{0}
    }
    + k^* M_{min}
\end{align}
Taking expectation conditioned on the initial iterate,
\begin{align}
    \Expc{\cl{C}_{k^*-1}}{\xio}
    \le
    \frac{
        \left(
            4 \sigma
            \sqrt{L (F(\xio) - F(\opt))}
            \left(\frac{1}{1 - \sqrt{r}}\right)
            + \sqrt{V_*} 
        \right)^2
    }{
        \epsilon^2 tol + (N-1)^{-1} \V[k^*-1]{0}
    }
    + \Expc{k^*}{\xio} M_{min}.
\end{align}
Using the following logarithm inequality  with $c/b+1>0$,
\begin{align}
    \frac{a}{b + c} 
    \le \frac{a}{c} \log\left(\frac{c}{b} + 1\right),
\end{align}
gives
\begin{align}
    \Expc{\cl{C}_{k^*-1}}{\xio}
    &\le
    \!\begin{multlined}[t][10.5cm]
    \frac{
        (N-1)
        \left(
            4 \sigma
            \sqrt{L (F(\xio) - F(\opt))}
            \left(\frac{1}{1 - \sqrt{r}}\right)
            + \sqrt{V_*} 
        \right)^2
    }{
        \V[k^*-1]{0}
    }
    \log\left(
        \frac{
            \V[k^*-1]{0}
        }{
            tol (N-1) \epsilon^2
        }
        +1
    \right) \\
    + \Expc{k^*}{\xio} M_{min}.
    \end{multlined}
\end{align}
Using Corollary~\ref{cor:convergence.seq} and \eqref{eq:rate_ineq_2} concludes the proof.
\end{proof}

\begin{cor}[Cost analysis of \texttt{SGD-A} on the finite sum case]
If the assumptions of Proposition~\ref{prp:pl_convergence} are satisfied, \texttt{SGD-A} finds an iterate $\bs{\xi}_{k^*}$ such that $\norm{\grad F(\bs{\xi}_{k^*})}^2 \le tol$ with expected gradient sampling cost
\begin{equation}
  \begin{aligned}
    \Exp{\cl{C}_{k^*-1}}
    \le {} &
    N
    \min\left\{
    1,
    \log\left(
        \frac{\frac{2 V^*}{tol} + 2 \sigma^2}{\epsilon^2(N-1)}
        +1
        \right)
    \right\}
    \\
    & \times
    \left(
    \max\left\{0,
    \frac{\kappa}{1-\epsilon^2}\log\!\left(tol^{-1}2L\Exp{F(\xio) - F(\opt)}\right)\right\}
    + \frac{\kappa}{1-\epsilon^2}
    \right).
  \end{aligned}
\end{equation}
\end{cor}

\begin{proof}
When using \texttt{SGD-A} to solve the finite sum minimization problem while taking into consideration that variance goes to zero as $M \rightarrow N$, the sample size at iteration $k$ is
\begin{equation}
    M_k = \left\lceil
    \frac{N}{N-1}
    \frac{V_k}{\epsilon^2 \norm{\grad F(\xik)}^2
    + \frac{V_k}{N-1}
    }
    \right\rceil,
\end{equation}
where $V_k = \Expc{\norm{\grad f(\xik, \tht) - \grad F(\xik)}^2}{\xik}$.
Thus, the total gradient sampling cost to reach iteration $k^*$ is
\begin{equation}
    \cl{C}_{k^*-1}
    \le \sum_{\ell=0}^{k^*-1}
    \frac{N}{N-1}
    \frac{V_\ell}{\epsilon^2 \norm{\grad F(\xil)}^2
    + \frac{V_\ell}{N-1}
    }.
\end{equation}
Using \eqref{eq:sgd_a_vl},
\begin{align}
    \cl{C}_{k^*-1}
    &\le \sum_{\ell=0}^{k^*-1}
    N
    \frac{
    2\sigma^2 \norm{\grad F(\xil)}^2 + 2 V^*
    }{
        \epsilon^2 \norm{\grad F(\xil)}^2 (N-1)
    + 2\sigma^2 \norm{\grad F(\xil)}^2 + 2 V^*
    } \\
    \label{eq:sgd_a_finite_1}
    &\le     N
    \frac{
    2\sigma^2 tol + 2 V^*
    }{
        tol(\epsilon^2 (N-1)
    + 2\sigma^2) + 2 V^*
    } k^*\\    
    &\le N k^*.
\end{align}
Another bound can be obtained from \eqref{eq:sgd_a_finite_1} as
\begin{align}
    \cl{C}_{k^*-1}
    &\le 
    N
    \frac{
    1
    }{
        \frac{
        \epsilon^2 tol (N-1)
        }{
        2\sigma^2 tol + 2 V^*
        }
    + 1
    } k^* \\
    &\le
    N
    \log\left(
    \frac{
        2\sigma^2 + 2 \frac{V^*}{tol}
        }{
        \epsilon^2 (N-1)
    }
    +1
    \right) k^*.
\end{align}
Taking expectation and using \eqref{eq:expect} concludes the proof.
\end{proof}

\begin{rmk}[More general $\bs{\theta}$ probability distributions]\label{rmk:theta_distribution}
  Although in Assumption~\ref{as:theta_distribution} we restricted our attention to the case where the probability distribution of $\bs{\theta}$, $\pi$, does not depend on $\bs{\xi}$, it is possible to use mappings to address more general cases. Indeed, let us consider the case where
  \begin{equation}
  \bs{\theta} = h(\tilde{\bs{\theta}}, \bs{\xi}),
  \end{equation}
  for some given smooth function $h$ and such that the distribution of $\tilde{\bs{\theta}}$, $\tilde \pi$, does not depend on $\bs{\xi}$.
  Then we can simply write, letting $\tilde f(\bs{\xi},\tilde{\bs{\theta}}) = f(\bs{\xi},h(\tilde{\bs{\theta}},\bs{\xi}))$,
  \begin{equation}
  F(\bs{\xi}) = \bb{E}\left[f(\bs{\xi},\bs{\theta})|\bs{\xi}\right] = \bb{E}\left[\tilde f(\bs{\xi},\tilde{\bs{\theta}})|\bs{\xi}\right]
  \end{equation}
  and, by sampling $\tilde{\bs{\theta}}$ instead of $\bs{\theta}$, we are back in the setup of Assumption~\ref{as:theta_distribution}.
  \end{rmk}

Having established theoretical convergence guarantees and gradient sampling cost bounds for \texttt{SGD-MICE}, we now turn to practical implementation aspects of the algorithm.

\section{\texttt{MICE} algorithm}\label{sec:mice_algorithm}

In this section, we describe the \texttt{MICE} algorithm and some of its practical implementation aspects.
Before we start, let us discuss the resampling technique used to build an approximated probability distribution for the norm of the gradient.

\begin{rmk}[Gradient resampling for calculating sample sizes]\label{rmk:grad_resampling}
To approximate the empirical distribution of $\norm{\Fk}$, we perform a jackknife \cite{efron1982jackknife} resampling of the approximate mean gradient using sample subsets for each iteration $\ell \in \cl{L}_k$.

First, for each element $\ell \in \cl{L}_k$, we partition the index set $\cl{I}_{\ell, k}$ in $n_\text{part}$ disjoint sets $\cl{I}^{(1)}_{\ell, k}, \cl{I}^{(2)}_{\ell, k},.., \cl{I}^{(n_\text{part})}_{\ell, k}$ with the same cardinality.
Then, we create, for each of these sets, their complement with respect to $\cl{I}_{\ell, k}$, i.e., $ \overline{\cl{I}}^{(i)}_{\ell, k} = \cl{I}_{\ell, k} \setminus \cl{I}^{(i)}_{\ell, k}$ for all $i=1,2,..,n_\text{part}$.
We use these complements to compute the average of these gradient differences excluding a portion of the data,
\begin{equation}
    \overline{\mu}_{\ell, k}^{(i)}
    =
    \left| \overline{\cl{I}}^{(i)}_{\ell,k} \right|^{-1} \sum_{\alpha\in\overline{\cl{I}}^{(i)}_{\ell,k}}\Delta_{\ell,k,\alpha},
\end{equation}
which we then sample for each $\ell \in \cl{L}_k$ to get a single sample of the mean gradient,
\begin{equation}\label{eq:resampling}
\nabla_{\bs{\xi}}\cl{F}_{k, \nu}
\coloneqq
\sum_{\ell\in\cl{L}_k}
\overline{\mu}_{\ell, k}^{\left(i_{\ell, \nu}\right)}
\end{equation}
by independently sampling $i_{\ell, \nu}$ from a categorical distribution with $n_\text{part}$ categories.
Sampling $\nabla_{\bs{\xi}}\cl{F}_{k, \nu}$ $n_{\text{samp}}$ times, we construct a set of gradient mean estimates $\{\nabla_{\bs{\xi}}\cl{F}_{k, \nu}\}_{\nu=1}^{n_\text{samp}}$.

Similarly, we set a right tail quantile $1-p_{\text{stop}}$ with $p_{\text{stop}} \le 0.5$ to define a gradient norm to be used as a stopping criterion.
We stop at $k$ if
\begin{equation}
    \norm{\nabla_{\bs{\xi}} \cl{F}_k^\text{stop}}^2 \le tol,
\end{equation}
where $\norm{\nabla_{\bs{\xi}} \cl{F}_k^\text{stop}}$ is the norm of the gradient respective to the $1 - p_{\text{stop}}$ quantile.

The resampling procedure requires storing $n_\text{part}$ partition-complement means per level, but it does not require additional gradient evaluations.

\begin{algorithm}[h]
\setstretch{1.35}
\caption{Gradient resampling (partition jackknife) for sizing and stopping.}
\label{alg:grad_resampling}
\begin{algorithmic}[1]
\Require Index set $\cl{L}_k$, level samples $\{\Delta_{\ell,k,\alpha}\}_{\ell\in\cl{L}_k,\alpha\in\cl{I}_{\ell,k}}$, integers $n_\text{part}$ and $n_\text{samp}$, quantiles $p_\text{re}\le 0.5$ and $p_\text{stop}\le 0.5$.
\For{$\ell \in \cl{L}_k$}
  \State Partition $\cl{I}_{\ell,k}$ into $\cl{I}^{(1)}_{\ell,k},\dots,\cl{I}^{(n_\text{part})}_{\ell,k}$ and define complements $\overline{\cl{I}}^{(i)}_{\ell,k}=\cl{I}_{\ell,k}\setminus \cl{I}^{(i)}_{\ell,k}$.
  \State Compute $\overline{\mu}_{\ell,k}^{(i)}=\left|\overline{\cl{I}}^{(i)}_{\ell,k}\right|^{-1}\sum_{\alpha\in\overline{\cl{I}}^{(i)}_{\ell,k}}\Delta_{\ell,k,\alpha}$ for all $i=1,\dots,n_\text{part}$.
\EndFor
\For{$\nu = 1,\dots,n_\text{samp}$}
  \State Sample independently $i_{\ell,\nu}\in\{1,\dots,n_\text{part}\}$ for each $\ell\in\cl{L}_k$.
  \State Form $\nabla_{\bs{\xi}}\cl{F}_{k,\nu}=\sum_{\ell\in\cl{L}_k}\overline{\mu}_{\ell,k}^{(i_{\ell,\nu})}$ and store $\norm{\nabla_{\bs{\xi}}\cl{F}_{k,\nu}}$.
\EndFor
\State Set $\norm{\nabla_{\bs{\xi}}\cl{F}_{k}^{\text{re}}}$ as the $p_\text{re}$-quantile of $\{\norm{\nabla_{\bs{\xi}}\cl{F}_{k,\nu}}\}_{\nu=1}^{n_\text{samp}}$.
\State Set $\norm{\nabla_{\bs{\xi}}\cl{F}_{k}^{\text{stop}}}$ as the $(1-p_\text{stop})$-quantile of $\{\norm{\nabla_{\bs{\xi}}\cl{F}_{k,\nu}}\}_{\nu=1}^{n_\text{samp}}$.
\State \textbf{Stop} if $\norm{\nabla_{\bs{\xi}}\cl{F}_{k}^{\text{stop}}}^2 \le tol$.
\end{algorithmic}
\end{algorithm}

The resampling step performs $n_\text{samp}$ aggregations of $|\cl{L}_k|$ vectors (plus norm evaluations), i.e., $\bigO{n_\text{samp}|\cl{L}_k|d_{\bs{\xi}}}$ arithmetic operations, with memory overhead proportional to $n_\text{part}|\cl{L}_k|d_{\bs{\xi}}$.
On the other hand, the confidence on the empirical distribution built by the resampling technique depends on the resampling sample size $n_\text{samp}$.
Conditional on the partition-complement means $\{\overline{\mu}_{\ell,k}^{(i)}\}$, the resampled norms $\{\norm{\nabla_{\bs{\xi}}\cl{F}_{k,\nu}}\}_{\nu=1}^{n_\text{samp}}$ in Algorithm~\ref{alg:grad_resampling} are i.i.d.\ draws (since the indices $i_{\ell,\nu}$ are sampled independently across $\nu$). Therefore, for the empirical CDF $\widehat{F}_{n_\text{samp}}$ of $\norm{\nabla_{\bs{\xi}}\cl{F}_{k,\nu}}$, the Dvoretzky--Kiefer--Wolfowitz inequality \cite{DKW1956,Massart1990} yields
\begin{equation}
\Probc{\sup_t\left|\widehat{F}_{n_\text{samp}}(t)-F(t)\right|\le \varepsilon}{\{\overline{\mu}_{\ell,k}^{(i)}\}}
\ge 1-2\exp\!\left(-2n_\text{samp}\varepsilon^2\right).
\end{equation}
In particular, with probability at least $1-\delta$ (conditionally), choosing $n_\text{samp}\ge \frac{1}{2\varepsilon^2}\log\!\frac{2}{\delta}$ guarantees that the empirical $(1-p_\text{stop})$ quantile used for stopping has tail miscoverage at most $p_\text{stop}+\varepsilon$ under the resampling distribution.

To control the work of the resampling technique, we measure the runtime needed to get a sample of $\nabla_{\bs{\xi}}\cl{F}_{k, \nu}$ and then set $n_\text{samp}$ so that the overall time does not exceed a fraction $\delta_\text{re}$ of the remaining runtime of \MICE[].
From our numerical tests, we recommend $n_\text{part}$ to be set between $3$ and $10$, $\delta_\text{re}$ between $0.1$ (for expensive gradients) and $1$, $n_{\text{samp}} \ge 10$, and $p_\text{re}=5 \%$.

\end{rmk}

In Algorithm~\ref{alg:mice}, we present the pseudocode for the \texttt{MICE} estimator, and in Algorithm~\ref{alg:indexset} we present the algorithm to update the index set $\cl{L}_k$ from $\cl{L}_{k-1}$ according to \S\ref{sec:opt_index_set}.
Two coupling algorithms for the multi-iteration stochastic optimizers are presented in Appendix~\ref{sec:algorithms}: \texttt{SGD-MICE} and \texttt{Adam-MICE}.

\begin{algorithm}
  \setstretch{1.5}
  \caption{Index-set update operator selection.}
  \label{alg:indexset}
  \begin{algorithmic}[1]
  \Procedure{\texttt{Index Set}}{$\cl{L}_{k-1}$, $\V{k}$, $\V{k}^{\text{drop}}$}
  \State $\cl{L}^{\text{add}}_{k}    \gets \cl{L}_{k-1}\cup\{k\}$
  \State $\cl{L}^{\text{drop}}_{k}   \gets \cl{L}_{k-1}\cup\{k\}\!\setminus\{k-1\}$
  \State $\cl{L}^{\text{rest}}_{k}   \gets \{k\}$
  \State Set $\cl{L}^{\text{clip}, \ell^*}_{k}$ as in \S\ref{sec:mice.clipping} with \Clip ``A'' or ``B''
  \State Set $\cl{L}_k \gets \cl{L}^{\text{add}}_{k}$
  \State Compute $\Delta \cl{W}_k(\cl{L}_k)$ and $\Delta \cl{W}_k(\cl{L}^{\text{drop}}_{k})$ using $\V{k}$, $\V{k}^{\text{drop}}$ and \eqref{eq:work_update_MICE}
  \If{$\Delta \cl{W}_k(\cl{L}^{\text{drop}}_{k}) \le (1+\delta_{\text{drop}})\Delta \cl{W}_k(\cl{L}_k)$}
    \State $\cl{L}_k \gets \cl{L}^{\text{drop}}_{k}$
  \EndIf
  \If{$\Delta \cl{W}_k(\cl{L}^{\text{rest}}_{k}) < (1+\delta_{\text{rest}})\Delta \cl{W}_k(\cl{L}_k)$}
    \State $\cl{L}_k \gets \cl{L}^{\text{rest}}_{k}$
  \EndIf
  \If{$\min_{\ell \in \cl{L}_{k-1}}\Delta \cl{W}_k(\cl{L}^{\text{clip,}\ell}_{k}) < \Delta \cl{W}_k(\cl{L}_k)$}
    \State Clip $\cl{L}_k$ at the minimizing $\ell$
  \EndIf
  \State \textbf{return} $\cl{L}_k$
  \EndProcedure
\end{algorithmic}
\end{algorithm}

\begin{algorithm}
\setstretch{1.5}
\caption{The \texttt{MICE} estimator.}
\label{alg:mice}
\begin{algorithmic}[1]
\Procedure{\texttt{MICE}}{}
\State $\cl{I}_k \gets \{\alpha\}_{\alpha=1}^{M_{min}}$
\State Sample $\bs{\theta}_\alpha \sim \pi \qquad\forall \alpha \in \cl{I}_k$
\State Compute $\grad f(\xik, \bs{\theta}_\alpha)$
\If{$k=0$}
\State Set $\cl{L}_k \gets \{0\}$ \Comment No differences at $k=0$
\State Compute $\V{k}=\Expc{\norm{\grad f(\xik, \bs{\theta}_\alpha) - \grad F(\xik)}^2}{\xik}$
\Else
\State Compute $\grad f(\bs{\xi}_{k-1}, \bs{\theta}_\alpha)$ and $\grad f(\bs{\xi}_{p_{k}(k-1)}, \bs{\theta}_\alpha)$
\State Compute
$
\begin{cases}
	    \V{k}=\Expc{\norm{\grad f(\xik, \bs{\theta}_\alpha) - \grad f(\bs{\xi}_{k-1}, \bs{\theta}_\alpha)}^2}{\xik, \bs{\xi}_{k-1}} \\
	    \V{k}^{\text{drop}} = \Expc{\norm{\grad f(\xik, \bs{\theta}_\alpha) - \grad f(\bs{\xi}_{p_{k}(k-1)}, \bs{\theta}_\alpha)}^2}{\xik, \bs{\xi}_{p_{k}(k-1)}}
\end{cases}
$
\State Use Algorithm~\ref{alg:indexset} to set $\cl{L}_k$
\EndIf
\While{$\sum_{\ell \in \cl{L}_k} \frac{V_{\ell, k}}{M_{\ell, k}} - \frac{V_{\ell, k}}{N}
    \ge 
    \epsilon^2 \norm{\grad F(\xik)}^2$}
    \Comment For expectation minimization, $\frac{V_{\ell, k}}{N}=0$.
\State Calculate $\left\{M^*_{\ell, k}\right\}_{\ell \in \cl{L}_k}$ from \eqref{eq:optimal.sampling.1} or Algorithm~\ref{alg:finite_sum} using 
    $\norm{\grad F(\xik)} \approx \norm{\nabla_{\bs{\xi}}\cl{F}_k^{\text{re}}}$
\For{$\ell \in \cl{L}_k$}
\State $\Delta M_{\ell, k} = \min\{(M^*_{\ell, k} - M_{\ell, k})_+, 2 \M{\ell}, N-M_{\ell, k}\}$
\Comment $(a)_+ \coloneqq \max\{a,0\}$
\State $\cl{I}'_\ell \gets \{\alpha\}_{\alpha=M_{\ell, k}+1}^{M_{\ell, k} + \Delta M_{\ell, k}}$
\State Sample $\bs{\theta}_\alpha \sim \pi \qquad\forall \alpha \in \cl{I}'_\ell$
\State Obtain $\Delta_{\ell, k, \alpha}$ from \eqref{eq:Delta.sampled} for each $\alpha \in \cl{I}'_\ell$
\State Calculate $V_{\ell,k}$ from \eqref{eq:Vl}
\State Get $\nabla_{\bs{\xi}}\cl{F}_k$ using \eqref{eq:mice}
\EndFor
\State $M_{\ell, k} \gets M^*_{\ell, k}$
\EndWhile
\State \textbf{return} $\nabla_{\bs{\xi}}\cl{F}_k = \sum_{\ell\in\cl{L}_k}\dfrac{1}{M_{\ell,k}^*}
\sum_{\alpha\in\cl{I}_{\ell,k}}\Delta_{\ell,k,\alpha}$ from \eqref{eq:mice}
\EndProcedure
\end{algorithmic}
\end{algorithm}

In general, keeping all gradient realizations for all iterations in memory may be computationally inefficient, especially for large-dimensional problems.
To avoid this unnecessary memory overhead, we use Welford's online algorithm to estimate the variances $V_{\ell,k}$ online. We keep in memory only the samples mean and second-centered moments and update them in an online fashion \cite{welford1962note}.
This procedure makes the memory overhead much smaller than naively storing all gradients and evaluating variances when needed.
Therefore, for each $\ell \in \cl{L}_k$ at iteration $k$, we need to store the mean gradient differences estimate, a vector of size $d_{\bs{\xi}}$; $V_{\ell, k}$, a scalar; and $M_{\ell, k}$, an integer.
Also, we store the gradient mean estimate in case we might clip the index set at $\ell$ in the future, and the respective sum of the variances component-wise, also using Welford's algorithm.
Thus, for first-order methods such as \texttt{Adam-MICE} and \texttt{SGD-MICE}, the memory overhead of \texttt{MICE} is of $2 | \cl{L}_k| (d_{\bs{\xi}}+2)$ floating-point numbers and $| \cl{L}_k|$ integers.
Thus, for large-scale problems, dropping iterations and restarting the index set are very important to reduce memory allocation.
Regarding the computational overhead, updating each $V_{\ell, k}$ using Welford's algorithm at iteration $k$ has complexity $\bigO{(M_{\ell, k} - M_{\ell, k-1}) \, d_{\bs{\xi}}}$.
Computing the sample sizes using \eqref{eq:optimal.sampling.1} or Algorithm~\ref{alg:finite_sum} requires a number of operations that is $\bigO{|\cl{L}_k| \, d_{\bs{\xi}}}$.
While sample sizes might be computed several times per iteration due to the progressive sample size increase, this cost does not increase with the dimensionality of the problem.
The resampling technique presented in \eqref{eq:resampling} increases the memory overhead by a factor $n_{\text{part}}$ and the computational work by a factor $\delta_{re}$.

\section{Numerical examples}\label{sec:numerics}

We validate \texttt{MICE} through three experiments of increasing complexity: (i) synthetic quadratic functions to study scaling with condition number and ablate operator contributions, (ii) a stochastic Rosenbrock problem to demonstrate coupling with \texttt{Adam}, and (iii) large-scale logistic regression comparing against state-of-the-art variance reduction baselines.

When using \texttt{SGD}, with or without \texttt{MICE}, we assume the constant $L$ to be known and use it to compute the step-size $\eta = 1/L$.
As a measure of the performance of the algorithms, we use the optimality gap, which is the difference between the approximate optimal value at iteration $k$ and the exact optimal value,
\begin{equation}\label{eq:opt_gap}
F(\bs{\xi}_k)-F(\bs{\xi}^*).
\end{equation}
In some examples, we know the optimal value and optimal point analytically; otherwise, we estimate numerically by letting optimization algorithms run for many iterations.

As for \texttt{MICE} parameters, when coupled with \SGD[], we use $\epsilon=\sqrt{1/3}$, and when coupled with \Adam[] we use $\epsilon=1$.
The other parameters are fixed for all problems, showing the robustness of \texttt{MICE} with respect to the tuning: $\delta_{\mathrm{drop}}=0.5$, $\delta_{\mathrm{rest}}=1.0$, $M_{min}$ is set to $5$ for general iterations and $500$ for restarts, and the maximum index set cardinality is set to $100$ except when noted (in the sensitivity studies below and in the logistic regression over the dataset \emph{HIGGS}.).
For the continuous cases, we use the clipping ``A'', whereas, for the finite case, we use clipping ``B''.
As for the resampling parameters, we use $n_\text{part}=5$, $\delta_\text{re}=1.0$, and $p_\text{re}=0.05$ with a minimum resampling size of $10$.
Note, however, that the current \MICE implementation is not optimized for performance, and could be much improved in this sense.
Regarding the stopping criterion, except in the first example, we do not define a $tol$. 
Instead, we define a fixed gradient sampling cost that, when reached, halts execution.
This choice allows us to better compare \SGD[-]\MICE with other methods.

  The Python implementation of \texttt{MICE} and baseline methods used to generate the data and figures presented in this section is available at GitHub\footnote{\url{https://github.com/agcarlon/mice}}.
  Moreover, \texttt{MICE} can be installed using PyPI; see the package documentation\footnote{\url{https://mice.readthedocs.io}}
  for more information.

\subsection{Random quadratic function}\label{ex:1}

\FloatBarrier
This problem is a simple numerical example devised to test the performance of \texttt{SGD-MICE} on the minimization of a strongly convex function.
The function whose expected value we want to minimize is
\begin{equation}
  f(\bs{\xi}, \theta) = \frac{1}{2} \bs{\xi} \cdot \bs{H}(\theta) \, \bs{\xi}
  - \bs{b} \cdot \bs{\xi},
\end{equation}
where
\begin{equation}
 \bs{H}(\theta) \coloneqq
 \bs{I}_2(1-\theta) +
 \begin{bmatrix}
  2 \kappa & 0.5 \\
  0.5      & 1
 \end{bmatrix}
 \theta,
\end{equation}
$\bs{I}_2$ is the identity matrix of size $2$, $\bs{b}$ is a vector of ones, and $\theta \sim \cl{U}(0, 1)$.
We use $\kappa=100$ and initial guess
$\bs{\xi}_0=(20, 50)$.
The objective function to be minimized is
\begin{equation}\label{eq:quad.example}
  F(\bs{\xi}) = \frac{1}{2} \bs{\xi} \cdot
  \bb{E}\left[\bs{H}(\theta)\right] \, \bs{\xi}
  - \bs{b} \cdot \bs{\xi},
\end{equation}
where
\begin{equation}
 \bb{E}\left[\bs{H}(\theta)\right] =
 \begin{bmatrix}
  \kappa+0.5 & 0.25 \\
  0.25       & 1
 \end{bmatrix}.
\end{equation}
The optimal point of this problem is $\bs{\xi}^* = \bb{E}\left[\bs{H}(\theta)\right]^{-1} \bs{b}$.
To perform optimization using \texttt{SGD-MICE} and \texttt{SGD}, we use the unbiased gradient estimator
\begin{equation}
  \nabla_{\bs{\xi}} f(\bs{\xi}, \theta) = \bs{H}(\theta) \, \bs{\xi} - \bs{b}.
\end{equation}
We use the eigenvalues of the Hessian of the objective function, $ \bb{E}\left[\bs{H}(\theta)\right]$, to calculate $L$ and thus define the step-size as $1/L$.
We set a stopping criterion of $tol=10^{-8}$.

\begin{figure}
  \includegraphics[width=.9\linewidth]{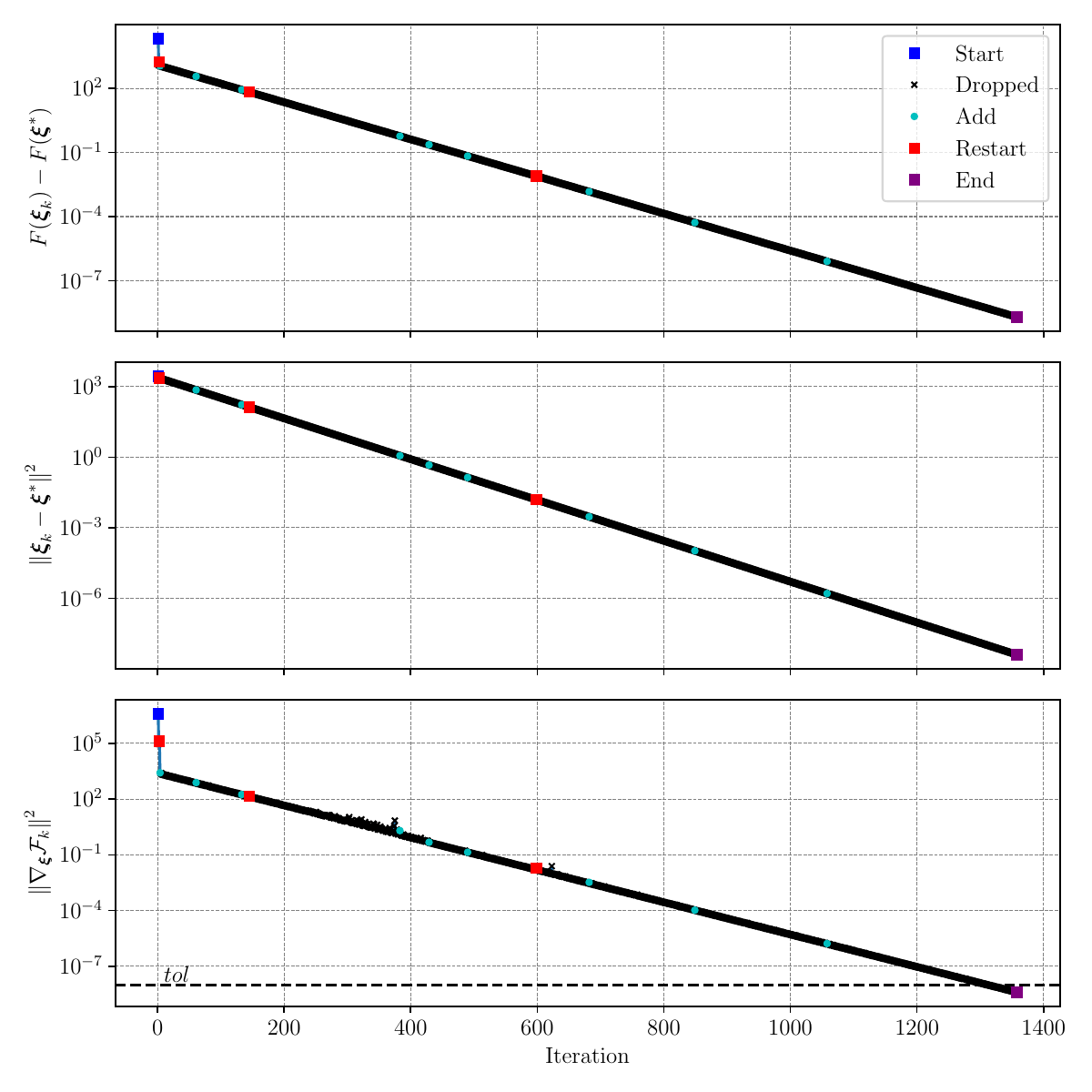}
  \caption{
  Single run, random quadratic example, Equation \eqref{eq:quad.example} with $\kappa=100$.
  Optimality gap (top), squared distance to the optimal point (center), and squared norm of gradient estimate (bottom) per iteration for \texttt{SGD-MICE}.
  The starting point, the restarts, and the end are marked respectively as blue, red, and purple squares, iterations dropped with black $\times$, and the remaining \texttt{MICE} points with cyan circles.
  \texttt{SGD-MICE} is able to achieve linear $L^2$ convergence as predicted in Proposition~\ref{prp:pl_convergence}.
  }
  \label{fig:quad_per_k}
\end{figure}

\begin{figure}
  \includegraphics[width=.9\linewidth]{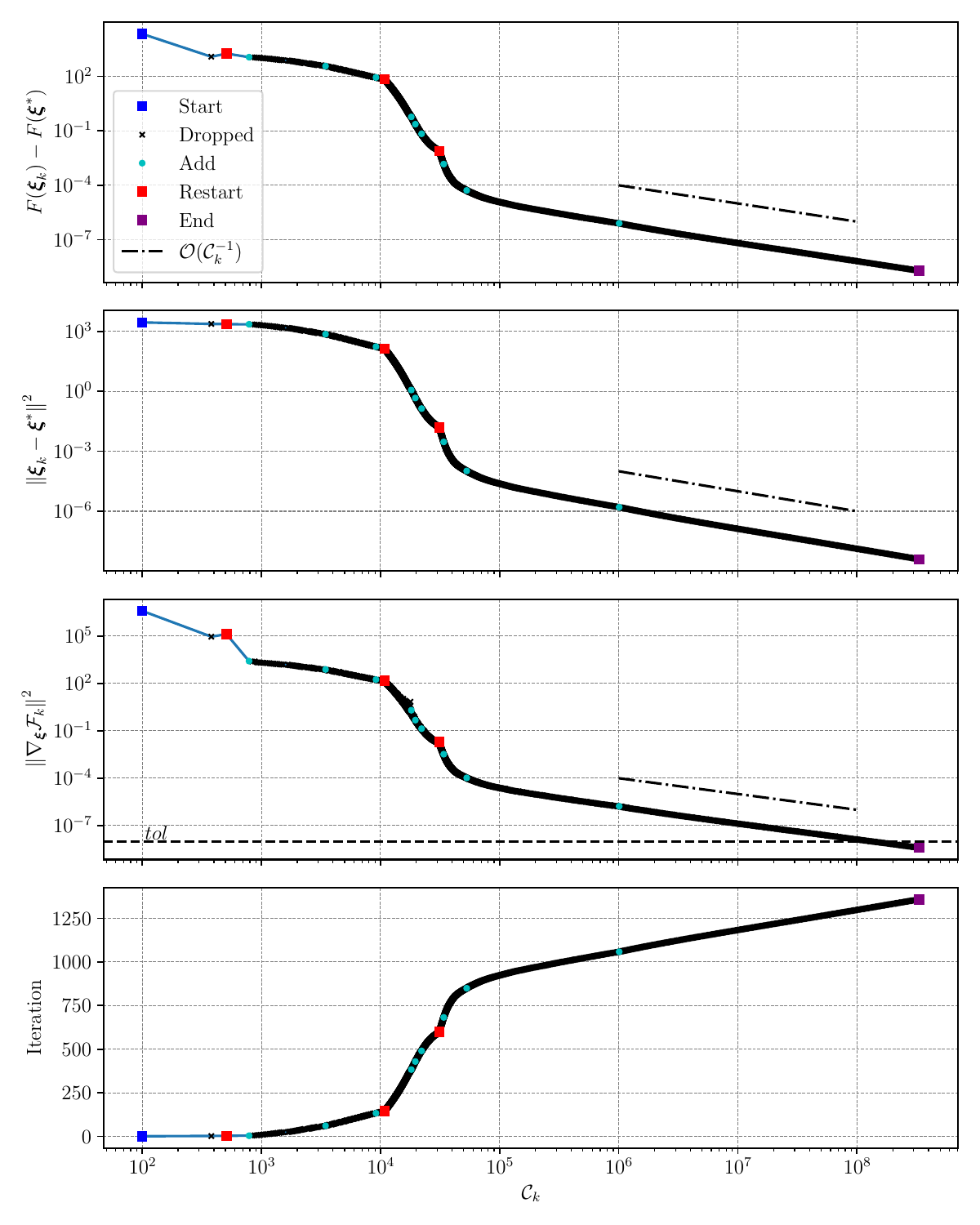}
  \caption{Single run, random quadratic example, Equation \eqref{eq:quad.example} with $\kappa=100$.
  Optimality gap (top), squared distance to the optimal point (center top), squared norm of gradient estimate (center bottom), and number of iterations (bottom) per number of gradient evaluations for \texttt{SGD-MICE}.
  The starting point, the restarts, and the end are marked respectively as blue, red, and purple squares, iterations dropped with black $\times$, and the remaining \texttt{MICE} points with cyan circles.
  The asymptotic convergence rate of $\cl{O}(\cl{C}_k^{-1})$ is presented when expected.
  }
  \label{fig:quad_per_grad}
\end{figure}

In Figures~\ref{fig:quad_per_k} and \ref{fig:quad_per_grad}, we present the optimality gap \eqref{eq:opt_gap}, the squared distance to the optimal point and the squared norm of the gradient estimate versus iteration and number of gradient sampling cost, respectively.
In Figure~\ref{fig:quad_per_grad} we also plot the iteration reached versus gradient sampling cost.
We mark the starting points, restarts, and ending points with blue, red, and purple squares, respectively; the dropped points with black $\times$, and the remaining iterations in the \texttt{MICE} index set with cyan dots.
In Figure~\ref{fig:quad_per_k}, one can observe that \texttt{SGD-MICE} attains linear convergence with a constant step-size, as predicted in Proposition~\ref{prp:pl_convergence}.
In Figure~\ref{fig:quad_per_grad}, we present the convergence plots versus gradient sampling cost, exhibiting numerical rates of $\cl{O}(\cl{C}_k^{-1})$.
These rates are expected as the distance to the optimal point converges linearly (see Corollary~\ref{cor:dist_to_opt_convergence}) and the cost of sampling new gradients per iteration grows as $\cl{C}_k = \cl{O}(\norm{\nabla F(\bs{\xi}_k)}^{-2})$, as shown in \eqref{eq:optimal.sampling.1}.
Finally, \texttt{SGD-MICE} was able to automatically decide whether to drop iterations, restart, or clip the index set to minimize the overall work required to attain the linear convergence per iteration.

\begin{figure}
  \includegraphics[width=.9\linewidth]{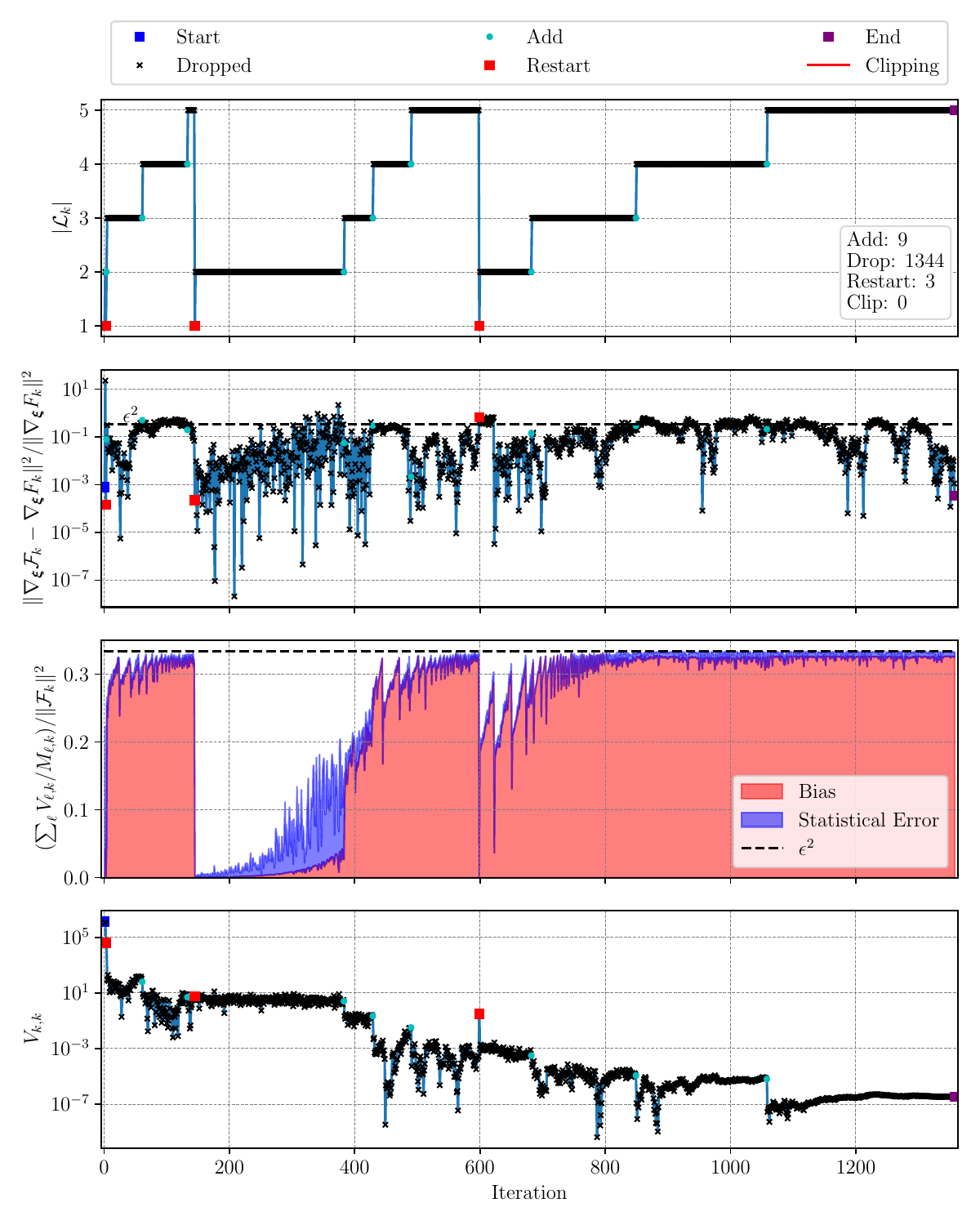}
  \caption{Single run, random quadratic example, Equation \eqref{eq:quad.example} with $\kappa=100$.
  From top to bottom, cardinality of the index set, true squared relative error, empirical relative error, and $\V{k}$ versus iteration.
  The starting point, the restarts, and the end are marked respectively as blue, red, and purple squares, iterations dropped with black $\times$, and the remaining \texttt{MICE} points with cyan circles.
  Dashed lines represent bounds used to control relative errors when applied.
  In the empirical relative error plot, we split the relative error between bias and statistical error.
  }
  \label{fig:quad_err}
\end{figure}

\begin{figure}[h]
  \centering
   \includegraphics[width=.75\linewidth]{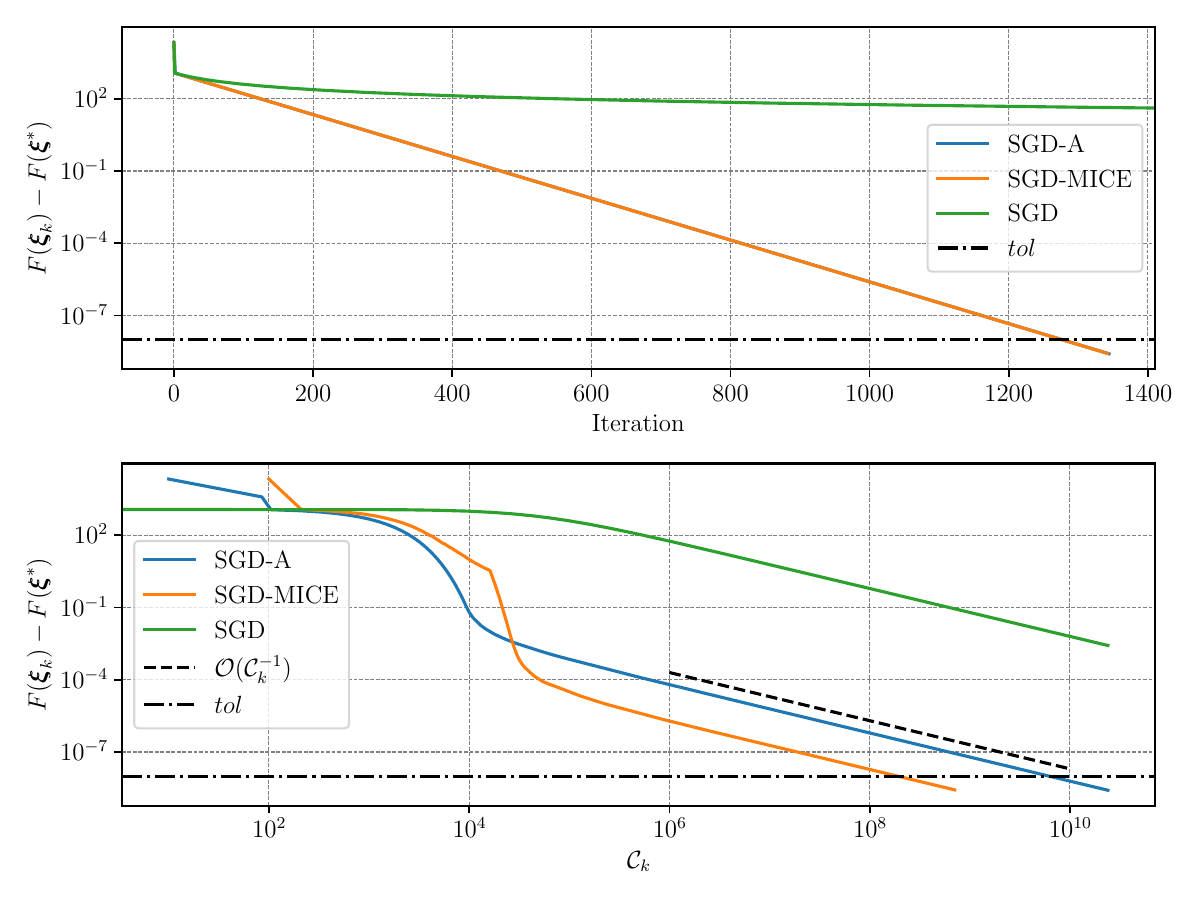}
   \caption{Single run, random quadratic example, Equation \eqref{eq:quad.example} with $\kappa=100$.
   Optimality gap versus iteration (top) and gradient sampling cost (bottom) for \SGD[-A], \SGD[-MICE], and \emph{vanilla }\SGD[].
   Dash-dotted lines represent $tol$, and the dashed line in the bottom plot illustrates the expected convergence rate of the optimality gap per cost, $\bigO{\cl{C}_k^{-1}}$.
   The top plot is limited to $1400$ iterations to illustrate \SGD[-A] and \SGD[-MICE] even though \SGD[] required close to $2.4 \times 10^6$ iterations.
   \SGD[-MICE] achieves $tol$ with less than $3\%$ of the sampling cost of \SGD[-A] and both achieve a much lower optimality gap than \SGD[] for the same cost.
   }
   \label{fig:quad_mice_mc}
 \end{figure}

\begin{figure}[h]
  \centering
   \includegraphics[width=.5\linewidth]{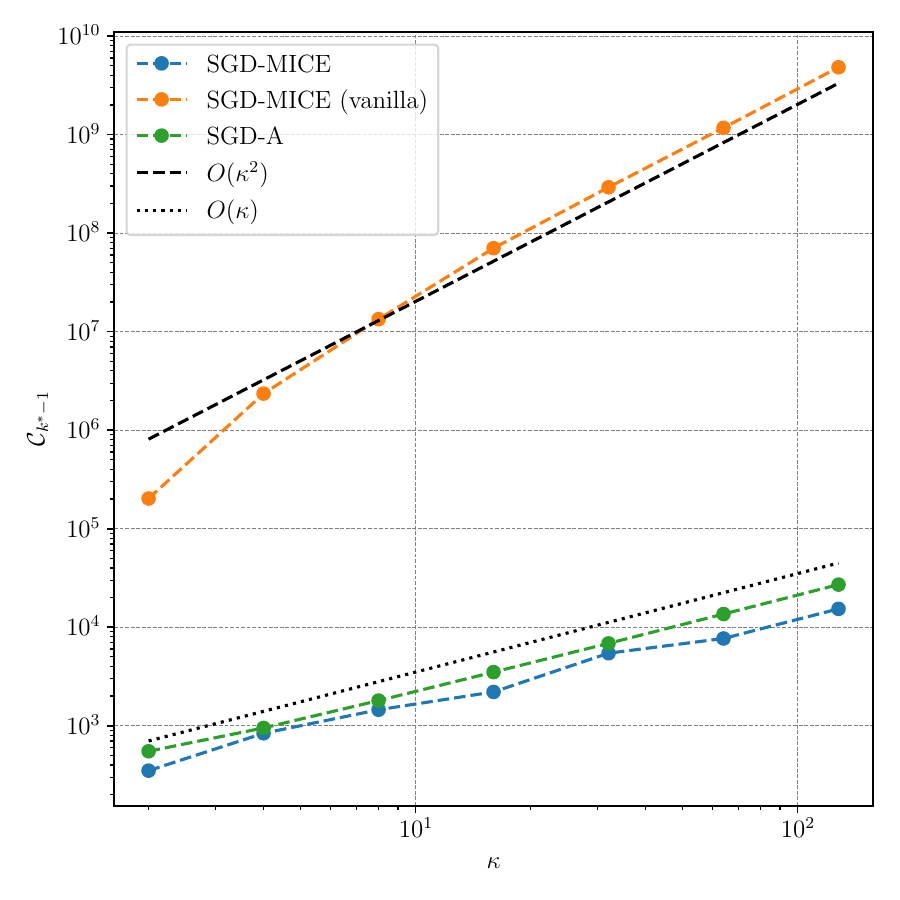}
   \caption{
    Gradient sampling cost versus condition number for \emph{vanilla} \texttt{SGD-MICE} (without \Restart[], \Drop[], or \Clip[]), \texttt{SGD-MICE} (with \Restart[], \Drop[], and \Clip[]), and \texttt{SGD-A}.
    The algorithms are run until they reach the stopping criterion defined as $\norm{\grad F(\xi_{k^*})}^2 < tol$.
    We also plot reference lines for $\bigO{\kappa^2}$ and $\bigO{\kappa}$.
    Note that \emph{vanilla} \texttt{SGD-MICE} cost increases as $\bigO{\kappa^2}$ as predicted in Corollary~\ref{cor:cost_sgd_mice} whereas \texttt{SGD-A} cost increases as $\bigO{\kappa}$, as predicted in Corollary~\ref{cor:cost_sgd_a}.
    Surprisingly, once the index set operators \Restart[], \Drop[], and \Clip[] are considered, \texttt{SGD-MICE} cost dramatically decreases, not only by a constant factor but effectively matching the rate of \SGD[-A] of $\bigO{\kappa}$.
   }
   \label{fig:quad_kappa_test}
 \end{figure}

In \S\ref{sec:sampling_cost}, we prove that, for expectation minimization, the gradient sampling cost necessary to reach a certain $\norm{\grad F(\xi_{k^*})}^2 < tol$ is $\bigO{\kappa^2 tol^{-1}}$ for \texttt{SGD-MICE} and $\bigO{\kappa tol^{-1} \log(tol^{-1})}$ for \texttt{SGD-A}.
To validate numerically the dependency of the cost with respect to the conditioning number, we evaluated both \texttt{SGD-MICE} and \texttt{SGD-A} with different condition numbers until the stopping criterion.
Moreover, we also tested \texttt{SGD-MICE} with and without the index set operators \Restart[], \Drop[], and \Clip[].
The reasoning for doing this test is that, in the analysis of Corollary~\ref{cor:cost_sgd_mice}, we consider the case where all iterates are kept in the index set.
However, in practice, one would expect \texttt{SGD-MICE} with the index set operators to perform better than both \emph{vanilla} \texttt{SGD-MICE} (without the operators) and \texttt{SGD-A}; in one extreme case where all iterates are kept, we recover \emph{vanilla} \texttt{SGD-MICE}, and in another extreme case we restart every iteration, resulting in \texttt{SGD-A}.
The gradient sampling cost versus $\kappa$ for these tests is presented in Figure~\ref{fig:quad_kappa_test}.

Figure~\ref{fig:quad_kappa_test} illustrates that enabling the index-set operators (notably \Restart[]) changes the observed scaling with the condition number from $\bigO{\kappa^2}$ (add-only, no operators) toward $\bigO{\kappa}$ (up to the fixed $\log(tol^{-1})$ factor for \SGD[-A]). This behavior is consistent with the fact that \SGD[-A] is a special case of \SGD[-MICE] obtained by restarting at every iteration. Our implementation selects at each iteration the operator that minimizes a per-iteration work proxy, so it can be interpreted as competing locally with this restart baseline. While this local comparison does not by itself yield a full global complexity guarantee for the adaptive operator selection (since future variance structure depends on past decisions), it provides intuition for why restarting can remove the additional $\kappa$ factor in the cost and why the operator-enhanced method empirically matches the $\bigO{\kappa}$ trend in this example.

\begin{figure}[h]
  \centering
  \includegraphics[width=.48\linewidth]{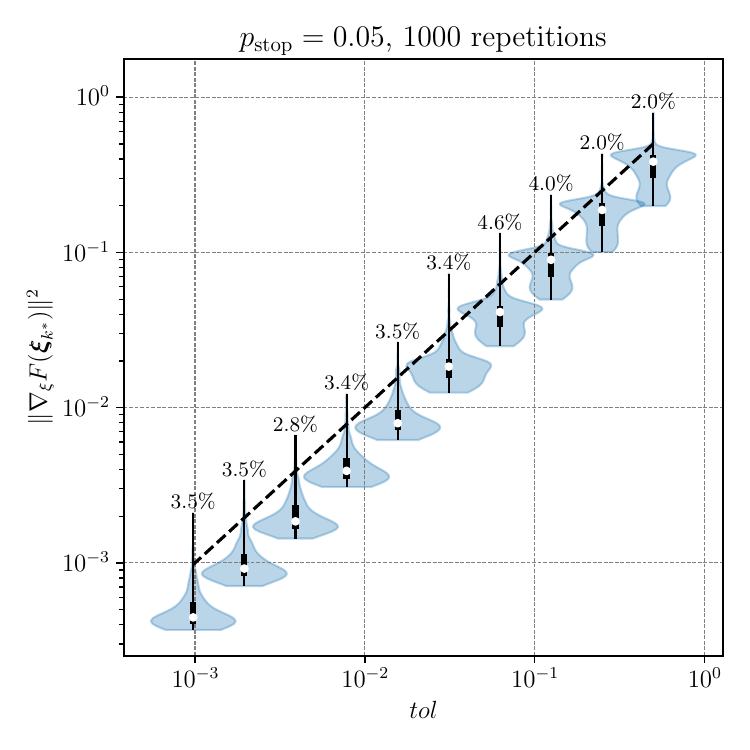}
  \includegraphics[width=.48\linewidth]{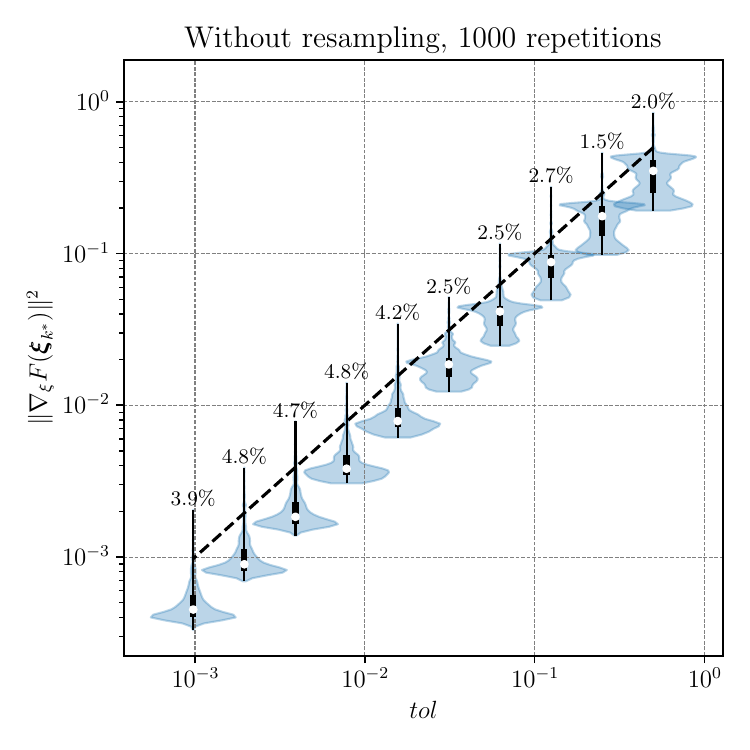}
  \caption{
  Random quadratic problem: consistency plot to validate the stopping criterion with the resampling technique (left) and without it (right).
  We present violin plots for the squared norm of the final gradient for different values of $tol$.
  The light blue shade represents an empirical pdf approximated using a Gaussian kernel density estimate, the thin hair represents the interval between maximum and minimum, the thick hair illustrates the quantiles between $0.25$ and $0.75$, and the white dot marks the median.
  A thousand independent runs were used to obtain the data presented.
  The dashed lines represent $tol$ and the percentage of runs with $\norm{\grad F(\xi_{k^*})}^2 > tol$ are presented for each $tol$.
  }
  \label{fig:consistenty}  
\end{figure}

In Remark~\ref{rmk:grad_resampling}, we present a resampling technique to take more informed decisions on stopping criterion and error control.
To validate our stopping criterion, we performed a thousand independent runs of \SGD[-MICE] for different values of $tol$, using the resampling to decide both sample sizes and the stopping criterion.
Figure~\ref{fig:consistenty} presents violin plots with approximations of empirical distributions of the squared gradient norms where optimization stopped and the percentage of times this quantity exceeded $tol$.
Moreover, we show both the case where we use the resampling technique and when we do not use it.
For lower tolerances, the resampling technique indeed reduced the percentage of premature stops, however, in both cases, a general trend of decrease following $tol$ is observed.
 
\FloatBarrier
\subsubsection{Benchmarks and ablations}
\begin{figure}[h]
  \centering
  \includegraphics[width=.7\linewidth]{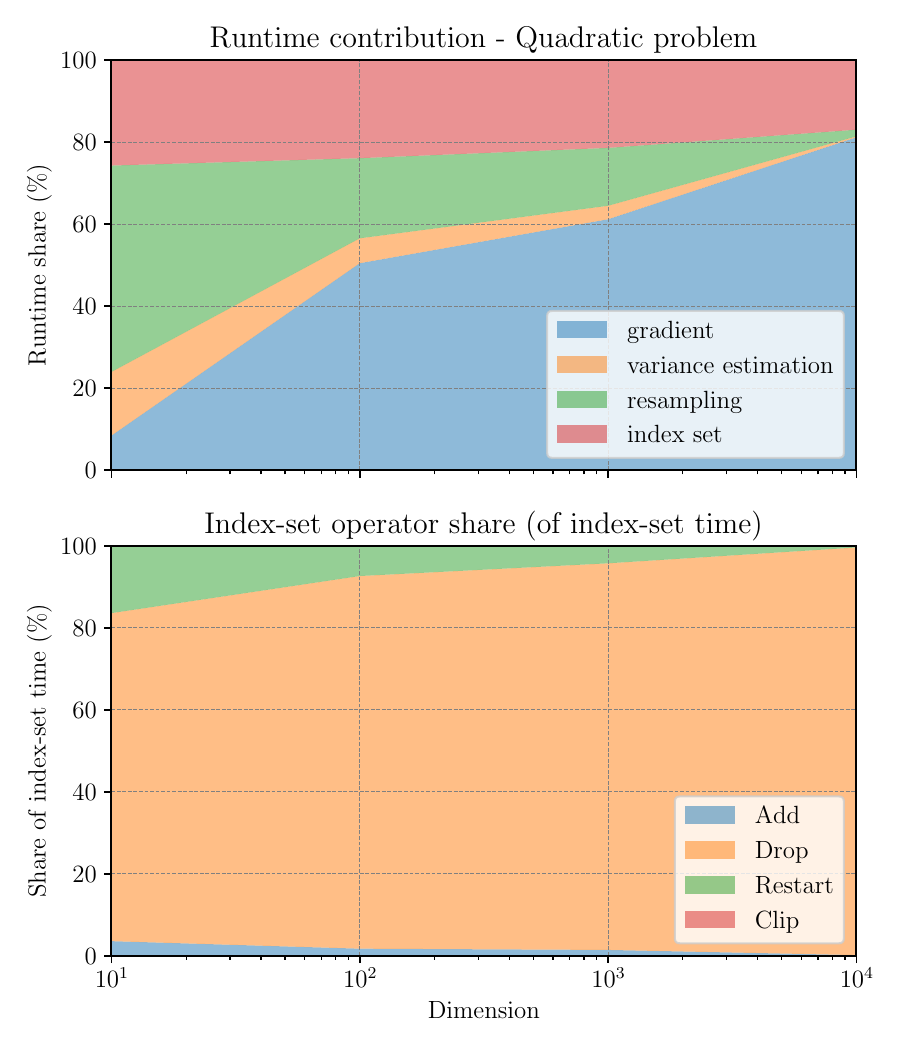}
  \caption{Quadratic problem benchmark: on top, the runtime breakdown for \texttt{SGD-MICE} versus dimensionality, and on the bottom, a breakdown of the runtime spent by each index set operator. The gradient evaluation cost dominates for dimensions $d \geq 100$, with operator overhead remaining under $10\%$ across all tested dimensions. Within the index-set operators, \Drop[] dominates the runtime because it requires gradient evaluations to determine whether to discard the most recent iterate.}
  \label{fig:benchmark}
\end{figure}

To evaluate the performance of \texttt{SGD-MICE} in higher-dimensional settings, we use a $d$-dimensional strongly-convex quadratic problem, which generalizes the two-dimensional setup of \S\ref{ex:1} (Eq.~\eqref{eq:quad.example}).
First, we benchmark the runtime of \texttt{SGD-MICE} and break it down into the contributions of gradient evaluations, variance estimation (using the Welford algorithm), resampling (as in Remark~\ref{rmk:grad_resampling}), and the index set operators defined in \S\ref{sec:opt_index_set}.
Moreover, we break the index set benchmark further down into the contributions of \Add[], \Drop[], \Restart[], and \Clip[].
The average runtime shares over $n=50$ runs (with a cap of $100{,}000$ gradient evaluations per run) are presented in Figure~\ref{fig:benchmark}.

To quantify the contribution of each index-set operator under a fixed gradient budget, we run an ablation on the same quadratic setup (dimension $d=100$, condition number $\kappa=100$) with a cap of $100{,}000$ gradient evaluations per run.
We repeat each configuration $n=50$ times and report the mean and standard deviation of the relative optimality (final optimality gap divided by initial gap) and the mean number of gradient evaluations used.
The configurations are: \Add only (no \Drop, \Restart, or \Clip); \Add+\Drop; \Add+\Drop+\Restart; and All+Clip (\Add, \Drop, \Restart, and \Clip).
Table~\ref{tab:quad_operator_ablation} summarizes the results; the mean number of times each operator was triggered (drop, restart, clip) is also reported.

\begin{table}[h]
    \centering
  \caption{Operator ablation on the random quadratic problem ($d=100$, $\kappa=100$): fixed budget of $100{,}000$ gradient evaluations, $n=50$ runs. Relative optimality is final gap/initial gap.}
  \label{tab:quad_operator_ablation}
  \begin{tabular}{lccccccc}
    \toprule
    Configuration & Mean (rel.\ opt.) & Std & Mean grad.\ evals & Drop & Restart & Clip \\
    \midrule
    \Add only      & $0.0120$ & $0.00230$ & $99{,}708$ & 0.0 & 0.0 & 0.0 \\
    \Add+\Drop     & $0.00974$ & $0.00077$ & $99{,}719$ & 184.7 & 0.0 & 0.0 \\
    \Add+\Drop+\Restart & $1.18 \times 10^{-9}$ & $9.48 \times 10^{-10}$ & $99{,}645$ & 947.6 & 2.9 & 0.0 \\
    All+\Clip     & $1.54 \times 10^{-9}$ & $1.33 \times 10^{-9}$ & $99{,}730$ & 940.9 & 1.3 & 5.8 \\
    \bottomrule
  \end{tabular}
  \end{table}

The operator ablation results in Table~\ref{tab:quad_operator_ablation} reveal that the \Restart[] operator is critical for achieving high-quality solutions under a fixed gradient budget: enabling \Restart[] reduces the mean relative optimality from approximately $10^{-2}$ (\Add+\Drop) to $10^{-9}$ (\Add+\Drop+\Restart), a dramatic improvement of seven orders of magnitude.
In contrast, \Drop[] alone yields only a modest improvement over the add-only baseline (from $0.0120$ to $0.00974$), and \Clip[] has minimal additional impact in this quadratic setting.
The frequency counts show that \Drop[] is invoked frequently (approximately $940$ times per run) once \Restart[] is enabled, whereas \Restart[] itself triggers only $1$--$3$ times and \Clip[] is rarely needed ($\sim 6$ times).
These findings confirm that our greedy index set policy effectively balances the trade-off between reusing accumulated gradient information (via \Add[] and selective \Drop[]) and periodically resetting the estimator (via \Restart[]).

Table~\ref{tab:quad_epsilon_sensitivity} shows sensitivity to the relative-error tolerance $\epsilon$ (All+Clip, same quadratic setup, budget $100{,}000$, $n=50$ runs per value).
Under a fixed gradient budget, larger values of $\epsilon$ yield better final optimality; increasing $\epsilon$ from $0.3$ to $1.0$ improves the mean relative optimality from $3.64 \times 10^{-9}$ to $7.94 \times 10^{-10}$.
This behavior is expected, since a looser error tolerance allows the algorithm to spend fewer gradients per iteration (smaller batch sizes at each level), thereby enabling more iterations and further progress toward the optimum within the same total gradient budget.
While the theoretical cost analysis in \S\ref{sec:sampling_cost} identifies $\epsilon = 1/\sqrt{3}$ as minimizing the leading-order sampling cost under the variance-sum sizing rule, the present fixed-budget experiments suggest that the optimal choice of $\epsilon$ depends on the problem at hand.

Table~\ref{tab:quad_delta_sensitivity} shows sensitivity to $\delta_{\mathrm{drop}}$ (left) and $\delta_{\mathrm{rest}}$ (right): relative optimality and mean number of drop/restart events (All+Clip; for the drop sweep $\delta_{\mathrm{rest}}=0$, for the restart sweep $\delta_{\mathrm{drop}}=0.5$). Same quadratic setup, budget $100{,}000$, $n=50$ runs per value.
The results indicate that performance is fairly robust to the choice of $\delta_{\mathrm{drop}}$: all tested values ($0.0$ to $1.0$) achieve similar relative optimality and trigger comparable numbers of drop events ($871$ to $947$ per run).
In contrast, the restart threshold $\delta_{\mathrm{rest}}$ exhibits a more pronounced effect: increasing $\delta_{\mathrm{rest}}$ from $0$ to $1.0$ improves mean relative optimality by roughly a factor of two (from $1.86 \times 10^{-9}$ to $8.59 \times 10^{-10}$) while modestly increasing the restart count from $1.4$ to $1.9$ events per run.
This suggests that a more aggressive restart policy (larger $\delta_{\mathrm{rest}}$) is beneficial in this quadratic setting, though the overall frequency of restarts remains low.

Table~\ref{tab:quad_max_index_sensitivity} shows sensitivity to the maximum index set cardinality $\max_k |\cl{L}_k|$ (All+Clip, same quadratic §setup, budget $100{,}000$, $n=50$ runs per value).
The results demonstrate robustness to the choice of maximum index-set size: increasing the cap from $100$ to $1000$ yields only a slight improvement in mean relative optimality (from $1.53 \times 10^{-9}$ to $1.44 \times 10^{-9}$).
This indicates that, for this problem, the greedy operator policy naturally maintains a relatively compact index set, so that even a modest cap of $100$ does not severely constrain performance.
In practice, imposing a cap remains useful to bound memory usage and per-iteration overhead, and the present experiments confirm that moderate caps (e.g., $100$--$500$) suffice for effective optimization on this quadratic benchmark.

\begin{table}[h]
    \centering
  \caption{Sensitivity to $\epsilon$ ($d=100$, $\kappa=100$, budget $100{,}000$, $n=50$ runs per value).}
  \label{tab:quad_epsilon_sensitivity}
  \begin{tabular}{lcc}
    \toprule
    $\epsilon$ & Mean (rel.\ opt.) & Std \\
    \midrule
    0.3   & $3.64 \times 10^{-9}$ & $2.16 \times 10^{-9}$ \\
    0.5   & $1.76 \times 10^{-9}$ & $1.21 \times 10^{-9}$ \\
    0.577 & $1.26 \times 10^{-9}$ & $7.18 \times 10^{-10}$ \\
    0.7   & $1.17 \times 10^{-9}$ & $9.02 \times 10^{-10}$ \\
    1.0   & $7.94 \times 10^{-10}$ & $6.16 \times 10^{-10}$ \\
    \bottomrule
  \end{tabular}
  \end{table}

\begin{table}[h]
    \centering
  \caption{Sensitivity to $\delta_{\mathrm{drop}}$ and $\delta_{\mathrm{rest}}$ ($d=100$, $\kappa=100$, budget $100{,}000$, $n=50$ runs per value).}
  \label{tab:quad_delta_sensitivity}
  \begin{tabular}{lccccc}
    \toprule
    \multicolumn{3}{c}{Sensitivity to $\delta_{\mathrm{drop}}$} & \multicolumn{3}{c}{Sensitivity to $\delta_{\mathrm{rest}}$} \\
    \cmidrule(lr){1-3} \cmidrule(l){4-6}
    $\delta_{\mathrm{drop}}$ & Rel.\ opt. & Mean drops & $\delta_{\mathrm{rest}}$ & Rel.\ opt. & Mean restarts \\
    \midrule
    0.0   & $1.30 \times 10^{-9}$ & 871.4 & 0.0   & $1.86 \times 10^{-9}$ & 1.4 \\
    0.25  & $1.25 \times 10^{-9}$ & 936.9 & 0.25  & $1.11 \times 10^{-9}$ & 1.6 \\
    0.5   & $1.45 \times 10^{-9}$ & 940.8 & 0.5   & $1.14 \times 10^{-9}$ & 1.6 \\
    0.75  & $1.63 \times 10^{-9}$ & 940.8 & 0.75  & $1.03 \times 10^{-9}$ & 1.7 \\
    1.0   & $1.50 \times 10^{-9}$ & 946.5 & 1.0   & $8.59 \times 10^{-10}$ & 1.9 \\
    \bottomrule
  \end{tabular}
  \end{table}

\begin{table}[h]
    \centering
  \caption{Sensitivity to maximum index set cardinality ($d=100$, $\kappa=100$, budget $100{,}000$, $n=50$ runs per value).}
  \label{tab:quad_max_index_sensitivity}
  \begin{tabular}{lcc}
    \toprule
    Max $|\cl{L}_k|$ & Mean (rel.\ opt.) & Std \\
    \midrule
    100  & $1.53 \times 10^{-9}$ & $1.20 \times 10^{-9}$ \\
    500  & $1.48 \times 10^{-9}$ & $1.07 \times 10^{-9}$ \\
    1000 & $1.44 \times 10^{-9}$ & $1.05 \times 10^{-9}$ \\
    \bottomrule
  \end{tabular}
  \end{table}

\FloatBarrier
\subsection{Stochastic Rosenbrock function}\label{ex:2}

The goal of this example is to test the performance of {\texttt{Adam-MICE}, that is, \texttt{Adam} coupled with our gradient estimator \texttt{MICE}}, in minimizing the expected value of the stochastic Rosenbrock function in \eqref{eq:rosenbrock}, showing that \texttt{MICE} can be coupled with {different} first-order optimization methods in a non-intrusive manner.
Here we adapt the deterministic Rosenbrock function to the stochastic setting, specializing our optimization problem \eqref{eq:underlying.problem} with
\begin{equation}
 f(\bs{\xi},\bs{\theta}) = \left(a - \xi_0 + \theta_{0}\right)^{2} + b \left(- \xi_0^{2} + \xi_1 + \theta_{0}^{2} - \theta_{1}^2\right)^{2},
\end{equation}
where $a=1$, $b=100$, $\theta_0,\theta_1\sim\cl{N}(0,\sigma_\theta^2)$.
The objective function to be minimized is thus
\begin{equation}\label{eq:rosenbrock}
F(\bs{\xi}) = (a - \xi_0)^2  + \sigma_\theta^{2} + b\left(4 \sigma_\theta^{4} + (\xi_1 - \xi_0^{2})^2\right),
\end{equation}
{and its} gradient {is given by}
\begin{equation}
 \nabla_{\bs{\xi}}F(\bs{\xi})=
 \begin{bmatrix}
  - 2 a + 4 b \xi_0^{3} - 4 b \xi_0 \xi_1 + 2 \xi_0 \\
  - 2 b \xi_0^{2} + 2 b \xi_1
 \end{bmatrix},
\end{equation}
which coincides with the gradient of the deterministic Rosenbrock function.
Therefore, the optimal point of the stochastic Rosenbrock is the same as the one of the deterministic: $\bs{\xi}^*=(a,a^2)$.
To perform the optimization, we sample the stochastic gradient
\begin{equation}
 \nabla_{\bs{\xi}}f(\bs{\xi}, \bs{\theta})=
 \begin{bmatrix}
  - 2 a + 4 b \xi_0\left(\xi_0^{2} - \xi_1 - \theta_0^2 + \theta_1^2 \right) + 2 \xi_0 - 2 \theta_0 \\
  2 b \left(-\xi_0^{2} + \xi_1 + \theta_0^2 - \theta_1^2 \right)
 \end{bmatrix}.
\end{equation}

Although this is still a low dimensional example, minimizing the Rosenbrock function poses a difficult optimization problem for first-order methods; these tend to advance slowly in the region where the gradient has near-zero norm.
Moreover, when noise is introduced in gradient estimates, their relative error can become large, affecting the optimization convergence.

We compare the convergence of the {classical} \texttt{Adam} algorithm against \texttt{Adam-MICE}.
To illustrate the effect of the dispersion of the random variable $\bs{\theta}$, two distinct noise levels are considered, namely $\sigma_\theta=10^{-4}$ and $\sigma_\theta=10^{-1}$.
As for the optimization setup, we set \texttt{Adam-MICE} with fixed step-size $0.3$ and \texttt{Adam} with a decreasing step-size $\eta_k = 0.01 / \sqrt{k}$, which we observed to be the best step-sizes for each method.
The stopping criterion for both algorithms is set as $10^7$ gradient evaluations.
For \texttt{Adam-MICE}, we use $\epsilon = 1$, whereas for \texttt{Adam} we use a fixed batch size of $100$.
In all cases, we start the optimization from $\bs{\xi}_0 = (-1.5, 2.5)$

\begin{figure}
 \includegraphics[width=.48\linewidth]{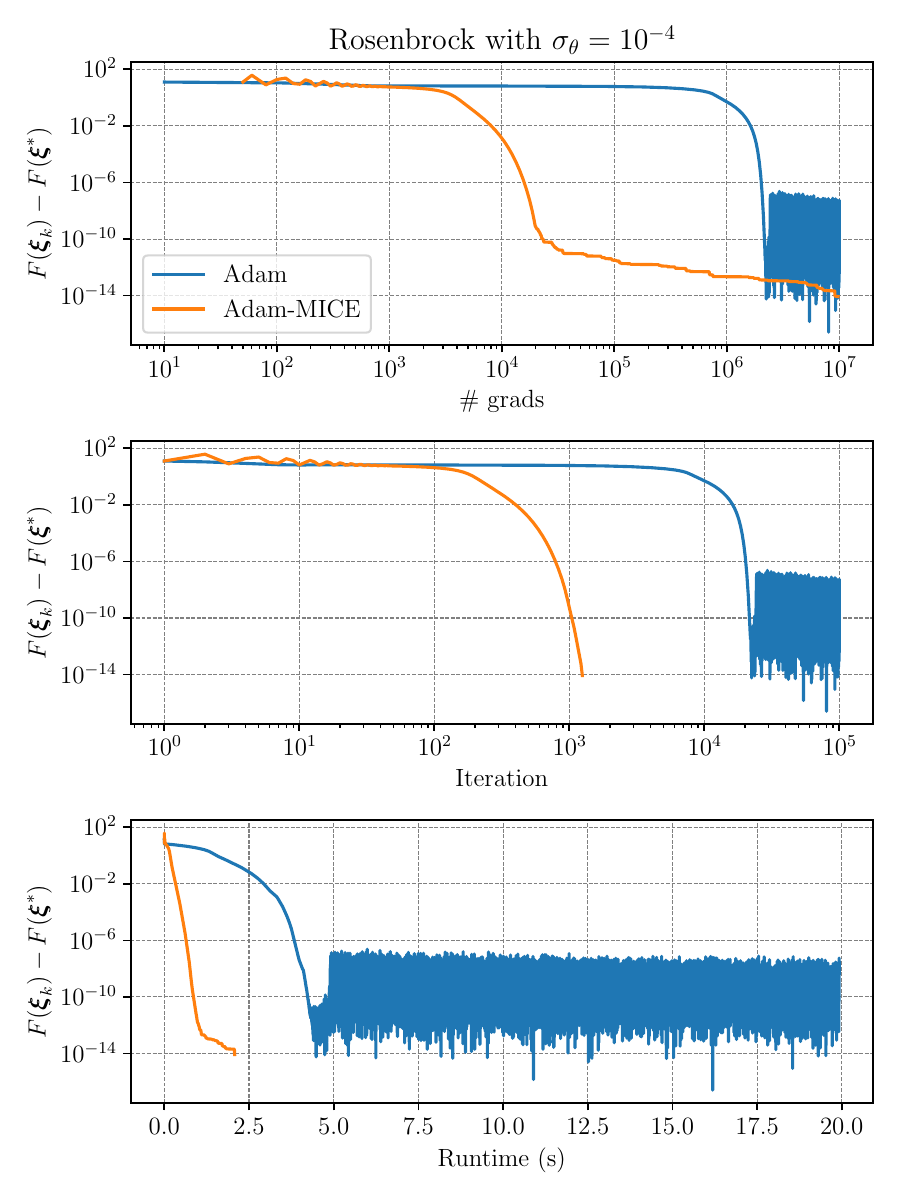}
 \includegraphics[width=.48\linewidth]{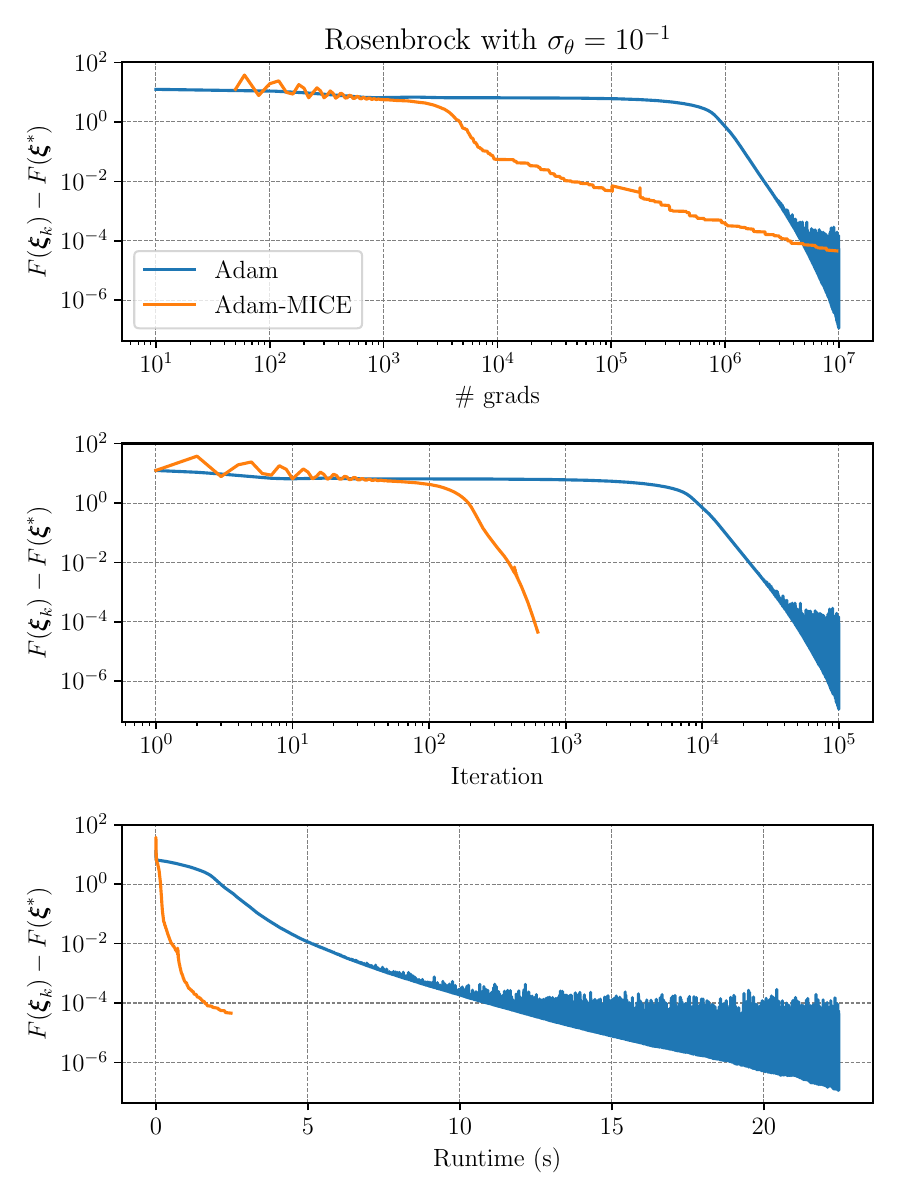}
 \caption{Single run, stochastic Rosenbrock function example, \eqref{eq:rosenbrock} with $\sigma_{\theta}=10^{-4}$.
 Optimality gap for \texttt{Adam} and \texttt{Adam-MICE} versus the number of gradient evaluations (top), iterations (center), and runtime in seconds (bottom).}
 \label{fig:rosen_1}
\end{figure}

In Figure~\ref{fig:rosen_1}, we present, for $\sigma_\theta$ of $10^{-4}$ (left) and $10^{-1}$ (right), the optimality gap for both \texttt{Adam} and \texttt{Adam-MICE} versus the number of gradients, iterations, and runtime in seconds.
It is clear that \texttt{Adam-MICE} is more stable than \texttt{Adam} as the latter oscillates as it approximates the optimal point in both cases.
The efficient control of the error in gradient estimates allows \texttt{Adam-MICE} to converge monotonically in the asymptotic phase.
Moreover, the number of iterations and the runtime are much smaller for \texttt{Adam-MICE} than for \texttt{Adam}.

As a conclusion, even though \Adam[] has its own mechanisms to control the statistical error of gradients, coupling it with \MICE[], for this example, has proven to be advantageous as it allows more evaluations to be performed simultaneously.
Moreover, as the gradient error is controlled, we can use \Adam[] with a fixed step-size.
Also, \MICE[] allows for a stopping criterion based on the gradient norm, which would not be possible for \emph{vanilla} \Adam[].

\FloatBarrier

\subsection{Logistic regression} \label{ex:3}
In this example, we train logistic regression models using \texttt{SGD-MICE}, \texttt{SAG} \cite{schmidt2017minimizing}, \texttt{SAGA} \cite{defazio2014saga}, \texttt{SARAH} \cite{nguyen2017sarah}, and \texttt{SVRG} \cite{Joh13} to compare their performances.
Here, we present a more practical application of \texttt{MICE}, where we can test its performance on high-dimensional settings with finite populations.
Therefore, we calculate the error as in \eqref{eq:finite_sum_problem} and use Algorithm~\ref{alg:finite_sum} to obtain the optimal sample sizes.
To train the logistic regression model for binary classification, we use the $\ell_2$-regularized log-loss function
\begin{equation}\label{eq:logistic}
 F(\bs{\xi}) = \frac{1}{N} \sum_{i=1}^N f\left(\bs{\xi}, \bs{\theta}_i= (\bs{x}_i, y_i)\right) = \frac{1}{N} \sum_{i=1}^N \log(1 + \exp(-y_i \bs{\xi} \cdot \bs{x}_i))
 + \frac{\lambda}{2} \norm{\bs{\xi}}^2,
\end{equation}
where each data point $(\bs{x}_i, y_i)$ is such that $\bs{x}_i \in \bb{R}^{d_{\bs{\xi}}}$ and $y_i \in \{-1,1\}$.
We use the datasets \emph{mushrooms}, \emph{gisette}, and \emph{HIGGS}, obtained from LibSVM\footnote{
 \url{https://www.csie.ntu.edu.tw/~cjlin/libsvmtools/datasets/binary.html}
}.
The size of the datasets $N$, number of features $d_{\bs{\xi}}$, and regularization parameters $\lambda$ are presented in Table~\ref{tab:datasets}.
\begin{table}[h]
 \caption{Size, number of features, and regularization parameters for the datasets used in the logistic regression example.}
 \label{tab:datasets}
 \begin{tabular}{lrrrrr}
  \toprule
  \multicolumn{1}{c}{Dataset} & 
  \multicolumn{1}{c}{Size} &
  \multicolumn{1}{c}{Features} &
  \multicolumn{1}{c}{$\lambda$} &
  \multicolumn{1}{c}{$\kappa$} &
  \multicolumn{1}{c}{epochs} \\
  \midrule
  \emph{mushrooms} & 8124     & 112                & $10^{-5}$ & 12316.30 & 100\\
  \emph{gisette}   & 6000     & 5000               & $10^{-4}$ & 1811.21  & 50\\
  \emph{HIGGS}     & 11000000 & 28                 & $10^{-4}$ & 765.76   & 10\\ \bottomrule
 \end{tabular}
\end{table}

When using \texttt{SGD-MICE} for training the logistic regression model, we use $\epsilon = 1/\sqrt{3}$.
For the other methods, we use batch sizes of size $10$.
Since we have finite populations, we use Algorithm~\ref{alg:finite_sum} to calculate the sample-sizes.
\texttt{SGD-MICE} step is based on the Lipschitz smoothness of the true objective function as presented in Proposition~\ref{prp:pl_convergence}. 
Conversely, the other methods rely on a Lipschitz constant that must hold for all data points, which we refer to as $\hat{L}$.
A maximum index set cardinality of $100$ is imposed on \texttt{SGD-MICE}; if $|\cl{L}_k|=100$, we restart the index set.
For the \emph{HIGGS} dataset, due to the smaller dimensionality (28 variables), we used a maximum index set cardinality of $1000$.
The step-sizes for \texttt{SAG}, \texttt{SAGA}, \texttt{SARAH}, and \texttt{SVRG} are presented in Table~\ref{tab:steps}.
These steps were chosen as the best performing for each case based on the recommendations of their original papers.
\begin{table}[h]
 \caption{Step-sizes chosen for each method for the logistic regression example.}
 \label{tab:steps}
 \begin{tabular}{lcccc}
  \toprule
  Method    & \texttt{SAG}                          & \texttt{SAGA}                        & \texttt{SARAH}              & \texttt{SVRG}               \\ \midrule
  Step-size & $\frac{1}{16(\hat{L} + \mu N)}$ & $\frac{1}{2(\hat{L} + \mu N)}$ & $\frac{1}{2 \hat{L}}$ & $\frac{1}{2 \hat{L}}$ \\ \bottomrule
 \end{tabular}
\end{table}

\begin{figure}[h]
  \includegraphics[width=.49\linewidth]{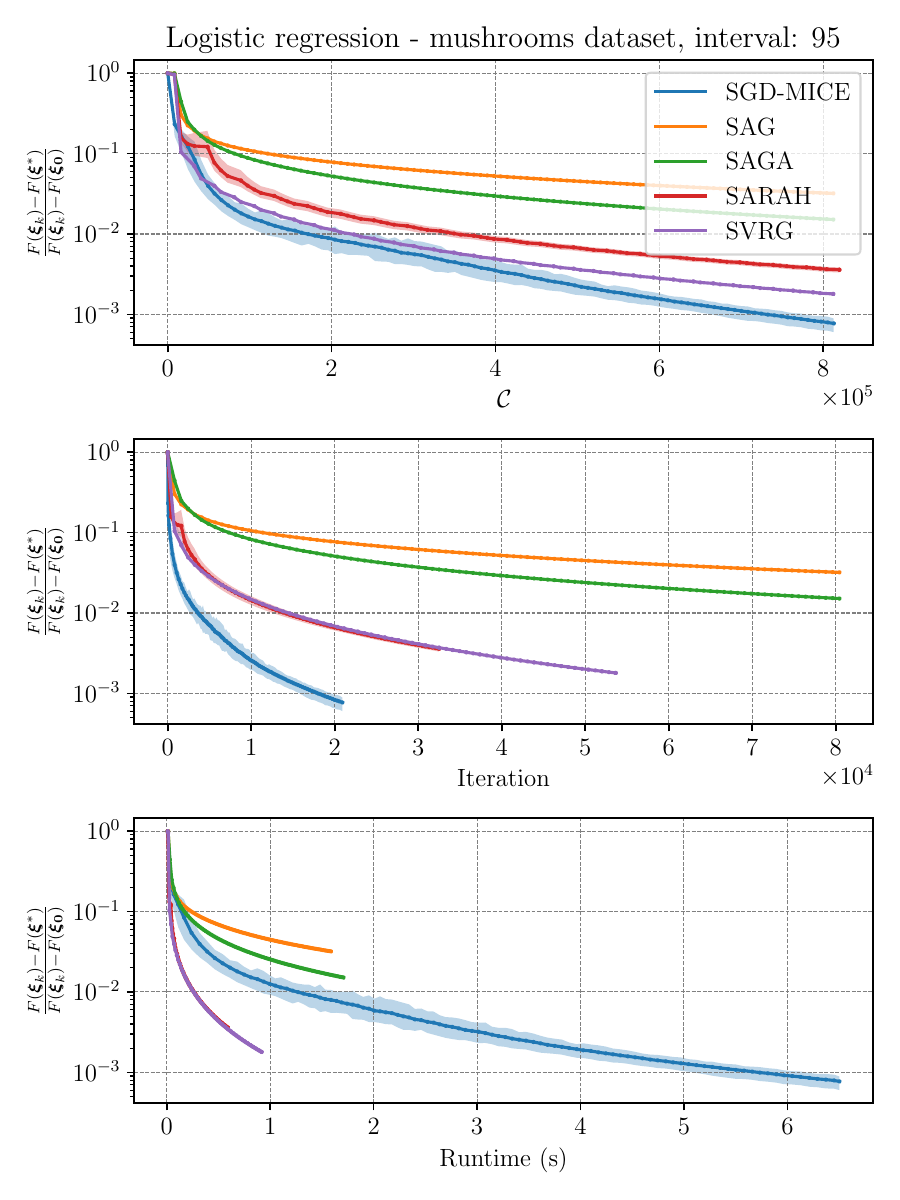}
  \includegraphics[width=.49\linewidth]{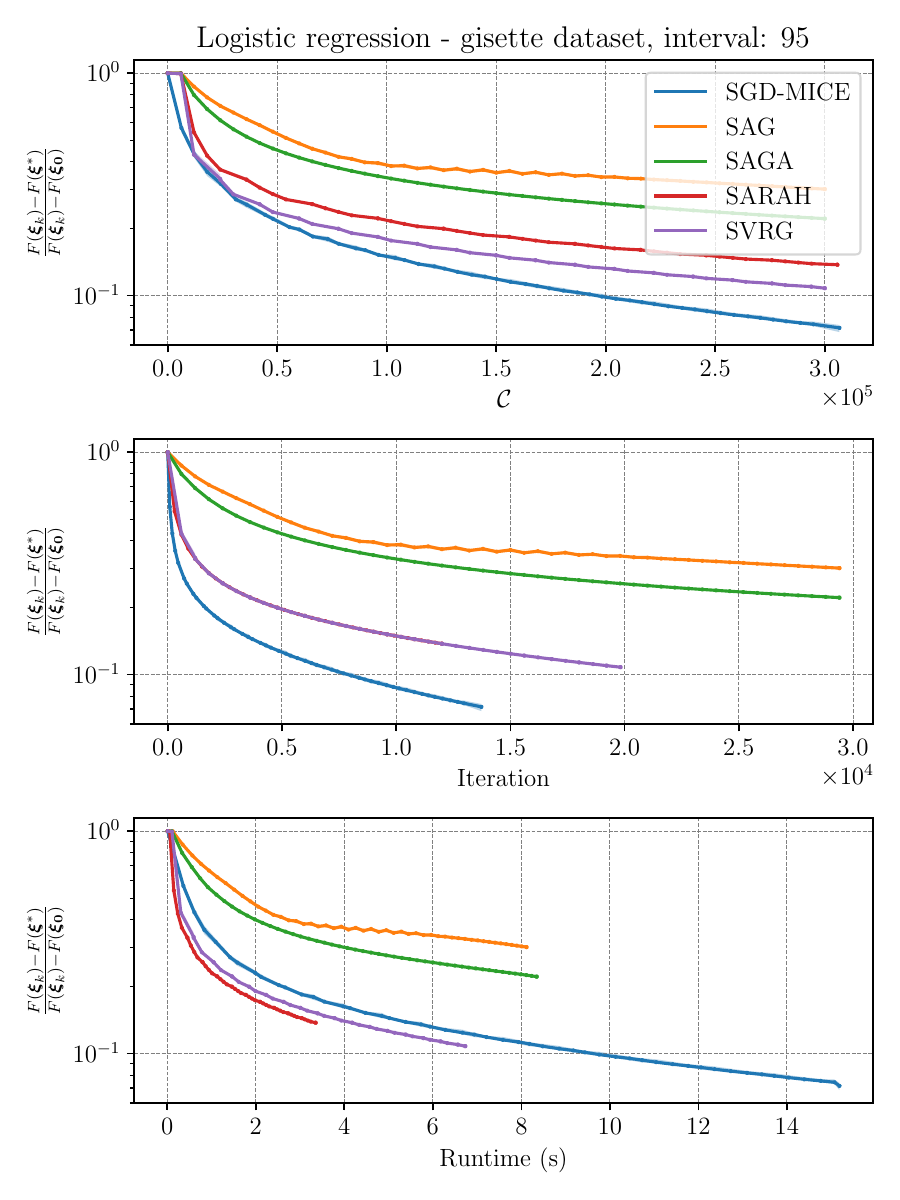}
  \includegraphics[width=.49\linewidth]{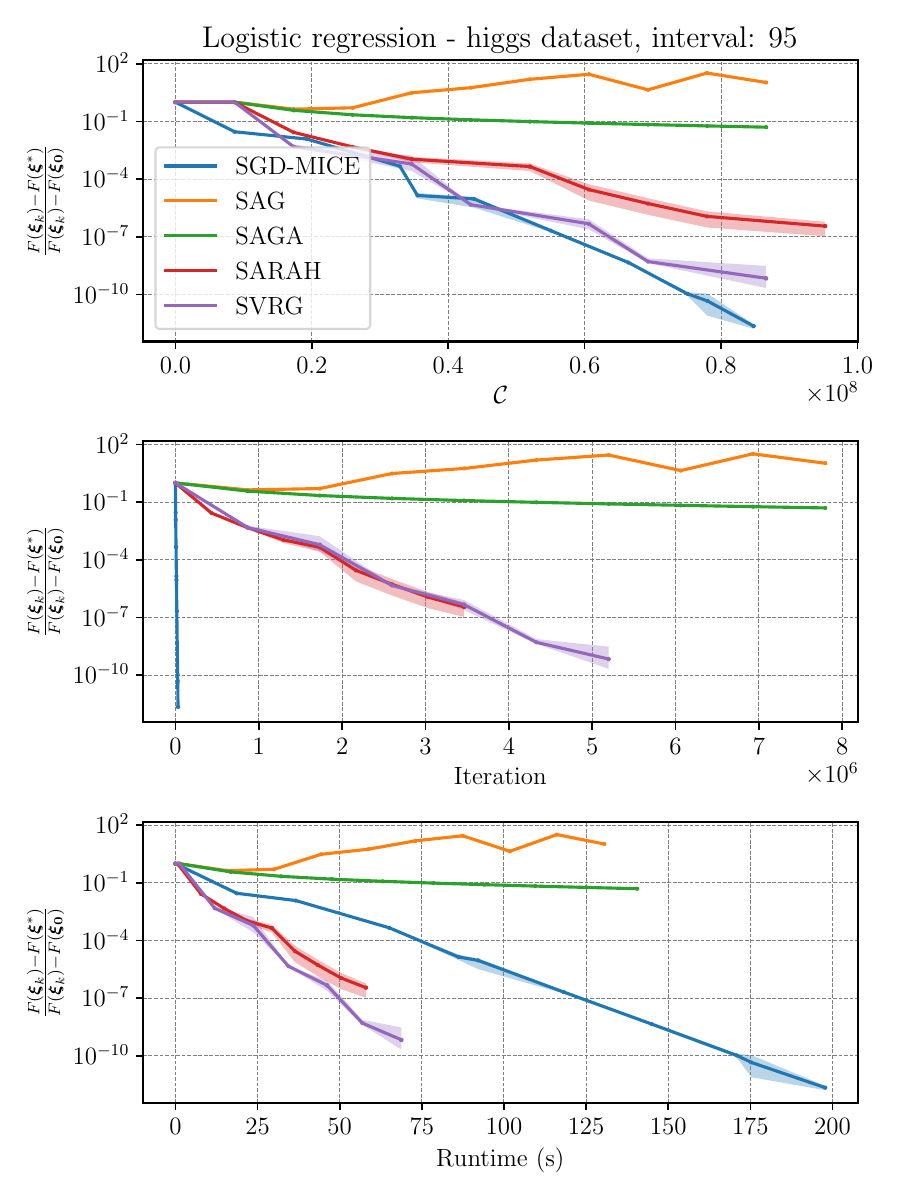}
  \caption{A hundred runs, logistic regression examples. 
  Relative optimality gap versus gradient sampling cost (top), iterations (center), and runtime in seconds (bottom) for \texttt{SGD-MICE}, \texttt{SAG}, \texttt{SAGA}, \texttt{SARAH}, and \texttt{SVRG}.
  The shaded regions represent confidence intervals between percentiles encompassing $95\%$ of values.
  }
  \label{fig:logreg}
\end{figure}

To evaluate the consistency of \SGD[-MICE] versus the other baseline methods, we perform $100$ independent runs of each method for each dataset.
Figure~\ref{fig:logreg} presents confidence intervals and medians of the relative optimality gap (the optimality gap normalized by its starting value) for the \emph{mushrooms}, \emph{gisette}, and \emph{HIGGS} datasets versus {the} number of gradient evaluations, iterations, and runtimes in seconds.

In the three cases studied, \texttt{SGD-MICE} performed better than the baseline methods in terms of optimality gap decrease per gradient sampling cost and number of iterations.
However, the runtime of \texttt{SGD-MICE} exceeded that of the other methods due to the \MICE[] estimator overhead.
We note, however, that the gradients in this logistic regression model are fairly cheap; in the case where gradient evaluations are more expensive, the \MICE[] overhead becomes negligible, cf. Figure~\ref{fig:benchmark}.

\begin{figure}
  \includegraphics[width=.8\linewidth]{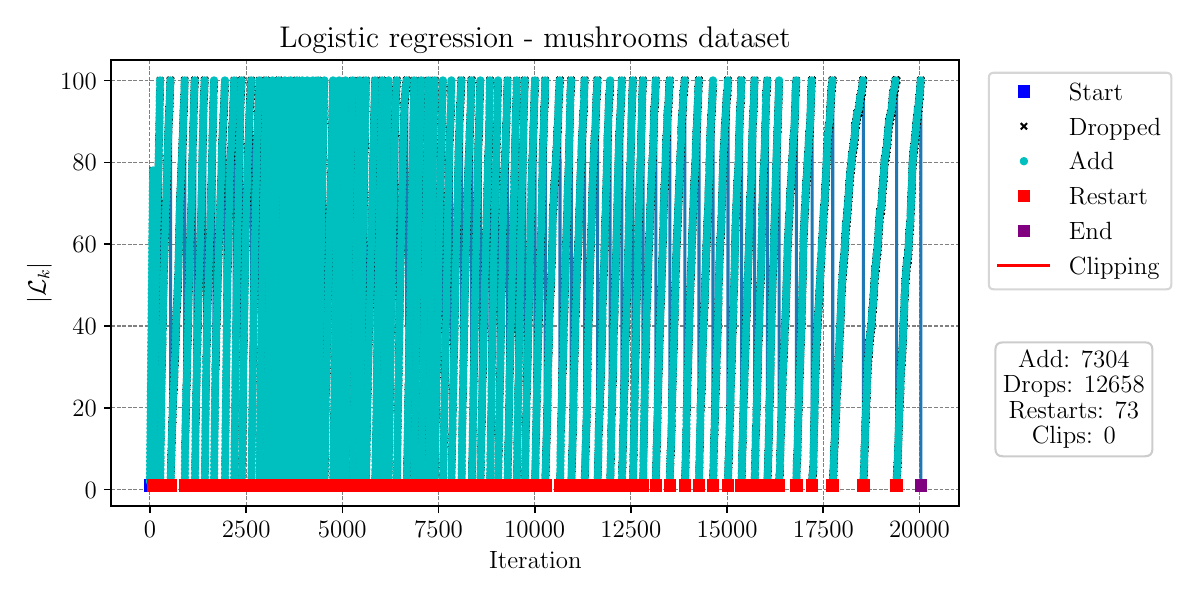}
  \includegraphics[width=.8\linewidth]{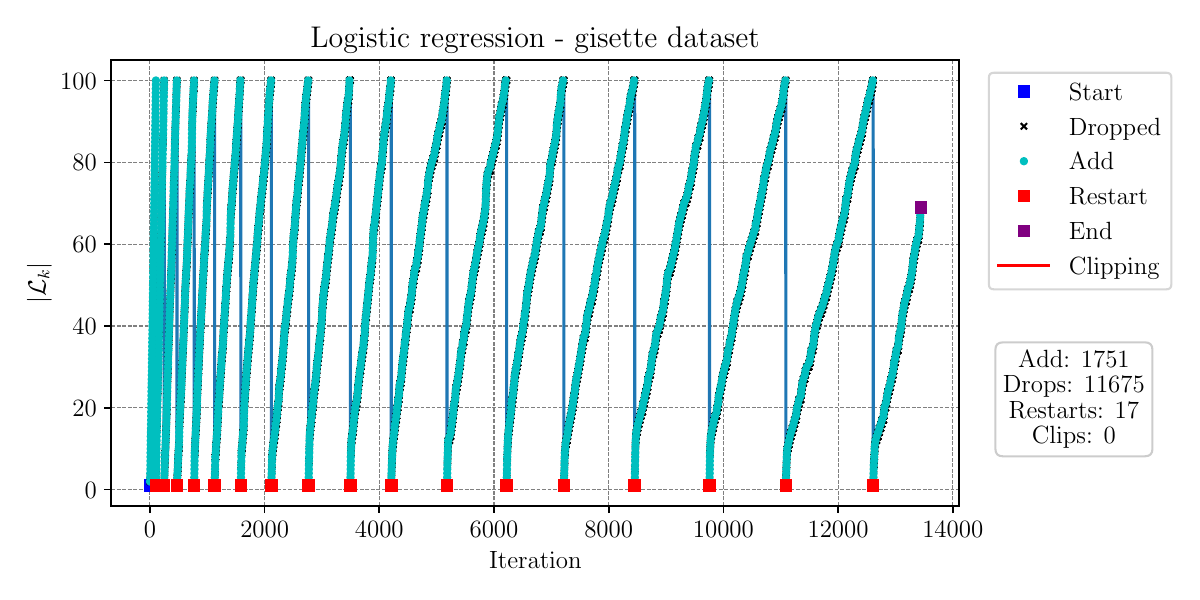}
  \includegraphics[width=.8\linewidth]{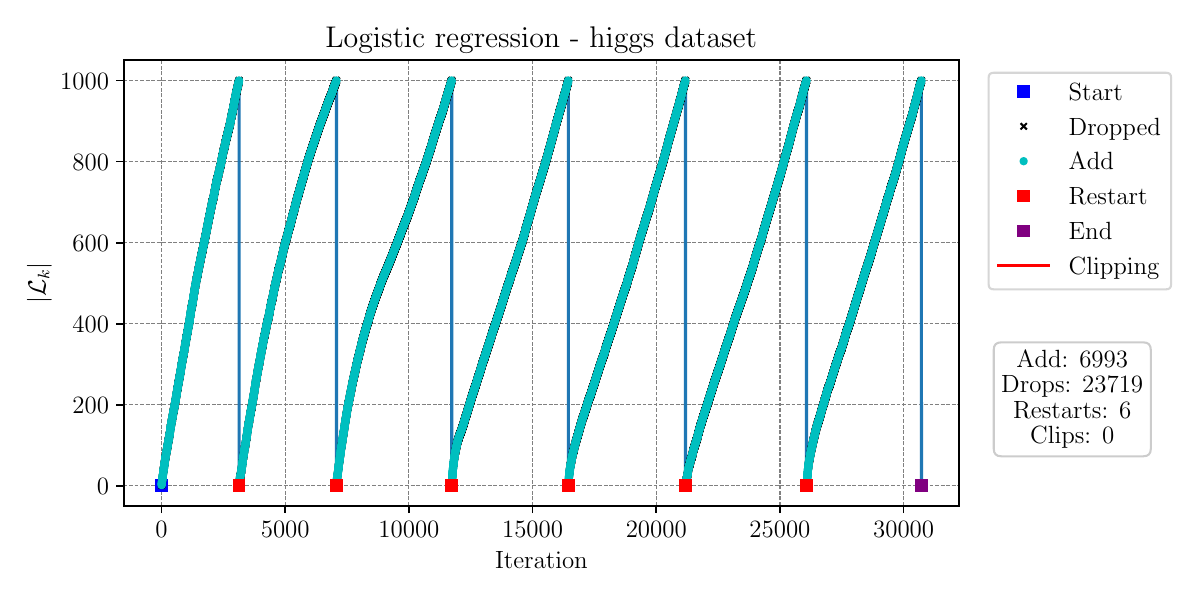}
  \caption{Index set cardinality versus iteration for the logistic regression of the \emph{mushrooms} dataset (top), \emph{gisette} dataset (center), and \emph{HIGGS} dataset, (bottom).
  We mark the dropped iteration as black $\times$'s, the iterations kept at the index set as cyan circles, restarts as red squares, and clippings as red lines.
  }
  \label{fig:logreg_length}
\end{figure}

Figure~\ref{fig:logreg_length} presents the index set cardinalities versus iterations of \texttt{SGD-MICE} for the three datasets.
Moreover, we present the iterations that were kept in the index set, the ones that were dropped, as well as restarts and clippings.

From the results of this example, we observe that \texttt{MICE} performs well in problems with a reasonably large number of parameters, for instance, $5000$ in the \emph{gisette}, and finite dataset populations ranging from the thousands to the millions.
One can conclude from the results obtained that \texttt{SGD-MICE} performs consistently better in terms of gradient sampling cost compared to the other methods.
Note that both \texttt{SAG} and \texttt{SAGA} need to decrease their step-sizes as the sample-size increases, and that \texttt{SARAH} and \texttt{SVRG} need to reevaluate the full-gradient after a few epochs to keep their convergence.

\section{Acknowledgments}

This publication is based upon work supported by the King Abdullah University of Science and Technology (KAUST) Office of Sponsored Research (OSR) under Award No. OSR-2019-CRG8-4033, the Alexander von Humboldt Foundation, and Coordination for the Improvement of Higher Education Personnel (CAPES).

Last but not least, we want to thank Dr. S\"{o}ren Wolfers and Prof. Jesper Oppelstrup. They both provided us with valuable ideas and constructive comments.

\section*{Statements and Declarations}
\subsection*{Conflict of interest}
The authors have no conflicts to disclose.

\subsection*{Data Availability}

The data that support the findings of this study are available from the corresponding author upon reasonable request. The Python implementation of \texttt{MICE} and the baseline methods used to generate the data and figures in this work are available at GitHub \url{https://github.com/agcarlon/mice}. Moreover, \texttt{MICE} can be installed using PyPI \url{https://pypi.org/project/mice/}. We use the datasets \emph{mushrooms}, \emph{gisette}, and \emph{HIGGS}, obtained from LibSVM \url{https://www.csie.ntu.edu.tw/~cjlin/libsvmtools/datasets/binary.html}.
The exact configurations and random seeds used to generate all figures are provided in the repository.

\section{Conclusion}

We introduced the \emph{Multi-Iteration Stochastic Optimizers}, a novel class of first-order stochastic optimization methods that use the \emph{Multi-Iteration stochastiC Estimator} (\texttt{MICE}). The \texttt{MICE} estimator utilizes successive control variates along optimization paths, efficiently reusing previously computed gradients to achieve accurate mean gradient approximations. At each iteration, it adaptively samples gradients to satisfy a user-specified relative error tolerance. Moreover, it employs a greedy strategy to determine which iterates to retain in memory and which to discard.

Thanks to its ability to control relative gradient error, \texttt{MICE} facilitates robust stopping criteria based on the gradient norm. Moreover, its nonintrusive design makes it readily integrable with existing first-order optimization methods, significantly expanding its applicability.

We provided a rigorous theoretical analysis for \texttt{SGD-MICE}, highlighting its strong performance across different classes of optimization problems. In particular, we proved exponential convergence in the $L^2$ sense for gradient-dominated functions under constant step-size conditions. For strongly convex problems, our results demonstrate that \texttt{SGD-MICE} achieves accuracy $tol$ with an average complexity of $\cl{O}(tol^{-1})$ gradient evaluations, outperforming conventional adaptive batch-size SGD methods, which require $\cl{O}(tol^{-1}\log(tol^{-1}))$ evaluations.

Numerically, we validated our theory through three illustrative examples, employing consistent \texttt{MICE} parameters across diverse scenarios. The tests ranged from a quadratic function with stochastic Hessian to a stochastic Rosenbrock problem solved via \texttt{Adam-MICE}, and finally logistic regression training on large-scale datasets. In the latter, \texttt{SGD-MICE} demonstrated competitive performance against established variance-reduction methods (\texttt{SAG}, \texttt{SAGA}, \texttt{SVRG}, and \texttt{SARAH}), confirming its practical efficiency and scalability.

Future research directions include extending \texttt{MICE} to quasi-Newton methods, exploring constrained optimization settings through standard techniques like projected gradients and active set methods, and investigating additional gradient estimation error sources, such as biases from numerical discretizations.
Exploring \texttt{MICE} in more realistic optimization scenarios, such as training neural networks \cite{xu2020second}, is another natural direction.

{
\footnotesize

\bibliography{references.bib}
}
\appendix

\section{Multi-iteration stochastic optimizers} \label{sec:algorithms}
In this section, we present the detailed algorithms for the multi-iteration stochastic optimizers using the \texttt{MICE} estimator for the mean gradient.
In Algorithms \ref{alg:sgd}, \ref{alg:adam}, we respectively describe the pseudocodes for \texttt{SGD-MICE} and \texttt{Adam-MICE}.

\begin{algorithm}[h]
\setstretch{1.5}
\caption{Pseudocode for \texttt{SGD-MICE} with fixed step-size.
\texttt{SGD-MICE} requires an unbiased estimator of the true gradient, $\nabla f$; 
a distribution from which $\bs{\theta}$ can be sampled, $\pi$;
a starting point, $\bs{\xi}_0$; 
and a tolerance on the squared gradient norm, $tol$.
}
\label{alg:sgd}
\begin{algorithmic}[1]
\Procedure{SGD-MICE}{$\nabla_{\bs{\xi}} f$, $\pi$, $\bs{\xi}_0$, $tol$}
\State $k \gets 0$
\While{$\norm{\grad \cl{F}_{k}^{\text{stop}}}^2 > tol$}
\Comment $\norm{\grad \cl{F}_{k}^{\text{stop}}}$ computed as in Remark~\ref{rmk:grad_resampling}
\State Evaluate $\nabla_{\bs{\xi}} \cl{F}_k$ using Algorithm~\ref{alg:mice}
\State $\bs{\xi}_{k+1} \gets \bs{\xi}_k - \eta \nabla_{\bs{\xi}} \cl{F}_k$
\State $k \gets k+1$
\EndWhile
\State \textbf{return} optimum approximation $\bs{\xi}_{k^*}$
\EndProcedure
\end{algorithmic}
\end{algorithm}

\begin{algorithm}[h]
\setstretch{1.5}
\caption{Pseudocode for \texttt{Adam-MICE} with fixed step-size.
\texttt{Adam-MICE} requires an unbiased estimator of the true gradient, $\nabla f$; 
a distribution from which $\bs{\theta}$ can be sampled, $\pi$;
a starting point, $\bs{\xi}_0$; 
and a tolerance on the squared gradient norm, $tol$.
Moreover, \texttt{Adam-MICE} requires the constants $\beta_1$, $\beta_2$, and $\epsilon_{\text{Adam}}$.
We use the values recommended by \cite{kingma2014adam}, $\beta_1=0.9$, $\beta_2=0.999$, and $\epsilon_{\text{Adam}}=10^{-8}$.
}
\label{alg:adam}
\begin{algorithmic}[1]
\Procedure{Adam-MICE}{$\nabla_{\bs{\xi}} f$, $\pi$, $\bs{\xi}_0$, $tol$} 
\State Initialize $\bs{m}_0$ and $\bs{v}_0$ as zero-vectors
\State $k \gets 0$
\While{$\norm{\grad \cl{F}_{k}^{\text{stop}}}^2 > tol$}
\Comment $\norm{\grad \cl{F}_{k}^{\text{stop}}}$ computed as in Remark~\ref{rmk:grad_resampling}
\State Evaluate $\nabla_{\bs{\xi}} \cl{F}_k$ using Algorithm~\ref{alg:mice}
\State $\bs{m}_{k+1} \gets \beta_1 \bs{m}_{k} + (1 - \beta_1) \nabla_{\bs{\xi}} \cl{F}_k$
\State $\bs{v}_{k+1} \gets \beta_2 \bs{v}_{k} + (1 - \beta_2) \nabla_{\bs{\xi}} \cl{F}_k^2$ \Comment{The gradient estimates are squared element-wise}
\State $\bs{\hat{m}}_{k+1} \gets \bs{m}_{k+1} / (1- \beta_1^{k+1})$
\State $\bs{\hat{v}}_{k+1} \gets \bs{v}_{k+1} / (1- \beta_2^{k+1})$
\State $\bs{\xi}_{k+1} \gets \bs{\xi}_k - \eta \bs{\hat{m}}_{k+1} / (\sqrt{\bs{\hat{v}}_{k+1}} + \epsilon_{\text{Adam}})$
\State $k \gets k+1$
\EndWhile
\State \textbf{return} optimum approximation $\bs{\xi}_{k^*}$
\EndProcedure
\end{algorithmic}
\end{algorithm}

Adapting stochastic optimization algorithms to use \MICE is as straight-forward as substituting the gradient estimator and the stopping criterion, as can be seen in Algorithms~\ref{alg:sgd} and \ref{alg:adam}.

\section{Variance reduction gradient estimators in the finite-sum case}
\label{sec:var_red_methods}
This short appendix subsection recalls the standard gradient estimators used by \texttt{SVRG}, \texttt{SARAH}, \texttt{SAG}, and \texttt{SAGA} in the finite-sum setting and explains how \texttt{SVRG} and \texttt{SARAH} correspond to representative choices of the \texttt{MICE} index set.

\subsection*{Finite-sum objective}
Consider
\begin{equation}\label{eq:finite_sum_appendix}
  \min_{\bs{\xi}\in\rset^{d_{\bs{\xi}}}}
  F(\bs{\xi})
  \coloneqq
  \frac{1}{n}\sum_{i=1}^n f_i(\bs{\xi}),
\end{equation}
with full gradient $\nabla F(\bs{\xi}) = \frac{1}{n}\sum_{i=1}^n \nabla f_i(\bs{\xi})$.
At iteration $k$, let $\cl{B}_k\subset\{1,\dots,n\}$ denote a mini-batch of indices.

\subsection*{SVRG}
\texttt{SVRG} \cite{Joh13} uses a control variate anchored at a snapshot point $\bs{\xi}_{\tilde{k}}$ (updated periodically every $m$ iterations in the classical implementation).
Given a snapshot iterate $\bs{\xi}_{\tilde{k}}$ and its full gradient $\nabla F(\bs{\xi}_{\tilde{k}})$, the \texttt{SVRG} estimator at $\bs{\xi}_k$ is
\begin{equation}\label{eq:svrg_estimator_appendix}
  v_k
  \coloneqq
  \frac{1}{|\cl{B}_k|}
  \sum_{i\in\cl{B}_k}\bigl(\nabla f_i(\bs{\xi}_k) - \nabla f_i(\bs{\xi}_{\tilde{k}})\bigr)
  + \nabla F(\bs{\xi}_{\tilde{k}}).
\end{equation}
In \texttt{MICE} terms, this corresponds to the two-level index set $\cl{L}_k=\{\tilde{k},k\}$: the snapshot provides an anchor (a ``level'' computed with a full batch), and the current iterate contributes a correction via gradient differences.
Updating the snapshot every $m$ iterations corresponds to a periodic restart of the anchor.

\subsection*{SARAH}
\texttt{SARAH} \cite{nguyen2017sarah} uses successive control variates between consecutive iterates.
With an initial full-gradient computation $v_{\tilde{k}}\coloneqq \nabla F(\bs{\xi}_{\tilde{k}})$, it updates recursively
\begin{equation}\label{eq:sarah_estimator_appendix}
  v_k
  \coloneqq
  \frac{1}{|\cl{B}_k|}
  \sum_{i\in\cl{B}_k}\bigl(\nabla f_i(\bs{\xi}_k) - \nabla f_i(\bs{\xi}_{k-1})\bigr)
  + v_{k-1}.
\end{equation}
This corresponds to the add-only / consecutive-differences regime of \texttt{MICE}: taking $\cl{L}_k=\{0,1,\dots,k\}$ yields a telescoping sum of level differences (see Remark~\ref{rmk:add_only} in the main text), while \texttt{MICE} further allows (i) cumulative sample growth at each level and (ii) adaptive sparsification of $\cl{L}_k$ via Drop/Restart/Clip operators.

\subsection*{SAG}
\texttt{SAG} \cite{schmidt2017minimizing} maintains a memory table of per-sample gradients.
Let $\{g_i\}_{i=1}^n$ store the most recently computed gradients and define their running average
\begin{equation}\label{eq:sag_average_appendix}
  \bar{g}_k \coloneqq \frac{1}{n}\sum_{i=1}^n g_i .
\end{equation}
Sampling $i_k\sim\mathrm{Unif}\{1,\dots,n\}$, \texttt{SAG} updates $g_{i_k}\leftarrow \nabla f_{i_k}(\bs{\xi}_k)$ and uses $\bar{g}_k$ as a direction.
Unlike \texttt{SVRG} and \texttt{SARAH}, this estimator is generally biased at finite $k$ because the table entries correspond to gradients evaluated at different past iterates; nevertheless, the method enjoys linear convergence under standard assumptions.
\texttt{SAG} is not naturally represented by an iteration-indexed set $\cl{L}_k\subset\{0,\dots,k\}$, since its control variates are indexed by \emph{data points}.

\subsection*{SAGA}
\texttt{SAGA} \cite{defazio2014saga} stores a memory table $\{g_i\}_{i=1}^n$ of the most recently computed per-sample gradients.
Sampling $i_k\sim\mathrm{Unif}\{1,\dots,n\}$, its estimator is
\begin{equation}\label{eq:saga_estimator_appendix}
  v_k
  \coloneqq
  \nabla f_{i_k}(\bs{\xi}_k) - g_{i_k} + \frac{1}{n}\sum_{i=1}^n g_i,
  \qquad
  g_{i_k}\leftarrow \nabla f_{i_k}(\bs{\xi}_k).
\end{equation}
As with \texttt{SAG}, \texttt{SAGA} is not directly captured by an iteration-indexed set $\cl{L}_k\subset\{0,\dots,k\}$ because its control variates are indexed by \emph{data points} rather than past iterates.
Nevertheless, both \texttt{SAGA} and \texttt{MICE} reduce variance by reusing previously computed gradient information: \texttt{SAGA} stores per-sample gradients, while \texttt{MICE} stores (and resamples) per-iteration gradient differences along the optimization path.

\section{Error decomposition of the \texttt{MICE} estimator}
\label{sec:error_decomposition}

The \MICE estimator has a conditional bias due to the reuse of previous information.
Here we prove that, if the statistical error of the estimator is controlled every iteration, then the bias is implicitly controlled as well.
Recall the \texttt{MICE} estimator is defined as
\begin{equation}
  \Fk = \sum_{\ell \in \cl{L}_k} \frac{1}{M_{\ell, k}} \sum_{\alpha \in \cl{I}_{\ell, k}} \Delta_{\ell, k, \alpha},
\end{equation}
where
\begin{equation}
  \Delta_{\ell, k, \alpha} =
  \begin{cases}
    \nabla f(\xil, \bs{\theta}_\alpha) - \nabla f(\xilpr, \bs{\theta}_\alpha) & \text{ if } \ell > \min\{\cl{L}_k\} \\
    \nabla f(\bs{\xi}_0, \bs{\theta}_\alpha) & \text{ if } \ell = \min\{\cl{L}_k\}.
  \end{cases}
\end{equation}

The squared $L^2$ error of the \MICE estimator can be decomposed as
\begin{align}
  \label{eq:l2_err_decomposed}
  \Exp{\norm{\Fk - \grad F(\xik)}^2}
  &= 
  \underbrace{
  \Exp{\norm{\Fk - \Expc{\Fk}{\xiset{k}}}^2}
  }_{\text{statistical error}}
  \\
  &\quad +
  \underbrace{
  \Exp{\norm{\Expc{\Fk}{\xiset{k}} - \grad F(\xik)}^2}
  }_{\text{bias term}}.
\end{align}
due to
\begin{multline}
  \Exp{
    \Expc{
      \Big\langle
        \Fk - \Expc{\Fk}{\xiset{k}}
      ,
        \Expc{\Fk}{\xiset{k}} - \grad F(\xik)
      \Big\rangle
    }{\xiset{k}}
  } \\
  =
  \Exp{
      \Big\langle
        \underbrace{
        \Expc{\Fk - \Expc{\Fk}{\xiset{k}}}{\xiset{k}}
        }_{=0}
      ,
      \Expc{\Fk}{\xiset{k}} - \grad F(\xik)
      \Big\rangle
  }
\end{multline}

Before we analyze the bias and statistical errors, let us analyze the conditional expectation of the \MICE estimator,
\begin{align}
  \Expc{\Fk}{\xiset{k}}
  &=
  \sum_{\ell \in \cl{L}_k} \frac{1}{M_{\ell, k}} \sum_{\alpha \in \cl{I}_{\ell, k}}
  \Expc{\Delta_{\ell, k, \alpha}}{\xiset{k}},
\end{align}
and noting that, for $\tilde{k} < k$, $\Delta_{\ell, \tilde{k}, \alpha} | \xiset{k}$  is deterministic,
\begin{equation}
  \Expc{\Delta_{\ell, \tilde{k}, \alpha}}{\xiset{k}} =
  \begin{cases}
    \Delta_{\ell, \tilde{k}, \alpha}  & \text{ if } \tilde{k} < k \\
    \nabla F(\xil) - \nabla F(\xilpr)               & \text{ if } \tilde{k} = k.
  \end{cases}
\end{equation}

Let $\cl{L}^{\cap}_k = \cl{L}_k \cap \cl{L}_{k-1}$.
Splitting the summands in \texttt{MICE} between the terms computed at $k$ and the previous ones,
\begin{equation}
  \label{eq:split_fk}
  \Fk =
  \underbrace{
    \left\{
      \sum_{\ell \in \cl{L}^{\cap}_k} \frac{1}{M_{\ell, k}}
      \sum_{\alpha \in \cl{I}_{\ell, k-1}}
      \Delta_{\ell, k, \alpha}
    \right\}
  }_{\text{previously computed}}
  +
  \underbrace{
    \left\{
      \frac{1}{M_{k, k}}
      \sum_{\alpha \in \cl{I}_{k, k}}
      \Delta_{k, k, \alpha}
      +
      \sum_{\ell \in \cl{L}^{\cap}_k}
      \frac{1}{M_{\ell, k}}
      \sum_{\alpha \in \cl{I}_{\ell, k} \setminus \cl{I}_{\ell, k-1}}
      \Delta_{\ell, k, \alpha}
    \right\}
  }_{\text{computed at $k$}},
\end{equation}
taking the expectation conditioned on $\xiset{k}$, and using $\nabla F(\bs{\xi}_{p_k(\min\{\cl{L}_k\})})=0$,
\begin{align}
  \Expc{\Fk}{\xiset{k}} 
  &=
  \!\begin{multlined}[t][12cm]
    \sum_{\ell \in \cl{L}^{\cap}_k} \frac{1}{M_{\ell, k}} \sum_{\alpha \in \cl{I}_{\ell, k-1}} \Delta_{\ell, k, \alpha}
    \\
  +
    \nabla F(\bs{\xi}_k)
    - \nabla F(\bs{\xi}_{k-1})
    + \sum_{\ell \in \cl{L}^{\cap}_k} \frac{M_{\ell, k} - M_{\ell, k-1}}{M_{\ell, k}}
    (\nabla F(\xil) - \nabla F(\xilpr))
  \end{multlined}\\
  &=
  \label{eq:split_fk_expectation_1}
  \!\begin{multlined}[t][10.5cm]
  \sum_{\ell \in \cl{L}^{\cap}_k} \frac{1}{M_{\ell, k}} \sum_{\alpha \in \cl{I}_{\ell, k-1}} \Delta_{\ell, k, \alpha}
  + \nabla F(\xik) -
  \sum_{\ell \in \cl{L}^{\cap}_k} \frac{M_{\ell, k-1}}{M_{\ell, k}}
  (\nabla F(\xil) - \nabla F(\xilpr))
  \end{multlined}
  \\
  &=
  \label{eq:split_fk_expectation}
  \grad F(\xik) +
    \sum_{\ell \in \cl{L}^{\cap}_k} \frac{M_{\ell, k-1}}{M_{\ell, k}} \muhat{\ell, k-1}
    - \sum_{\ell \in \cl{L}^{\cap}_k} \frac{M_{\ell, k-1}}{M_{\ell, k}} \bs{\mu}_{\ell, k-1},
\end{align}
where in \eqref{eq:split_fk_expectation_1} we used $\sum_{\ell \in \cl{L}^{\cap}_k} (\nabla F(\xil) - \nabla F(\xilpr)) = \grad F(\bs{\xi}_{k-1})$.

Next, we investigate the bias of the \texttt{MICE} estimator conditioned on the current iterate $\xik$ and its contribution to the squared $L^2$ error.
\begin{prop}[Bias of the \MICE estimator in expectation minimization]
  \label{prp:bias}
  Let the bias of the \MICE estimator be defined as
  \begin{equation}
    \bias \coloneqq \Expc{\Fk}{\xiset{k}} - \nabla F(\xik).
  \end{equation}
  Then, the bias is
  \begin{equation} \label{eq:bias_new}
    \bias
    =
    \sum_{\ell \in \cl{L}^{\cap}_k}
    \frac{M_{\ell, k-1}}{M_{\ell, k}}
    \left(
      \muhat{\ell, k-1}
      - \bs{\mu}_{\ell, k-1}
    \right),
  \end{equation}
  and its contribution to the squared $L^2$ error is
  \begin{align}
    \Exp{\norm{\bias}^2}
    &
    \label{eq:bias_true}
    =
    \Exp{
    \sum_{\ell \in \cl{L}^{\cap}_k}
    \frac{M_{\ell, k-1}}{M_{\ell, k}^2}
      V_{\ell, k-1}}.
  \end{align}
\end{prop}

\begin{proof}

Equation \eqref{eq:bias_new} follows directly in \eqref{eq:split_fk_expectation}.
Now, let's investigate $\Exp{\norm{\bias}^2}$.
\begin{align}
  \Exp{\norm{\bias}^2}
  &=
  \Exp{
      \norm{
        \sum_{\ell \in \cl{L}^{\cap}_k}
        \frac{M_{\ell, k-1}}{M_{\ell, k}}
        \left(
        \muhat{\ell, k-1}
        - \bs{\mu}_{\ell, k-1}
        \right)
    }^2
  }\\
  =&
  \!\begin{multlined}[t][12cm]
   \bb{E}\Bigg[
    \sum_{\ell \in \cl{L}^{\cap}_k}
    \bigg(
    \frac{M_{\ell, k-1}^2}{M_{\ell, k}^2}
    \Expc{
      \norm{
        \muhat{\ell, k-1}
        - \bs{\mu}_{\ell, k-1}
      }^2}{\xiset{\ell}}
    \\
    +
    2 \sum_{j \in \cl{L}^{\cap}_k : j > \ell}
    \frac{M_{\ell, k-1}}{M_{\ell, k}}
    \frac{M_{j, k-1}}{M_{j, k}}
    \Expc{
      \inner{
        \muhat{\ell, k-1}
        - \bs{\mu}_{\ell, k-1}
      }{
        \muhat{j, k-1}
        - \bs{\mu}_{j, k-1}
      }
    }{\xiset{j}}
    \bigg)\Bigg].
  \end{multlined}
\end{align}

Using Lemma~\ref{lem:inner_product} and $\Expc{\norm{\muhat{\ell, k-1} - \bs{\mu}_{\ell, k-1}}^2}{\xiset{\ell}} = \V[k-1]{\ell}\M[k-1]{\ell}^{-1}$ concludes the proof.
\end{proof}

Note from \eqref{eq:bias_new} that $\Exp{\bias} = \bs{0}$.

\begin{cor}[Bias of the \MICE estimator in finite sum minimization]
  \label{cor:bias_finite}
  The bias $\bias$ of the \MICE estimator in finite sum minimization is similar to the expectation minimization one, with the consideration of the finite population correction factor,
  \begin{align}
    \Exp{\norm{\bias}^2}
    &
    =
    \Exp{
    \sum_{\ell \in \cl{L}^{\cap}_k}
    \left(
    \frac{N-\M[k-1]{\ell}}{N}
    \right)
    \frac{M_{\ell, k-1}}{M_{\ell, k}^2}
      V_{\ell, k-1}}.
  \end{align}
\end{cor}
\begin{proof}
  The proof follows exactly as in Proposition~\ref{prp:bias}, except the finite population correction factor is used in the centered second moment of $\muhat{\ell, k-1}$, $\Expc{\norm{\muhat{\ell, k-1} - \bs{\mu}_{\ell, k-1}}^2}{\xiset{\ell}} = 
  (N - \M[k-1]{\ell})N^{-1}
  \V[k-1]{\ell}\M[k-1]{\ell}^{-1}$.
\end{proof}

\begin{prop}[Statistical error of the \MICE estimator in expectation minimization]
  \label{prp:statistical_error}
  The statistical error of the \MICE estimator in the case of expectation minimization is
  \begin{align}
    \Exp{\norm{\Fk - \Expc{\Fk}{\xiset{k}}}^2}
    &=
    \Exp{
    \sum_{\ell \in \cl{L}^{\cap}_k}
    \frac{
      (M_{\ell, k} - M_{\ell, k-1}) V_{\ell, k}
    }{
      M_{\ell, k}^2
    }}
    + 
    \Exp{\frac{V_{k, k}}{M_{k, k}}}
  \end{align}
\end{prop}

\begin{proof}
  From \eqref{eq:l2_err_decomposed}, we can use Lemma~\ref{lem:mice_expected_error} and Proposition~\ref{prp:bias} to get
  \begin{align}
    \Exp{\norm{\Fk - \Expc{\Fk}{\xiset{k}}}^2}
    &= \Exp{\norm{\Fk - \grad F(\xik)}^2} - 
    \Exp{\norm{\Expc{\Fk}{\xiset{k}} - \grad F(\xik)}^2} \\
    &= \Exp{\sum_{\ell \in \cl{L}_k} \frac{\V{\ell}}{\M{\ell}}}
    - 
    \Exp{
    \sum_{\ell \in \cl{L}^{\cap}_k}
    \frac{M_{\ell, k-1}}{M_{\ell, k}^2}
      V_{\ell, k-1}}\\
    &= 
    \Exp{\sum_{\ell \in \cl{L}^{\cap}_k} \frac{\V{\ell}}{\M{\ell}}
    - 
    \sum_{\ell \in \cl{L}^{\cap}_k}
    \frac{M_{\ell, k-1}}{M_{\ell, k}^2}
      V_{\ell, k-1}}
    +
    \Exp{\frac{\V{k}}{\M{k}}}.
  \end{align}
  Using $\V[k]{\ell}=\V[k-1]{\ell}$ for $\ell \in \cl{L}^{\cap}_k$ concludes the proof.
\end{proof}

\begin{cor}[Statistical error of the \MICE estimator in finite sum minimization]
  \label{cor:statistical_error_finite}
  The statistical error in the finite sum minimization case is
    \begin{multline}
      \Exp{\norm{\Fk - \Expc{\Fk}{\xiset{k}}}^2}
      = \\
      \Exp{
      \sum_{\ell \in \cl{L}^{\cap}_k}
      \left(
        \frac{N-\M{\ell}}{N}
      \right)
      \frac{
        (M_{\ell, k} - M_{\ell, k-1}) V_{\ell, k}
      }{
        M_{\ell, k}^2
      }}
      + 
      \left(
        \frac{N-\M{k}}{N}
      \right)
      \Exp{\frac{V_{k, k}}{M_{k, k}}}
  \end{multline}
\end{cor}
\begin{proof}
  The proof follows exactly as in Proposition~\ref{prp:statistical_error}, except Remark~\ref{rmk:err_finite_sum} and Corollary~\ref{cor:bias_finite} are used instead of Lemma~\ref{lem:mice_expected_error} and Proposition~\ref{prp:bias}.
\end{proof}

\section{A high-probability bound for add-only \texttt{MICE}}
\label{sec:hp_mice}
In this appendix we give a simple high-probability control of the realized \texttt{MICE} gradient estimation error in an add-only regime, which is sufficient to prove Corollary~\ref{cor:hp_pl_mice}.

\subsection*{Setup (add-only / fixed-anchor regime)}
Recall $\Fk=\nabla_{\bs{\xi}}\cl{F}_k$ and define the realized gradient estimation error
\begin{equation}
  \label{eq:hp_err_def}
  e_k \coloneqq \Fk - \grad F(\xik).
\end{equation}
We assume add-only index sets $\cl{L}_{k-1}\subseteq \cl{L}_k$ and thus the previous rule $p_k(\ell)=\ell-1$, so that once a level $\ell$ is introduced it remains anchored to the fixed pair $(\bs{\xi}_\ell,\bs{\xi}_{\ell-1})$ until a restart occurs.
In this regime, each level contribution is a sample mean of i.i.d.\ centered increments at that anchored pair.
Moreover, conditional on the iterate history, the sample sets used at distinct levels are independent by construction, so the level contributions are conditionally independent.
We work under the tail assumption in Assumption~\ref{as:hp_coord_subg}.

\begin{lem}[A single good event for \texttt{MICE} error control]
\label{lem:hp_mice_single_event}
Assume Assumption~\ref{as:hp_coord_subg} holds for every level.
Assume that, conditional on the iterate history, the sample sets used at distinct levels are independent.
Fix a summable schedule $(\delta_k)_{k\ge 0}$ with $\sum_{k\ge 0}\delta_k\le \delta$.
Assume that at each iteration $k$, the sample sizes $\{M_{\ell,k}\}_{\ell\in\cl{L}_k}$ are chosen based only on the past (in particular, before drawing any new samples at iteration $k$).
Define the event $\Omega_\delta^{\mathrm{MICE}}$ by requiring that, for all $k\ge 0$,
\begin{equation}
  \label{eq:hp_mice_error_bound}
  \norm{e_k}
  \le
  \sqrt{2d_\xi\,\log\!\Big(\frac{2d_\xi}{\delta_k}\Big)}
  \left(\sum_{\ell\in\cl{L}_k}\frac{\sigma_\ell^2}{M_{\ell,k}}\right)^{1/2}.
\end{equation}
Then $\bb{P}[\Omega_\delta^{\mathrm{MICE}}]\ge 1-\delta$.
\end{lem}
\begin{proof}
Fix $k\ge 0$ and a coordinate $j\in\{1,\ldots,d_\xi\}$.
Condition on the iterate history and on the chosen sample sizes $\{M_{\ell,k}\}_{\ell\in\cl{L}_k}$, so that the sample sizes can be treated as fixed when applying tail bounds to the samples.
In the add-only regime,
\begin{equation}
  (e_k)_j
  =
  \sum_{\ell\in\cl{L}_k}\frac{1}{M_{\ell,k}}\sum_{i=1}^{M_{\ell,k}}(Z_{\ell,i})_j.
\end{equation}
By Assumption~\ref{as:hp_coord_subg}, each coordinatewise sample mean is sub-Gaussian with proxy variance $\sigma_\ell^2/M_{\ell,k}$, and by independence across levels the sum $(e_k)_j$ is sub-Gaussian with proxy variance
\(
  v_k \coloneqq \sum_{\ell\in\cl{L}_k}\sigma_\ell^2/M_{\ell,k}.
\)
Therefore, for all $t>0$,
\begin{equation}
  \bb{P}\big[|(e_k)_j|\ge t\big]\le 2\exp\!\left(-\frac{t^2}{2v_k}\right).
\end{equation}
Taking $t=\sqrt{2v_k\log(2d_\xi/\delta_k)}$ gives $\bb{P}[|(e_k)_j|\ge t]\le \delta_k/d_\xi$.
A union bound over $j=1,\ldots,d_\xi$ yields that with probability at least $1-\delta_k$,
\(
  \norm{e_k}_\infty\le \sqrt{2v_k\log(2d_\xi/\delta_k)}.
\)
Since $\norm{e_k}\le \sqrt{d_\xi}\,\norm{e_k}_\infty$, the bound \eqref{eq:hp_mice_error_bound} holds with probability at least $1-\delta_k$.
Finally, a union bound over $k\ge 0$ and $\sum_{k\ge 0}\delta_k\le \delta$ yields $\bb{P}[\Omega_\delta^{\mathrm{MICE}}]\ge 1-\delta$.
\end{proof}

\begin{cor}[Uniform relative error control under the variance-sum constraint]
\label{cor:hp_mice_rel_err}
Under the conditions of Lemma~\ref{lem:hp_mice_single_event}, if the sample sizes satisfy \eqref{eq:hp_variance_sum_constraint} for all $k\ge 0$, then with probability at least $1-\delta$,
\begin{equation}
  \norm{e_k}\le \epsilon\norm{\grad F(\xik)}, \qquad \forall\,k\ge 0.
\end{equation}
Consequently, Theorem~\ref{thm:hp_pl_event} applies and \eqref{eq:hp_linear_rate} holds for all $k\ge 0$ with probability at least $1-\delta$.
\end{cor}

\begin{lem}[High-probability iteration bound to reach $tol$]
\label{lem:hp_iter_to_tol}
Assume the conditions of Corollary~\ref{cor:hp_pl_mice}.
In addition, assume $F$ satisfies Assumption~\ref{as:convex}.
Let $r_{\mathrm{hp}}\in(0,1)$ be the contraction factor in \eqref{eq:hp_linear_rate} and define
\begin{equation}
  \label{eq:K_tol_hp}
  K(tol)
  \coloneqq
  \left\lceil
  \frac{\log\!\big(2L(F(\xio)-F(\opt))/tol\big)}{\log(1/r_{\mathrm{hp}})}
  \right\rceil.
\end{equation}
Then, with probability at least $1-\delta$, we have $k^*(tol)\le K(tol)$, where $k^*(tol)$ is defined in \eqref{eq:stopping_criterion}.
\end{lem}
\begin{proof}
By Corollary~\ref{cor:hp_pl_mice}, with probability at least $1-\delta$ we have for all $k\ge 0$,
\begin{equation}
  F(\xik)-F(\opt)\le r_{\mathrm{hp}}^k\big(F(\xio)-F(\opt)\big).
\end{equation}
Since $F$ is convex and $L$-smooth, \eqref{eq:grad_bound} holds (as used earlier in Corollary~\ref{cor:grad_norm_convergence}).
Therefore, on the same event,
\begin{equation}
  \|\nabla F(\xik)\|^2
  \le
  2L\,r_{\mathrm{hp}}^k\big(F(\xio)-F(\opt)\big).
\end{equation}
By the definition of $K(tol)$, for all $k\ge K(tol)$ we have $2L\,r_{\mathrm{hp}}^k(F(\xio)-F(\opt))\le tol$, hence $\|\nabla F(\xik)\|^2\le tol$ and therefore $k^*(tol)\le K(tol)$.
\end{proof}

\begin{cor}[High-probability sampling-cost bound under variance-sum sizing]
\label{cor:hp_cost_to_tol}
Assume the conditions of Corollary~\ref{cor:hp_pl_mice} and Lemma~\ref{lem:hp_iter_to_tol}.
In particular, work in the add-only regime.
Fix a summable schedule $(\delta_k)_{k\ge 0}$ with $\sum_{k\ge 0}\delta_k\le \delta$ and assume it is nonincreasing.
For each iteration $k$, define
\begin{equation}
  \label{eq:Rk_hp}
  R_k
  \coloneqq
  \frac{\epsilon^2}{2d_\xi\,\log(2d_\xi/\delta_k)}\,\norm{\nabla F(\xik)}^2.
\end{equation}
Assume the sample sizes are chosen (predictably) to satisfy the variance-sum constraint
\begin{equation}
  \sum_{\ell\in\cl{L}_k}\frac{\sigma_\ell^2}{M_{\ell,k}}
  \le
  R_k,
  \qquad \forall\,k\ge 0,
\end{equation}
and that, at each iteration $k$, they are (approximately) work-optimal for the continuous relaxation
\begin{equation}
  \min_{\{M_{\ell,k}>0\}}\ \sum_{\ell\in\cl{L}_k} \cost{\ell}\,M_{\ell,k}
  \quad \text{s.t.}\quad
  \sum_{\ell\in\cl{L}_k}\frac{\sigma_\ell^2}{M_{\ell,k}}\le R_k,
\end{equation}
whose KKT optimizer is
\begin{equation}
  \label{eq:KKT_hp}
  M_{\ell,k}^\star
  =
  \frac{\sigma_\ell}{\sqrt{\cost{\ell}}}\,
  \frac{\sum_{j\in\cl{L}_k}\sigma_j\sqrt{\cost{j}}}{R_k},
  \qquad \ell\in\cl{L}_k.
\end{equation}
In the add-only regime with cumulative sampling, the cumulative gradient-evaluation cost $\cl{C}_{k}$ (defined as the telescoping sum of increments in \eqref{eq:cost_of_iter_k}) coincides with the cost of the final cumulative sample sizes at iteration $k$, i.e., $\cl{C}_{k}=\sum_{\ell\in\cl{L}_{k}} \cost{\ell}\,M_{\ell,k}$.
Then, with probability at least $1-\delta$,
\begin{equation}
  \label{eq:hp_cost_bound}
  \cl{C}_{k^*(tol)-1}
  \le
  \frac{2d_\xi\,\log(2d_\xi/\delta_{K(tol)-1})}{\epsilon^2\,tol}
  \left(\sum_{\ell\in\cl{L}_{K(tol)-1}}\sigma_\ell\sqrt{\cost{\ell}}\right)^2
  +
  \sum_{\ell\in\cl{L}_{K(tol)-1}}\cost{\ell},
\end{equation}
where $K(tol)$ is defined in \eqref{eq:K_tol_hp}.
\end{cor}
\begin{proof}
Under the KKT sizing \eqref{eq:KKT_hp}, the corresponding (continuous) sampling cost satisfies
\begin{equation}
  \sum_{\ell\in\cl{L}_k}\cost{\ell}\,M_{\ell,k}^\star
  =
  \frac{\left(\sum_{\ell\in\cl{L}_k}\sigma_\ell\sqrt{\cost{\ell}}\right)^2}{R_k}.
\end{equation}
Rounding to integers by $M_{\ell,k}=\lceil M_{\ell,k}^\star\rceil$ increases the left-hand side by at most $\sum_{\ell\in\cl{L}_k}\cost{\ell}$.
On the event of Lemma~\ref{lem:hp_iter_to_tol} we have $k^*(tol)\le K(tol)$.
At the last pre-stopping iterate $k=k^*(tol)-1$ we have $\|\nabla F(\xik)\|^2>tol$, hence from \eqref{eq:Rk_hp},
\begin{equation}
  \frac{1}{R_k}
  \le
  \frac{2d_\xi\,\log(2d_\xi/\delta_k)}{\epsilon^2\,tol}
  \le
  \frac{2d_\xi\,\log(2d_\xi/\delta_{K(tol)-1})}{\epsilon^2\,tol},
\end{equation}
where the last inequality uses that $(\delta_k)$ is nonincreasing and $k\le K(tol)-1$.
Finally, using add-only monotonicity $\cl{L}_{k^*(tol)-1}\subseteq \cl{L}_{K(tol)-1}$ yields \eqref{eq:hp_cost_bound}.
\end{proof}

\end{document}